\documentclass[reqno,10pt,a4paper]{amsart}
\usepackage[utf8]{inputenc}

\theoremstyle{definition}

\theoremstyle{remark}

\numberwithin{equation}{section}

\usepackage{amssymb,eucal,mathrsfs,array,setspace,geometry,enumitem,cite,tensor,amsmath,wrapfig,amscd,mathptmx,graphicx,bm,bbm,multirow,multicol,adjustbox}
\usepackage[normalem]{ulem}
\usepackage[dvipsnames]{xcolor}
\usepackage{tikz}
\usetikzlibrary{knots,calc}
\usetikzlibrary{decorations.markings}
\usepackage[centertableaux]{ytableau}
\usepackage{txfonts}            
\usepackage{newtxtext}            
\usepackage{mathtools} 
\usepackage{stmaryrd} 

\usepackage{caption}
\usepackage{subcaption}
\usepackage{xcolor}

\usetikzlibrary{positioning}
\usetikzlibrary{decorations.markings}

\usepackage[colorlinks=true,citecolor=red,linkcolor=blue]{hyperref} 
\usepackage[capitalise,noabbrev]{cleveref}

\usepackage{tikz-cd}
\usetikzlibrary{decorations.pathmorphing}

\newcommand{\correct}[1]{#1}

\usetikzlibrary{knots}

\DeclareSymbolFont{largesymbols}{OMX}{zplm}{m}{n} 

\setitemize{leftmargin=*}   
\setenumerate{leftmargin=*} 
\setlist[enumerate,1]{label=\textup{(\arabic*)}, 
  ref=\arabic*} 
\setlist[enumerate,2]{label = \textup{(\roman*)},
  ref = \theenumi.\roman*}

\let\originalleft\left     
\let\originalright\right
\renewcommand{\left}{\mathopen{}\mathclose\bgroup\originalleft}
\renewcommand{\right}{\aftergroup\egroup\originalright}

\newcolumntype{C}{>{$}c<{$}} 

\numberwithin{equation}{section}
\allowdisplaybreaks

\renewcommand{\ge}{\geq}
\renewcommand{\le}{\leq}

\DeclarePairedDelimiter{\brac}{\lparen}{\rparen} 
\DeclarePairedDelimiter{\sqbrac}{\lbrack}{\rbrack} 
\DeclarePairedDelimiter{\set}{\lbrace}{\rbrace}
\newcommand{\st}{\mspace{5mu} \vert \mspace{5mu}} 

\DeclarePairedDelimiter{\floor}{\lfloor}{\rfloor}
\DeclarePairedDelimiter{\ang}{\langle}{\rangle}

\DeclarePairedDelimiterX{\comm}[2]{\lbrack}{\rbrack}{#1 , #2}  
\DeclarePairedDelimiterX{\acomm}[2]{\lbrace}{\rbrace}{#1 , #2} 
\DeclarePairedDelimiterX{\super}[2]{\lparen}{\rparen}{#1 \delimsize\vert \mathopen{} #2} 

\newcommand{\pair}[2]{\ang*{#1,#2}} 

\DeclareMathOperator{\id}{id}
\DeclareMathOperator{\im}{im}
\DeclareMathOperator{\diag}{diag}

\newcommand{\pd}{\partial}     

\newcommand{\dd}{\mathrm{d}}   
\newcommand{\ii}{\mathfrak{i}} 
\newcommand{\ee}{\mathsf{e}}   
\newcommand{\epi}{\mathbf{e}}

\newcommand{\wun}{\mathbf{1}}  
\newcommand{\mdw}{p}
\newcommand{\mdwt}{\mdw_0} 
\newcommand{\mdwtt}{\mdw_1}
\newcommand{\mdwttt}{\mdw_2} 

\newcommand{\eis}[1]{\mathcal{G}_{#1}}

\DeclareMathOperator{\Res}{Res}
\DeclareMathOperator{\cspn}{span}
\newcommand{\spn}[1]{\cspn_{\CC}\set*{#1}}                    

\DeclareMathOperator{\Hom}{Hom}
\newcommand{\Homgrp}[3]{\Hom_{#1}\brac*{#2,#3}}

\newcommand{\fld}[1]{\mathbb{#1}}    
\newcommand{\alg}[1]{\mathfrak{#1}}  
\newcommand{\grp}[1]{\mathsf{#1}}    
\newcommand{\categ}[1]{\mathscr{#1}} 

\newcommand{\ZZ}{\fld{Z}}
\newcommand{\NN}{\fld{N}}

\newcommand{\RR}{\fld{R}}
\newcommand{\CC}{\fld{C}}
\newcommand{\HH}{\fld{H}}

\newcommand{\SLG}[2]{\grp{#1} \brac[\big]{#2}}       
\newcommand{\SLTZ}{\grp{SL\brac*{2,\ZZ}}}
\newcommand{\PSLTZ}{\grp{PSL\brac*{2,\ZZ}}}
\newcommand{\SLTZp}[1]{\grp{SL}\brac*{2,\ZZ_{#1}}}
\newcommand{\GL}[1]{\grp{GL}\brac*{#1,\CC}}
\newcommand{\Bthree}{\grp{B}_3}

\newcommand{\SLA}[2]{\alg{#1} \brac*{#2}}                 

\DeclareMathOperator{\wt}{wt}                                               

\newcommand{\catC}{\categ{C}}           

\newcommand{\rep}[1]{\mathsf{Rep}\brac*{#1}}
\DeclareMathOperator{\adj}{ad}
\newcommand{\addcat}[1]{\categ{#1}_{\adj}}

\newcommand{\coup}[4]{\lambda^{#1}_{(#2,#3)#4}}
\newcommand{\dcoup}[4]{\Upsilon^{#1}_{(#2,#3)#4}}

\DeclarePairedDelimiter{\ket}{\lvert}{\rangle}
\DeclarePairedDelimiterX{\braket}[2]{\langle}{\rangle}{#1 \delimsize\vert \mathopen{} #2}
\DeclarePairedDelimiterX{\bracket}[3]{\langle}{\rangle}{#1 \delimsize\vert \mathopen{} #2 \delimsize\vert \mathopen{} #3}

\newcommand{\modS}{\mathsf{S}}                        
\newcommand{\modT}{\mathsf{T}}                        
\newcommand{\Rmat}{\mathsf{R}}                        
\newcommand{\Fmat}{\mathsf{F}}
\newcommand{\Gmat}{\mathsf{G}}

\newcommand{\Smat}[5]{\modS^{\brac*{#5}}_{#1 #2,#3 #4}}
\newcommand{\Tmat}[5]{\modT^{\brac*{#5}}_{#1 #2,#3 #4}}
\newcommand{\twist}[1]{\theta_{#1}}

\newcommand{\cft}{conformal field theory}

\newcommand{\voa}{vertex operator algebra}

\newcommand{\PBW}{Poincar\'{e}-Birkhoff-Witt}

\newcommand{\lhs}{left side}
\newcommand{\rhs}{right side}

\newcommand{\cfin}[1]{\(C_{#1}\)-cofinite}
\newcommand{\ctwo}{\cfin{2}}
\newcommand{\mtc}{modular tensor category}
\newcommand{\mtcs}{modular tensor categories}

\renewcommand{\epsilon}{\varepsilon}
\newcommand{\vvmf}{vector-valued modular form}

\newcommand{\refThm}{Theorem }
\newcommand{\refProp}{Proposition }
\newcommand{\refCor}{Corollary }
\newcommand{\refLem}{Lemma } 
\newcommand{\refEq}{Equation }
\newcommand{\refSec}{Section }
\newcommand{\refTab}{Table }
\newcommand{\refRmk}{Remark }

\theoremstyle{plain}
\newtheorem{thm}{Theorem}[section]
\newtheorem{prop}[thm]{Proposition}
\newtheorem{lem}[thm]{Lemma}

\newtheorem*{thm*}{Theorem}

\theoremstyle{definition} 

\newtheorem{defn}[thm]{Definition}

\Crefname{thm}{Theorem}{Theorems}
\Crefname{prop}{Proposition}{Propositions}
\Crefname{lem}{Lemma}{Lemmas}
\Crefname{cor}{Corollary}{Corollaries}
\Crefname{defn}{Definition}{Definitions}
\Crefname{tab}{Table}{Tables}

\newcommand{\intset}{\Xi}
\newcommand{\ptf}{\psi}    
\newcommand{\vvptf}{\Psi}    

\begin{document}

\title{The modular properties of \(\SLA{sl}{2}\) torus \(1\)-point functions}

\author{Matthew Krauel}
\address{
  Department of Mathematics and Statistics \\
  California State University, Sacramento \\
    6000 J.\ Street \\
  Sacramento, California, USA, 95819.
}
\curraddr{}
\email{krauel@csus.edu}
\thanks{}
\author{Jamal Noel Shafiq}
\address{
School of Mathematics \\
Cardiff University \\
Cardiff, United Kingdom, CF24 4AG.
}
\curraddr{}
\email{shafiqjn@cardiff.ac.uk}
\thanks{}
\author{Simon Wood}
\address{
School of Mathematics \\
Cardiff University \\
Cardiff, United Kingdom, CF24 4AG.
}
\curraddr{}
\email{woodsi@cardiff.ac.uk}
\thanks{}


\date{}

\dedicatory{}

\begin{abstract}
  Conformal field theory and its axiomatisation in terms of \voa{s} or
    chiral algebras are most commonly considered on the    
    Riemann sphere. However, an important constraint in physics and an
    interesting source of mathematics is the fact that conformal field
    theories are expected to be well defined on any Riemann surface. To this
    end, a thorough understanding of chiral torus \(1\)-point functions, 
    ideally including explicit formulae, is a prerequisite for a detailed
    understanding of higher genera. These are distinguished from characters or
    vacuum torus \(1\)-point functions because the insertion point is explicitly allowed to be
    labelled by any module over the \voa{} rather than just the \voa{} itself.
    
  Compellingly, chiral torus \(1\)-point functions exhibit interesting modular properties, which we explore here in
    the example of the simple affine $\SLA{sl}{2}$ \voa{s} at non-negative integral levels. 
    We determine the dimension of the space spanned by such functions, choose
    a natural basis to construct vector-valued modular forms and describe the
    congruence properties of these forms. In particular, we find explicit generators
    for the spaces of all vector-valued modular forms of dimension at most
    three, when the insertion comes from a simple module other than the \voa{}.
    Finally, we use the fact that categories of modules over rational \voa{s}
    are modular tensor categories to give explicit formulae
    for the action of the modular group on chiral torus \(1\)-point functions entirely in terms of
    categorical data. The usual modular \(\modS\) and \(\modT\) matrices of
    characters are known not to be complete invariants of modular tensor
    categories, so these generalised modular data are good candidates for more
    fine-grained invariants.
\end{abstract}

\maketitle


\section{Introduction}
  An important theme of conformal field theory is well definedness,
  not only
  on the Riemann sphere, but on all Riemann surfaces in the sense that
  \(n\)-point 
  correlation functions should be well defined on all these surfaces
  \cite{Cardy-Operator, MSMTC1089}.
The special case of vacuum torus \(1\)-point functions has been extensively
studied
in the literature. For example, their congruence had long been conjectured \cite{Moore-Symmetry,
  Eholzer-Classification, Eholzer-Uniqueness, DM-MoonshineSurvey,
  Bauer-Comments}
and was later proved for a large class of \voa{s} in
\cite{Dong-Congruence}. Specifically,
evaluating such vacuum torus \(1\)-point functions at the
vacuum vector gives characters (or graded dimensions) which are central to many important research strands, such as classifying rational (bulk) conformal field
theories built from a given \voa{} \cite{Cappelli87, Mathur88-ClassificationRCFT}, the origins of Moonshine \cite{ConwayNorton79, Gannon-MoonsineBeyond}, classifying certain families of \voa{s} \cite{Schellekens93, FrancMason20-3dimClass, Mason21-2dimClass}, and tackling number theoretic and combinatorial problems \cite{Milas04-Ramanujan,Cook07,Milas08}. However, general torus \(1\)-point functions have
received far less attention. To the best of the authors' knowledge, in the context of \voa{s} the only case considered in the literature is the family of Virasoro minimal models
  \cite{intvvmf18}. Thus, one goal of this paper is to detail another family
  of examples: the important class of simple affine \voa{s} \(L(k,0)\)
  constructed from \(\SLA{sl}{2}\) at non-negative integral levels \(k\), also known
  as $\SLG{SU}{2}$ Wess–Zumino–Witten models in the physics
  literature. We note that some relevant, but different work is performed in
  \cite{EtingofKirillov,Hara03} where generalisations of Jack symmetric polynomials are
  constructed in the context of affine Lie algebras. Meanwhile, although highest weight vectors govern the theory
  developed in \cite{intvvmf18}, this is not the case for the \voa{s} here,
  where the story is much more subtle. Therefore, a dual goal of this paper is
  to improve on the tools that are available for studying torus \(1\)-point
  functions, which includes the incorporation of categorical notions and resources surrounding modular tensor categories.
In particular, we show that
these \(1\)-point functions are a source of infinite families of
non-congruence representations of the modular group, we derive
explicit formulae for \(q\)-series expansions, and we contrast these analytic
number theoretic data with the categorical data of the modular tensor categories
formed by \(L(k,0)\)-modules.

We provide a little more context for readers unfamiliar with
\(n\)-point correlation functions in \cft{}. On higher genus Riemann surfaces
these are in practice constructed from those on the sphere by gluing together
points to add handles. For each pair of points glued in this way, the number of points in the
correlation function decreases by two and the genus of the surface increases
by one. The configuration of these points determines the complex structure of
the resulting surface with many different configurations giving equivalent
complex structures. All configurations giving an equivalent complex structure
are 
famously related by the actions of mapping class groups. Due to conformal
invariance being closely related to the existence of complex structure, one
may be tempted to conclude that a well defined conformal field theory should
not be able to distinguish Riemann surfaces with equivalent complex
structures. However, this is only true for bulk or full conformal field
theory. For chiral conformal field theory, which is the focus here (specifically
its algebraic axiomatisation in the form of \voa{s}), one
merely has that the mapping class groups act on the spaces of chiral correlation
functions, as opposed to this action being trivial. All considerations from here on will be purely chiral and so
henceforth correlation function or \(n\)-point function will refer to the
chiral version. In the special case of the torus with either \(0\) or \(1\) points,
the mapping class groups are, respectively,
\(\PSLTZ\) and \(\Bthree\) (the Braid group
on three strands). Recall that \(\Bthree\) is the universal central extension
of \(\PSLTZ\). It turns out that the action of \(\Bthree\) on torus \(1\)-point functions
can always be rescaled using multiplier systems to yield an action by its
quotient \(\SLTZ\). Hence the properties of torus correlation functions are
commonly presented in terms of \(\SLTZ\). The groups \(\PSLTZ\) and \(\SLTZ\)
are somewhat confusingly both commonly referred to as ``the modular group'' in
the literature and so one speaks of modular invariance of torus \(1\)-point functions.

  While \cite{MSMTC1089} gives a compelling motivation for the role of
  modular invariance, this has only been rigorously established for conformal
  field theories constructed from rational \ctwo{} \voa{s}.
  In this setting, torus \(n\)-point functions where the \(n\) points only
  take 
  insertions from the \voa{} can be constructed as traces of a product of
  \(n\) copies of the \voa{} action on some module \(M\). We call these
  \(n\)-point functions \emph{vacuum torus \(n\)-point functions} because the
  \voa{} is sometimes also called the \emph{vacuum module} (note that the
  insertions need not be the vacuum vector of the \voa{}). The special case of
  \(n=1\) with the insertion being the vacuum vector (this can also be
  thought of as a \(0\)-point function) is called the character of \(M\). In \cite{Zhu}, Zhu
  proved the modular invariance of such vacuum torus \(n\)-point functions and in
  particular showed that vacuum torus \(1\)-point functions are closed under the action of \(\SLTZ\).
  The properties of the \(\SLTZ\) representations arising from vacuum torus
  \(1\)-point functions have been heavily studied.
  Much of this work, for example the congruence property \cite{Dong-Congruence} mentioned above, rests on using the theory
of tensor categories \cite{Huang-Rigidity} and Verlinde's formula
\cite{Verlinde,HuaVer08}. An alternative route to studying congruence now exists in the case of
characters due to developments in number theory and a proof of the unbounded
denominator conjecture \cite{calegari2023unbounded}. However, again, this does not apply to general torus \(1\)-point functions such as those studied here.

  The transition from vacuum torus \(n\)-point functions to general torus
  \(n\)-point functions requires the replacement of \voa{} actions by
  intertwining operators. This case, as mentioned above, has
  so far received far less attention within the literature and is the focus here. For rational \ctwo{} 
  \voa{s} the modular invariance of general torus \(1\)-point functions was shown in
  \cite{miyamoto2000intertwining}. This was generalised to orbifolds in \cite{Yamauchi-IntertwinerModularity} and to
  torus \(n\)-point functions in \cite{Huang-IntModInv}. However, neither insights
  from Verlinde's formula nor the unbounded denominator property apply
  here. For example, it was observed in \cite{intvvmf18} that both congruence and non-congruence
  representations of \(\SLTZ\) appear for the Virasoro minimal models and we
  will show below that this is also the case for simple affine $\SLA{sl}{2}$ \voa{s} at non-negative integral levels.

  This paper is organised as follows. In \cref{sec:intmodprops} we review the
  analytic number theory required for studying vector-valued modular forms and
  how \voa{s} can be used to generate examples of vector-valued modular
  forms. The components (entries) of these vectors are torus
  \(1\)-point functions.

  In \cref{sec:1ptSpace} we review and develop general tools to characterise
  the space of all torus \(1\)-point functions (as modules over the algebra of
  holomorphic modular forms and the algebra \(\mathcal{R}\) of modular differential operators) obtained by varying the insertion
  vector over an entire simple \voa{} module.
  The main results are \cref{prop:RmoduleResult}, which gives sufficient
    conditions for the span of torus \(1\)-point functions obtained from
    Virasoro descendants of certain vectors to be a cyclic
    \(\mathcal{R}\)-module, and \cref{thm:generaldimformula} which
    gives 
    sufficient conditions for the span of all torus \(1\)-point functions
    to be a cyclic \(\mathcal{R}\)-module.

  In \cref{sec:sl2} we introduce the simple affine
  \voa{} constructed from \(\SLA{sl}{2}\) at non-integral levels.
  The main result of the section is the multi-part \cref{thm:sl2alldims},
  which collects the most important general results surrounding the analysis
  of torus $1$-point functions. These include finding vectors giving non-zero
  torus $1$-point functions, establishing linear independence among a certain set
  of these functions, obtaining that vectors generated from these functions
  are weakly holomorphic vector-valued modular forms, and providing necessary
  and sufficient conditions for when these are holomorphic. 
  
  In \cref{sec:analyzing} we study spaces of vector-valued modular forms
  associated to the \(\SLTZ\) representations that arise in
  \cref{sec:sl2}. We show the representations for forms of dimensions one
  and two are always congruence, while for dimension three there exists an
  infinite family of non-congruence representations. For all these 
  dimensions we provide explicit formulae for
  certain distinguished vector-valued modular forms from which all others can
  be generated. We also highlight the precise levels for which the space of
  all holomorphic vector-valued modular forms is obtained from torus
  \(1\)-point functions (the first few terms in the expansions of the
  distinguished \vvmf{s} at these levels are recorded in \cref{tab:2dtable,tab:3dtable}). For dimension four we describe the space of all torus
    \(1\)-point functions in those cases that a relevant space of holomorphic
    vector-valued modular forms is a cyclic module over the algebra of modular
    differential operators. 
  For general dimensions we identify levels of
  affine \(\SLA{sl}{2}\) for which the representation is non-congruence, if it
  is irreducible. 
  
  In \cref{sec:mtcs} we use the fact that categories of modules over rational
  \ctwo{} \voa{s} are modular tensor categories. The categorical counterpart
  to torus \(1\)-point functions are \(3\)-point coupling
  spaces and we study the action of \(\Bthree\) on these.
  The action on what corresponds to vacuum torus \(1\)-point functions is well
  studied and yields so called modular \(\modS\) and \(\modT\) matrices. 
  These are invariants of \correct{\mtcs{}} that are, however, known not to be complete
  invariants \cite{ModData17}. The full action of \(\Bthree\) on all torus \(1\)-point
  functions therefore has the potential to be a more fine-grained invariant.
  We derive
  explicit formulae for this action in terms categorical data (specifically
  twists and fusing matrices) in \cref{thm:BKmodaction}. We use these formulae
  to complement the results of \cref{sec:analyzing} by showing that a certain
  representation of dimension four is irreducible.
  
\subsection*{Acknowledgements}

The authors thank Luca Candelori, Thomas
Creutzig, Cameron Franc, Terry Gannon, Simon Lentner, Christopher
  Marks, Geoffrey Mason,
and Ingo Runkel for helpful discussions. 
SW's research is supported by the Engineering and Physical Sciences Research
Council (EPSRC) grant EP/V053787/1
and by the Alexander von Humboldt Foundation. 
  
\section{Modular properties of traces of intertwining operators}
\label{sec:intmodprops}

\subsection{Vector-valued modular forms}

We begin by fixing notation and recalling some facts about modular
forms. Let $\mathbb{N}$ denote the set of positive integers and $\mathbb{N}_0 = \mathbb{N}\cup \{0\}$.
Here and throughout, \(\tau\) will always lie in the complex upper half-plane
\(\HH=\set{z\in \CC\st \im(z)>0}\) and for a variable $x$ we denote
$\epi(x)=\ee^{2\pi \ii x}$. In particular, \(q=\epi(\tau)\) so that \(q\) lies in the
interior of the complex unit disk.
Let $\mathcal{M}$ denote the $\CC$-algebra of integral weight
holomorphic modular forms, and denote by $\mathcal{M}_k$ the subspace
of modular forms of weight $k$. Recall that, for $k\in \NN$, the
Eisenstein series of weight $2k$ are given by
 \begin{equation}
 \eis{2k} (\tau)= -\frac{B_{2k}}{(2k)!} +\frac{2}{(2k-1)!} \sum_{n = 1}^\infty \frac{n^{2k-1}q^n}{1-q^n},
\end{equation}
where $B_{\ell}$ denotes the $\ell$th Bernoulli number given by the generating series
\begin{equation}
 \sum_{\ell=0}^\infty B_\ell \frac{x^\ell}{\ell !} = \frac{x}{\ee^x-1}.
\end{equation}
For \(k\ge2\), \(\eis{2k}\) is an example of a holomorphic modular form of weight \(2k\). Moreover, \(\mathcal{M}\) is freely generated by \(\eis{4}\) and
\(\eis{6}\), that is, \(\mathcal{M}=\CC[\eis{4},\eis{6}]\).

Recall that, $\SLTZ$, the group of
integral \(2\times2\) matrices of unit determinant, is generated by
$\modS =\left(\begin{smallmatrix} 0 & -1 \\ 1 & 0\end{smallmatrix}\right)$
and $\modT =\left(\begin{smallmatrix} 1 & 1 \\ 0 &
    1 \end{smallmatrix}\right)$ and admits the presentation $\SLTZ =
\ang{\modS,\modT\st \modS^4= 1,\ (\modS\modT)^3=\modS^2}$. 
This group is furnished with a
natural action on \(\HH\) via
 \begin{equation}
 \gamma \tau = \frac{a\tau +b}{c\tau +d},
\end{equation}
where $\gamma = \left(\begin{smallmatrix} a&b \\
    c&d \end{smallmatrix}\right)\in \SLTZ$. Additionally,
for $k\in \RR$ we set
\begin{equation}
 j_k(\gamma ;\tau) = (c\tau +d)^k.
\end{equation}

 \begin{defn}
 A function $\nu \colon \SLTZ \to \{r\in \CC \mid \lvert r \rvert =1\}$
 is called a \emph{multiplier system for $\SLTZ$ of weight $k\in\RR$} if for any $A,B\in \SLTZ$ it satisfies
 \begin{equation}
   \nu (AB)j_k(AB;\tau) = \nu (A)\nu (B) j_k(A;B\tau)j_k(B;\tau).
   \label{eq:multrel}
\end{equation}
The \emph{cusp parameter} of $\nu$ is the unique $m\in \RR$ such that $0\leq m<12$ and $\nu (T)=\mathbf{e}(m/12)$.
\end{defn}
Note that the relation \eqref{eq:multrel} implies that multiplier
systems are uniquely characterised by their values on the generators
\(\modS\) and \(\modT\). The purpose of multiplier systems is to redefine projective representations of \(\SLTZ\) (specifically representations  of the braid group on three
strands \(\grp{B}_3=\ang{ \modS, \modT\st (\modS \modT)^3= \modS^2}\), where
\(\modS^4\) acts as a phase) so that they are no longer projective, as we shall see shortly.

Note that if \(\nu\) is a multiplier system of weight
\(k\), then it is also one of weight \(k+n\) for any \(n\in \ZZ\). The remainder of this subsection does not depend on
the choice of multiplier system one wishes to consider. Nevertheless, for use later we note that for each \(r\in \RR\) there exists a multiplier system \(\nu_r\) of weight
  \(r\) satisfying
  \begin{equation}
  \label{thm:multsys}
    \nu_r(\modT)=\epi\brac*{\frac{r}{12}},\qquad
    \nu_r(\modS)=\epi\brac*{\frac{-r}{4}},\qquad
    \nu_r(\modS\modT)=\epi\brac*{\frac{-r}{6}}
  \end{equation}
  (see, for example, \cite[\refProp 2.3.2]{Marks2010}).
  The multiplier system $\nu_r$ is precisely that which makes
  $\eta^{2r}$ transform
  as a modular form of weight $r$, where $\eta$ is Dedekind's eta function
\begin{equation}
  \eta (\tau) = q^{\frac{1}{24}} \prod_{n=1}^\infty \left(1-q^n\right).
\end{equation}
Note that $r$ is the cusp parameter for $\nu_r$ if and only if
$0\leq r <12$.

Multiplier systems allow us to define an action of \(\SLTZ\) on tuples
of functions (called vectors) on \(\HH\) via the following. 

\begin{defn}
 Let $\rho \colon \SLTZ \to \GL{d}$ be a $d$-dimensional representation of
 $\SLTZ$ 
 and consider
 holomorphic functions $f_1,\dots ,f_d\colon \HH \to \CC$ arranged into a vector
 \begin{equation}
   F= \left(f_1,\dots ,f_d\right)^t,
 \end{equation}
 where $\mathbf{x}^t$ denotes the transpose of a vector $\mathbf{x}$.
 \begin{enumerate}
 \item The vector \(F\) is a $d$-dimensional \emph{weakly holomorphic vector-valued modular form of weight
     $k\in \RR$} on $\SLTZ$ for the representation $\rho$ and a multiplier
   system $\nu$, if the following hold.
   \begin{enumerate}
   \item Each $f_j$ is meromorphic at the cusp $\ii\infty$.
     \label{itm:meratinf}
   \item For each
     $\gamma = \left(\begin{smallmatrix} a & b \\ c &
         d \end{smallmatrix}\right) \in \SLTZ$ we have
     \begin{equation}
       \label{eqn:VVMFTransformationProperty}
       F\vert_{k}^\nu \gamma = \rho (\gamma)F,
     \end{equation}
     where we define $\vert_{k}^\nu \gamma$ on each $f_j$ by
     \begin{equation}
       \left(f_j\vert_{k}^\nu \gamma \right)(\tau) = \nu(\gamma)^{-1}
       j_k(\gamma; \tau)^{-1} f_j\left(\gamma \tau\right)
       \label{eq:multaction}
     \end{equation}
     and extend the definition of $\vert_{k}^\nu$ component wise to $F$.
     \label{itm:action}
   \end{enumerate}
   \label{itm:merdef}
 \item
   The vector $F$ is a \emph{holomorphic vector-valued modular form} if it is a weakly holomorphic vector-valued modular form for which each $f_j$ is holomorphic at the cusp $i\infty$.
 \end{enumerate}
 \label{def:whvvmf}
\end{defn}

Here we see that if the multiplier system were omitted from the action
\eqref{eq:multaction}, then the \(\modS^4=1\) relation of \(\SLTZ\)
would not necessarily hold. This would therefore define a projective
action, or alternatively, an action of \(\grp{B}_3\).

For a fixed representation \(\rho \colon \SLTZ \to \GL{d}\) and multiplier
system \(\nu\) of weight
\(k\in\RR\), the corresponding vector spaces of weakly holomorphic and
holomorphic vector-valued modular forms of weight $k$ for representation
$\rho$ and multiplier system $\nu$ are denoted $\mathcal{M}^!(k,\rho,\nu)$ and
$\mathcal{H}(k,\rho ,\nu)$, respectively. As noted above, multiplier systems
only determine weights up to shifts by integers and these weights are always
in the same integer coset as the cusp parameter \(m\) of \(\nu\). We therefore denote by
\(\mathcal{M}^!(\rho,\nu)=\bigoplus_{n\in \ZZ} \mathcal{M}^!(m+n,\rho,\nu)\) and
\(\mathcal{H}(\rho ,\nu)=\bigoplus_{n\in \ZZ}\mathcal{H}(m+n,\rho ,\nu)\),
respectively, the spaces of all weakly holomorphic and
holomorphic vector-valued modular forms for the pair \((\rho,\nu)\). Further,
\(\mathcal{H}(\rho ,\nu)\) always admits a minimal weight
\(\mdwt\in\RR\) such that \(\mathcal{H}(\rho ,\nu)=\bigoplus_{n\in
  \NN_{0}}\mathcal{H}(\mdwt+n,\rho ,\nu)\) and \(\mathcal{H}(\mdwt-\ell,\rho
,\nu)=0\) for all \(\ell\in \NN\).

We will always assume that
$\rho (\modT)$ is diagonal with
 \begin{equation}
  \rho (\modT) = \diag\left(\epi\left(r_1\right),\dots
    ,\epi\left(r_d\right)\right),
  \label{eq:rhoT}
 \end{equation}
 for real numbers $r_1,\dots ,r_d$, in particular, \(\rho(T)\) is a unitary
   matrix.
This assumption can always be made for \voa{s} with a semisimple
representation theory (which is the case we will specialise to shortly) as the
intertwining operators can then be chosen without loss of generality to take values in simple  modules.
   A further simplifying assumption that is common in the number theory
   literature is that \(\rho(\modS^2)\) is a scalar matrix. Since \(\modS^2\)
   generates the centre of \(\SLTZ\) and has finite order (in the standard
   realisation of \(\SLTZ\) as integral \(2\times 2\) matrices with unit
   determinant, we have \(\modS^2=-1\)), the vector space
   underlying the representation \(\rho\) always admits a direct sum
   decomposition with \(\rho(\modS^2)\) acting as a scalar on each summand. This
   assumption therefore primarily simplifies the presentation of certain
   theorems and hence is not necessary. In the context of \voa{s},
   \(\modS^2\) carries the additional interpretation of being the
   \emph{charge conjugation involution} (the functor which assigns a module to
   its dual), so \(\rho(\modS^2)\) cannot be
   diagonal in a basis of intertwining operators which only take values in
   simple modules, if there are modules which are not self dual.
   For an account of the role of these assumptions in number theory one can consult, for example, \cite{Gannon-TheoryOfVVMF} or \cite{Marks2011}.
 In the latter, it is also assumed that $0\leq r_1,\dots ,r_d<1$, but this
 will not be required here. With these assumptions, if $F\in
 \mathcal{M}^!(k,\rho, \nu)$ we can replace
 \cref{def:whvvmf}.\ref{itm:meratinf}  with the condition that
 each $f_j$ has a Fourier expansion of the form 
 \begin{equation}
 \label{eq:VVMFq-expansion}
 f_j(\tau) = q^{\lambda_j} \sum_{n=0}^\infty a_n q^n,
 \end{equation}
 for some real numbers $\lambda_j$ (see, for example, \cite{Gannon-TheoryOfVVMF} for more discussion). As described in \cite{Marks2011}, a holomorphic vector-valued modular form $F$ requires each
 $f_j$ to have an expansion
 \eqref{eq:VVMFq-expansion} where
 \begin{equation}
 \label{eq:MinimalAdmissibleRelation}
  0\leq \lambda_j \equiv r_j +\frac{m}{12} \pmod{\ZZ},
 \end{equation} 
 and $m$ is the cusp parameter of the multiplier system $\nu$.
 \begin{defn}
   Let \(\rho \colon \SLTZ\to \GL{d}\) be a representation such that \(\rho(T)\) is
   diagonal and unitary as in \eqref{eq:rhoT} and \(\nu\) a multiplier system.
   A set of non-negative real numbers $\{\lambda_1,\dots ,\lambda_d\}$
   satisfying \eqref{eq:MinimalAdmissibleRelation} is called an
   \emph{admissible set} for $(\rho,\nu)$. The \emph{minimal admissible set}
   for $(\rho ,\nu)$ is the unique admissible set which
   additionally satisfies $\lambda_j <1$ for each $j$.
   \label{def:adm}
 \end{defn}
  As pointed out in \cite{Marks2011}, and which follows from
  \eqref{eq:VVMFq-expansion} and \eqref{eq:MinimalAdmissibleRelation} above,
  for a minimal admissible set $\{\lambda_1,\dots ,\lambda_d\}$ we have that
  every non-zero $F\in \mathcal{H}(\rho,\nu)$ has the form
 \begin{equation}
 \label{eq:FormOfVVMF}
 \begin{pmatrix} q^{\lambda_1 +n_1} \sum_{n=0}^\infty a_1(n)q^n 
 \\ \vdots
 \\  q^{\lambda_d +n_d} \sum_{n=0}^\infty a_d(n)q^n
 \end{pmatrix},
 \end{equation}
 where for each \correct{$j=1,\dots,d$} we have $a_j(0)\not =0$ and $n_j$ are non-negative integers. \correct{Note that \(a_j(n)\) refers to the \(n\)th coefficient for the \(j\)th entry in the column vector.}

 To provide conditions for a lower bound on the minimal weight $\mdwt$, among other things, we require the modular derivative in weight $k\in \RR$, which is defined as
 \begin{equation}
 \label{eqn:ModularDerivative}
 \partial =\partial_k = \frac{1}{2\pi i } \frac{d}{d\tau} +k\eis{2}(\tau),
 \end{equation}
on weight \(k\) (vector-valued) modular forms and is then extended
  linearly.
 The modular derivative increments the weight of (vector-valued) modular forms by two. In
   particular, the homogeneous subspaces of  the ring of
  integral weight holomorphic modular forms are related by \(\partial_k
  \mathcal{M}_k \subset \mathcal{M}_{k+2}\).
For $n\in \mathbb{N}$ and $\phi \in \mathcal{M}_k$ we let $\partial^n \phi$ denote the composition of operators $\partial_{k+2n}\circ \partial_{k+2(n-1)} \circ \cdots \circ \partial_{k+2} \circ \partial_k \phi$. This allows us to consider an order $n\in \mathbb{N}$ monic modular differential equation in weight $k\in \mathbb{R}$, which is an ordinary differential equation of the form
 \begin{equation}
 \left(\partial_k^n + \sum_{j=0}^{n-2} \phi_{2(n-j)} \partial_k^{j}\right) f =0
 \end{equation}
 in the disk $\lvert q\rvert <1$, where $\phi_j \in \mathcal{M}_j$ for each $j$. For more details about monic modular differential equations see, for example, \cite{Marks2011, Mason-VVMFandDiffOps}.
  
The modular derivative can be adjoined to the algebra of integral weight
  modular forms \(\mathcal{M}\) to form
 \begin{equation}
   \mathcal{R} = \left\{ \phi_0+\phi_1\partial + \cdots +\phi_n \partial^n \mid \phi_i \in \mathcal{M},n\geq 0\right\},
   \label{eq:moddiffops}
 \end{equation}
 the skew polynomial ring of modular differential operators, where addition is
 defined component wise, and multiplication is characterised
 by $\partial
 \cdot \phi = \phi \partial +\partial_k \phi$ for $\phi \in
 \mathcal{M}_k$. In fact, for any $F \in \mathcal{M}^!(k,\rho,\nu)$, defining
 $\partial F$ to be $\partial$ applied component wise, we find $\partial F\in
 \mathcal{M}^!(k+2,\rho ,\nu)$. Similarly, for $\phi \in \mathcal{M}_k$ we let
 $\phi F$ denote the vector $F$ with $\phi$ multiplied component wise.
Thus both \(\mathcal{M}^!(\rho,\nu)\) and \(\mathcal{H}(\rho,\nu)\) are left
$\mathcal{R}$-modules with $\mathcal{H}(\rho ,\nu)$ as an \(\mathcal{R}\)-submodule.
  
 We conclude this section with convenient criteria for determining when the
 components of a holomorphic vector-valued modular form span the solution space of a monic
 modular differential equation and for determining the minimal weight \(\mdwt\)
 of the space \(\mathcal{H}(\rho,\nu)\) of holomorphic vector-valued modular forms.
\begin{thm}[Mason {\cite{Mason-VVMFandDiffOps}}, Marks {\cite[\refThm 2.8]{Marks2011}}]
  Let \(\rho \colon \SLTZ\to \GL{d}\) be a representation such that \(\rho(T)\) is
   diagonal and unitary as in \eqref{eq:rhoT}, \(\nu\) a multiplier system,
   and $\{\lambda_1,\dots ,\lambda_d\}\subset \RR$ the minimal admissible set for
   \((\rho,\nu)\).
   Consider a holomorphic vector-valued modular form
   $F\in \mathcal{H}(\correct{p},\rho,\nu)$ which must therefore have an expansion of
   the form \eqref{eq:FormOfVVMF} such that the leading exponent of the
   \(j\)th component is \(\lambda_j+n_j\), \(n_j\in \NN_0\).
   If the components of \(F\) are linearly independent over \(\CC\), then the
   weight \(\mdw\) is bounded below by the inequality 
 \begin{equation}
  \mdw \geq \frac{12\left(\sum_j \left(\lambda_j +n_j\right)\right)}{d} +1-d,
\label{eq:bigweightineq}
\end{equation}
 and equality holds if and only if the components of $F$ span the solution
 space of a monic modular differential equation.
 In particular, the minimal weight $\mdwt$ of \(\mathcal{H}(\rho,\nu)=
 \bigoplus_{n\in \NN_0}\mathcal{H}(\mdwt+n,\rho,\nu)\) 
 satisfies
 \begin{equation}
   \mdwt \geq \frac{12\left(\sum_j \left(\lambda_j \right)\right)}{d} +1-d.
   \label{eq:weightineq}
 \end{equation}
\label{thm:LowestWeightOfVVMF}
\end{thm}

\subsection{Modularity of torus 1-point functions}

The purpose of this section is to fix notation relating to \voa{s} and to
introduce torus 1-point functions, which can be defined in a number of ways. 
Here we shall define and construct them as suitable traces of
intertwining operators. We mostly follow the conventions of \cite{FreVer01}, however, since few textbooks cover
intertwining operators, we introduce these in greater detail. See
\cite{HuaLog} for an exhaustive account of intertwining operators and the
tensor structures that arise from them.

Let \((V,Y,\wun,\omega)\) be a vertex operator algebra, where \(V\)
denotes the underlying vector space, \(Y\) the field map, \(\wun\) the
vacuum vector, and \(\omega\) the conformal vector. The central charge
of the Virasoro algebra generated by the field expansion of
\(Y(\omega,z)=\sum_{n\in\ZZ} L_n z^{-n-2}\) will be denoted \(\mathbf{c}\). Further, let
$(U,Y_U)$ be a $V$-module, with \(U\) the underlying vector space
and \(Y_U\) the field map (or action) representing the \voa{} \(V\). 
We will always assume that modules are graded by
generalised \(L_0\) eigenvalues, that is,
\begin{equation}
  U=\bigoplus_{n\in\CC} U_n,\qquad U_n=\set{u\in U\st \exists m\in
    \NN, (L_0-n)^mu=0}.
  \label{eq:confgrading}
\end{equation}
For a homogeneous element $u\in U_n$, we denote the \emph{conformal weight}
  \(n\) of \(u\) by \(\operatorname{wt}(u)=n\).

\begin{defn}
    Let \((V,Y,\wun,\omega)\) be a \voa{} and let \((U_i,Y_{U_i}),\ i=1,2,3\) be \(V\)-modules. Consider a linear map
  \begin{align}
    \mathcal{Y} \colon U_1\otimes U_2&\to U_3\{z\}[\log (z)]\nonumber\\
    u_1\otimes u_2&\mapsto \mathcal{Y}(u_1,z) u_2=\sum_{\substack{s\in \CC\\
    t\ge0}} (u_1)_{s,t} u_2 z^{-s-1}\log( z)^t,
  \end{align}
  where \(\{z\}\) denotes unbounded power series with arbitrary complex
  exponents and \(\log (z)\) is a formal variable satisfying the relation
  \(\frac{\dd}{\dd z} \log (z)=1/z\). Such a map is called an \emph{intertwining operator of type \(\binom{U_3}{U_1\ U_2}\)}, if it satisfies the following conditions.
  \begin{enumerate}
  \item Truncation: For fixed \(u_1\in U_1\), \(u_2\in U_2\) and \(s\in\CC\), \(t\ge0\),
    \begin{equation}
      (u_1)_{s+\ell,t} u_2=0
    \end{equation}
    for sufficiently large \(\ell\in \ZZ\).
  \item \(L_{-1}\)-derivative: For any \(u_1\in U_1\) and \(u_2\in U_2\),
    \begin{equation}
      \dfrac{\dd}{\dd z}\mathcal{Y}(u_1,z)u_2=\mathcal{Y}(L_{-1}u_1,z)u_2.
    \end{equation}
  \item Jacobi identity: For any \(v\in V, \ u_1\in U_1,\ u_2\in U_2\),
    \begin{align}
      z_0^{-1} \delta \brac*{\frac{z_1-z_2}{z_0}} Y_{U_3}(v,z_1)\mathcal{Y}(u_1,z_2)u_2
     &= z_0^{-1} \delta\brac*{\frac{-z_2 +z_1}{z_0}} \mathcal{Y}(u_1,z_2) Y_{U_2}(v,z_1)u_2 \nonumber\\
      &\quad + z_2^{-1} \delta \brac*{\frac{z_1 - z_0}{z_2}} \mathcal{Y}(Y_{U_1}(v,z_0)u_1,z_2)u_2,
        \label{eq:jacID}
    \end{align}
    where
    \begin{equation}
      \delta\brac*{\frac{a-b}{c}}=
      \sum_{\substack{r\in \ZZ\\ s\ge 0}}
      \binom{r}{s}(-1)^s a^{r-s} b^{s} c^{-r}
    \end{equation}
    is the algebraic Dirac delta function.
  \end{enumerate}
\end{defn}
We will soon specialise to rational \ctwo{} \voa{s}. Among other
helpful properties, such \voa{s} only admit modules for which \(L_0\) acts
semisimply with finite-dimensional eigenspaces, and all eigenvalues
are rational, bounded below, and discrete. Furthermore, intertwining operators do
not contain any \(\log (z)\) terms.
Given two \(V\)-modules \((U,Y_U)\) and \((W,Y_W)\), we are particularly
interested in intertwining operators \(\mathcal{Y}\) of type
\(\binom{W}{U\ W}\). Note that we use the symbol \(\binom{W}{U\ W}\) both to denote the type of an intertwining operator and to denote the vector space of all intertwining operators of that type.
For any \(L_0\) eigenvector \(u\in U_{\wt (u)}\),
we define the \(\mathcal{Y}\) zero mode \(o^{\mathcal{Y}}(u)\) of \(u\) to be the coefficient
of \(z^{-\wt (u)}\) in the expansion of \(\mathcal{Y}(u,z)\). This zero
mode preserves generalised \(L_0\) eigenvalues, that is,
\(o^{\mathcal{Y}}(u)(W_m)\subset W_m\) for all generalised
\(L_0\) eigenvalues \(m\in\CC\).

\begin{defn}
  Let \((V,Y,\wun,\omega)\) be a \voa{}, \((U,Y_{U}),\ (W,Y_{W})\) be
  \(V\)-modules, and \(\mathcal{Y}\) be an intertwining operator of type
  \(\binom{W}{U\ W}\). The \textit{torus 1-point function associated to
  \(\mathcal{Y}\)} is the trace
  \begin{equation}
    \ptf^{\mathcal{Y}} (u,\tau) =
 \operatorname{tr}_{W} o^{\mathcal{Y}}(u)
 q^{L(0)-\frac{\mathbf{c}}{24}},
 \qquad u\in U.
 \label{eq:q-expansionOfTraceFunction}
  \end{equation}
  If the intertwining operator \(\mathcal{Y}\) is clear from context, we will omit \(\mathcal{Y}\) as a label for zero modes and
  torus 1-point functions.
\end{defn}

Beyond the standard expansions of fields, we will also need to consider
transformed expansions
\begin{equation}
\label{eqn:SquareBracket}
  Y[a,z] = Y(a,\ee^z-1)\ee^{z \wt (a)}=\sum_{n\in \mathbb{Z}}a_{[n]}z^{-n-1},\qquad a\in V_{\wt (a)},
\end{equation}
and extended linearly,
which implies the formula
\begin{equation}
  a_{[n]}=\Res_z \brac*{Y(a,z)(\log(1+z))^n(1+z)^{\wt (a)-1}}.
\end{equation}
For example, for $a\in V_{\wt (a)}$ we have
 \begin{equation}
 a_{[0]} = \sum_{j=0}^\infty \binom{\wt (a)-1}{j} a_j ,
 \end{equation}
 so that if $a\in V_1$ we obtain $a_{[0]}=a_{0}$.
 In fact, the map \eqref{eqn:SquareBracket} gives $V$ another structure of a
 \voa{} of central charge $\mathbf{c}$ with the same vacuum vector and
 conformal element $\tilde{\omega}=\omega - \frac{\mathbf{c}}{24} \wun$
(see \cite[\refSec 4]{Zhu} for details). Similar to above, defining $L_{[n]}$ via $Y[\tilde{\omega},z]=\sum_{n\in \mathbb{Z}} L_{[n]}z^{-n-2}$ gives us a
\emph{square bracket grading} 
\begin{equation}
  U=\bigoplus_{n\in\CC} U_{[n]},\qquad U_{[n]}=\set{u\in U\st \exists m\in
    \NN, (L_{[0]}-n)^mu=0}.
  \label{eq:ModulesqGrading}
\end{equation}
If $u\in U_{[n]}$ we write $\operatorname{wt}[u]=n$.
In the case where the \voa{} \(V\) is rational and \ctwo{}, we again have
that \(L_{[0]}\) acts semisimply with finite-dimensional eigenspaces, and that all eigenvalues
are rational, bounded below, and discrete. Additionally, if we assume
that \(U\) is simple then the minimal \(L_0\) and \(L_{[0]}\)
eigenvalues will be equal and denoted \(h_U\), which is a
rational number \cite{DLM-Orbifold}, and is called the \emph{conformal weight
  of the module \(U\)}.
Further, for \(n\in \mathbb{N}_0\) we have
\(U_{[h_U+n]}\subset \bigoplus_{m=0}^n U_{h_U+m}\) and \(U_{h_U+n}\subset
\bigoplus_{m=0}^n U_{[h_U+m]}\).

For later use we prepare some helpful identities involving torus 1-point
functions and square bracket expansions.

\begin{prop}[Zhu {\cite{Zhu}}, Miyamoto {\cite[\refProp 3.1 and 3.3]{miyamoto2000intertwining}}]
  Let \(V\) be a \voa{}, \(U,W\) be \(V\)-modules, and \(\mathcal{Y}\)
  an intertwining operator of type \(\binom{W}{U\ W}\). Then for any
  $a\in V$ and $u \in U$ we have
 \begin{equation}
   \ptf^{\mathcal{Y}} (a_{[0]}u,\tau)=0
   \label{eq:1pt[0]result}
 \end{equation}
 and 
 \begin{equation}
 \ptf^{\mathcal{Y}} (a_{[-1]}u,\tau)= \operatorname{tr} \vert_{W} o(a)o(u) q^{L_0 -\frac{\mathbf{c}}{24}}
 +\sum_{\ell =1}^\infty \eis{2\ell}(\tau) \ptf^{\mathcal{Y}} (a_{[2\ell -1]}u,\tau).
 \label{eq:1pt[-1]result}
\end{equation}
 \label{thm:1ptrels}
\end{prop}

The identity \eqref{eq:1pt[0]result} can be specialised and refined as follows.

\begin{prop}
  Let $V$ be a \voa{} and $U$ a $V$-module with a decomposition
  into generalised \(L_{[0]}\) eigenspaces as in
  \eqref{eq:ModulesqGrading}. Suppose that for any $x,y\in V_1$ the binary operation \([x,y]=\Res_z Y(x,z)y=x_0y\) furnishes $V_1$ with the structure of a finite-dimensional
  reductive Lie algebra and thus each homogeneous space $U_{[m]}$ is a
  module over \(V_1\). \correct{If \(U_{[m]}\) is semisimple over \(V_1\), let} \(U_{[m]}=U_{[m]}^{\text{triv}}\oplus
  U_{[m]}^{\text{non-triv}}\) be the unique decomposition into
  the maximal trivial submodule \(U_{[m]}^{\text{triv}}\) and its complement
  \(U_{[m]}^{\text{non-triv}}\) containing all non-trivial simple summands.
  Then for any \(u
  \in U_{[m]}^{\text{non-triv}}\), any $V$-module $W$, and any intertwining operator $\mathcal{Y}$ of type $\binom{W}{U\ W}$, we have
 \begin{equation}
   \ptf^{\mathcal{Y}} (u,\tau) = 0.
 \end{equation}
  \label{prop:trivialmodulerequiredfornonzerotrace}
\end{prop}

\begin{proof}
  Recall that for \(a\in V_1\) the square \correct{bracket} and non-square \correct{bracket} zero modes
  coincide, that is, \(a_{[0]}=a_0\).
  Since by assumption $V_1$ is a reductive Lie algebra, we can
  decompose $U_{[m]}$ into a direct sum of irreducible $V_1$-modules
  and group them into trivial and non-trivial ones.
  In particular
  \begin{equation}
    U_{[m]}^{\text{triv}}=\set{u\in U_{[m]}\st \forall a\in V_1, a_0
      u=0},\qquad
    U_{[m]}^{\text{non-triv}}=\set{a_0 u\st a\in V_1, u\in U_{[m]}}.
  \end{equation}
  Therefore, if \(u\) lies in a
  simple non-trivial submodule of \(U_{[m]}\), then it also lies in
  the image of some \(a\in V_1\) and thus by \cref{thm:1ptrels} the
  result follows.
\end{proof}

Note that a sufficient condition for \(V_1\) being a Lie algebra under
the bracket given above is \(\dim (V_0)=1\) and \(\dim (V_{n})=0\) for \(n<0\).
If in addition to being rational \(V\) is \ctwo{} and \(L_{1} V_1=0\), then
\(V_1\) is reductive \correct{and all conformal weight spaces are finite-dimensional
semisimple modules over \(V_1\)} \cite{DM-Effective,MasCIr14}, which is the case for the examples we shall consider later.
 
For the remainder of this article we assume that $V$ is a \voa{} for which the conformal weights are bounded below by \(0\), the conformal weight \(0\) space \(V_0\) is $1$-dimensional, $V$ is \ctwo{}, the contragredient or graded dual
\(V^\ast\) satisfies \(V^\ast\cong V\), and $V$ is rational, i.e., the
category of \correct{finitely generated 
modules with conformal weights bounded below},
\(\rep{V}\), is semisimple.
The \ctwo{ness} of $V$ implies that \(\rep{V}\) admits only a finite
number of inequivalent simple $V$-modules
$V=W^1,\dots ,W^{d_V}$ for some $d_V\in \NN$ \cite{Zhu} (see also \cite{DLM-Twisted}). In this notation, for each $\mu \in \{1,\dots, d_V\}$, we let
$h_\mu$ denote the conformal weight of $W^\mu$, where $h_1=0$.
Further, \ctwo{ness} also implies that the central charge and
conformal weights of all modules are rational \cite[Corollaries 5.10 and 5.11]{Miyamoto-C2} (see also \cite[\refThm 1.1]{DLM-Orbifold}). Note that the assumptions \(\dim
(V_0)=1\), \(V\cong V^\ast\), and that the conformal weights are bounded below
by \(0\) are not required for the  space of torus \(1\)-point functions 
to be closed under the action of the modular group \cite{Zhu} (see also
\cite{DLM-Orbifold}). 
They are necessary to prove that \(\rep{V}\) is rigid and, additionally, 
a modular tensor category  \cite{HuaVer08}.

Given a simple \(V\)-module \(W^\lambda\) we introduce the vector space of intertwining operators
\begin{equation}
  \mathcal{I}_\lambda=\bigoplus_{\mu}\binom{W^\mu}{W^\lambda\ W^\mu}.
  \label{eq:iopspace}
\end{equation}
Let $N_{\lambda,\mu}^\mu$ denote the dimension of the space of intertwining operators of type $\binom{W^\mu}{W^\lambda\ W^\mu}$ so that \(\dim (\mathcal{I}_\lambda) = \sum_{\mu} N^\mu_{\lambda,\mu}\). Note
 that in the special case when \(W^\lambda\) is the \voa{} \(V\), we have
 \(N^\mu_{\lambda,\mu}=1\) for all labels \(\mu\) and there is a
 distinguished basis given by the field maps \(Y_{W^\mu}\), that is, the
 action of \(V\) on \(W^\mu\).

The space of \(1\)-point functions, \(\mathcal{C}_1 (W^{\lambda})\), with insertion from a simple
module \(W^\lambda\) admits a number of characterisations in various level of generality, 
however, for rational \ctwo{} \voa{s} this space can always be realised as the span 
 \begin{equation}
 \mathcal{C}_1 \left(W^{\lambda}\right) = \set{ \ptf^{\mathcal{Y}} (- ,\tau)
   \st  \mathcal{Y}\in \mathcal{I}_\lambda},
\end{equation}
see \cite[\refThm 5.1]{Yamauchi-IntertwinerModularity} for details, which we use as the definition here, for simplicity.
Further, we define the space of torus 1-point functions
evaluated at $u\in W^{\lambda}$ to be
 \begin{equation}
 \mathcal{C}_1^u \left(W^{\lambda}\right) = \set{ \ptf^{\mathcal{Y}} (u ,\tau) \st  \mathcal{Y}\in \mathcal{I}_\lambda}.  
\end{equation}
Bounds on dimensions are then given by
 \begin{equation}
 \dim \left(\mathcal{C}_1^u \left(W^{\lambda}\right)\right) \leq
  \dim \left(\mathcal{C}_1 \left(W^{\lambda}\right)\right)\leq
 \dim \left(\mathcal{I}_\lambda\right) = \sum_{\mu =1}^{d_V} N_{\lambda ,\mu}^{\mu}.
  \label{eq:DimOf1ptSpaceBound}
\end{equation}

\begin{thm}[{Miyamoto \cite[\refThm 5.1]{miyamoto2000intertwining}, Yamauchi
    \cite[\refThm 5.1]{Yamauchi-IntertwinerModularity}, Huang
    \cite[\refThm 7.3]{Huang-IntModInv}}]
  Recall the multiplier system \(\nu_{r}\) given in \eqref{thm:multsys}. Let \(V\) be a rational \ctwo{} \voa{.} Then for any simple module \(W^\lambda\)
  and any homogeneous vector
  \(u\in W^\lambda_{[\wt [u]]}\), every torus 1-point function
  \(\ptf(u,\tau)\in \mathcal{C}_1^u (W^{\lambda})\)
  evaluated at \(u\) is a holomorphic function on \(\HH\). 
  For \(\gamma\in \SLTZ\) and \(u\in W^\lambda_{[\wt [u]]}\), 
  \begin{equation}
    \ptf^{\mathcal{Y}} (u, \tau)|_\gamma=\nu_{\wt [u]} (\gamma)^{-1}j_{\wt [u]}(\gamma ;\tau)^{-1} \ptf^{\mathcal{Y}} (u, \gamma \tau)  \in \mathcal{C}_1^u \left(W^{\lambda}\right)
  \end{equation}
  defines an action of \(\SLTZ\) on \(\mathcal{C}_1^u
  (W^{\lambda})\) which lifts to an action on \(\mathcal{C}_1
  (W^{\lambda})\).
  \label{thm:modinv}
\end{thm}

Recall that the weight of a multiplier system can be freely shifted by
integers and that for \(u\) as in the theorem above \(\wt
[u]-h_\lambda\in\ZZ\). We can therefore define the right action of \(\SLTZ\) on \(\mathcal{C}_1^u (W^{\lambda})\) to be
\begin{equation}
 \left(\ptf^{\mathcal{Y}} \vert_{\operatorname{wt}[u]}^\lambda \gamma \right)(u,\tau)= \nu_{h_\lambda} (\gamma)^{-1} j_{\operatorname{wt}[u]} (\gamma ;\tau)^{-1} \ptf^{\mathcal{Y}}(u, \gamma \tau),
\end{equation}
thus making the multiplier system independent of the choice of vector \(u\in W^\lambda\).

To transition from considering right \(\SLTZ\) actions on \(\mathcal{C}_1^u
(W^{\lambda})\) to vector-valued modular forms, we need to
  choose elements in \(\mathcal{C}_1^u
(W^{\lambda})\) to form the components of a vector, which we shall
now do.
Let $\mathcal{B}_\lambda^\mu$ be a basis for the space of intertwining operators of type $\binom{W^\mu}{W^\lambda\ W^\mu}$ which is the empty set if \(\lambda,\mu\) are such that \(N^\mu_{\lambda,\mu}=0\) as the corresponding space of intertwining operators will vanish. Let
\begin{equation}
  \intset_\lambda = \bigcup_\mu \mathcal{B}_\lambda^\mu
  \label{eq:intset}
\end{equation}
be the union of all these bases for fixed \(\lambda\), thus forming a basis
of \(\mathcal{I}_\lambda\). The space of intertwining operators
\(\mathcal{I}_\lambda\) and the spaces of torus 1-point functions \(\mathcal{C}_1
(W^{\lambda})\) and \(\mathcal{C}_1^u (W^{\lambda})\)
are related by the linear maps
\begin{align}
  \operatorname{tr}^\lambda \colon \mathcal{I}_\lambda&\to \mathcal{C}_1\left(W^{\lambda}\right),&
  \operatorname{ev}_u \colon \mathcal{C}_1\left(W^{\lambda}\right)&\to \mathcal{C}_1^u \left(W^{\lambda}\right),\nonumber\\
  \mathcal{Y}&\mapsto \operatorname{tr} o^{\mathcal{Y}}(-) q^{L(0)-\frac{\mathbf{c}}{24}}&
                                                                                           f(-,\tau)&\mapsto f(u,\tau)
 \label{eq:intmaps}                                                                                                 
\end{align}
that is, the first map is the taking of traces and the second is evaluation at
the vector \(u\in W^\lambda\). These maps are surjective by construction, hence
their composition is too, yet they need not be injective. In particular,
certain choices of \(u\in W^\lambda\) can lead to large kernels.
For example, if $V$ is  the simple affine vertex operator algebra constructed from 
 $\SLA{sl}{3}$ at level $3$, with $W^\lambda$ chosen to be $V$, and $u = \wun$, then we have $\lvert \intset_\lambda \rvert =10$
 while \correct{$\dim (\mathcal{C}^{\textbf{1}}_1(V))=6$}. That is, there are 10 simple
 modules up to equivalence, yet the span of characters is only $6$-dimensional. Indeed,
 \cref{prop:trivialmodulerequiredfornonzerotrace} shows that there can exist non-zero
 \(u\in W^\lambda\) for which \(\operatorname{ev}_u\) is the zero map. 

 Setting  \(\delta(\lambda)=\dim
 \mathcal{C}_1 (W^{\lambda})\), let \(\Delta_\lambda\subset \operatorname{tr}^\lambda(\intset_\lambda)=\{\ptf^1,\dots
 \ptf^{\delta(\lambda)}\}\) be a
 linearly independent subset of the image of the basis \(\intset_\lambda\) and
 hence a basis of \(\mathcal{C}_1 (W^{\lambda})\), and define the
 vector \(\vvptf_\lambda=(\ptf^1,\dots, \ptf^{\delta(\lambda)})^t\). Then
 \cref{thm:modinv} can be restated as follows.
\begin{thm}
  For any $u\in W^\lambda_{[\operatorname{wt}[u]]}$ the space
  $\mathcal{C}_1^{u} (W^{\lambda})$ carries a $\modT$-unitarisable
  representation $\rho_\lambda \colon \SLTZ \to \operatorname{GL}({\delta(\lambda)},\CC)$ such that $\vvptf_\lambda (u,\tau)$
  is a $\delta(\lambda)$-dimensional weakly holomorphic vector-valued modular form of weight $\wt [u]$, representation $\rho_\lambda$, and multiplier system $\nu_{h_\lambda}$. That is, $\vvptf_\lambda (u,\tau) \in \mathcal{M}^!(\wt [u],\rho_\lambda , \nu_{h_\lambda})$. 
  For each component \(\ptf^j \) of \(\vvptf_\lambda \) let \(\mu_{j}\) be the corresponding module label, that is, \(\ptf^j\) is the trace of an intertwining operator of type \(\binom{W^{\mu_{j}}}{W^\lambda\ W^{\mu_{j}}}\). Then additionally
    \begin{equation}
  \label{eqn:RhoT}
  \rho_\lambda \left(T\right) = \diag\left\{\mathbf{e}\left(r_{1}\right),
    \dots ,\mathbf{e}\left(r_{\delta(\lambda)}\right) \right\},
  \end{equation}
  where
  \begin{equation}
    r_j = h_{\mu_j} -\frac{\mathbf{c}}{24} - \frac{h_\lambda}{12},\qquad
    1\le j\le \delta(\lambda).
    \label{eq:Texponents}
  \end{equation}
  Moreover, if $h_{\mu_j} - \mathbf{c}/24 \geq 0$ for all $1\le j\le \delta(\lambda)$, then $\vvptf_\lambda (u,\tau) \in \mathcal{H}(\wt [u],\rho_\lambda , \nu_{h_\lambda})$. 
 \label{thm:TraceFunctionIsVVMF}
\end{thm}
 
\begin{proof}
  By construction every basis element \(\ptf^j\in \Delta_\lambda\) is the image of an
  intertwining operator that takes values in a simple module. The exponents of
  \(q\) in the series expansion of \(\ptf^j\) will therefore only differ by
  integers and hence the matrix for \(\modT\) will be diagonal and all
  diagonal entries will be complex numbers of modulus 1. 
 That $\vvptf_\lambda (u,\tau)$ is a $\lvert \intset_\lambda\rvert$-dimensional weakly holomorphic vector-valued modular
 form of weight $\wt [u]$,
 representation $\rho_\lambda$, and multiplier system $\nu_{h_\lambda}$
 follows from Theorem \ref{thm:modinv} (for additional details, see
 \cite[\refProp 2.5.2]{Marks2010}). 
  Taking $\gamma =\modT$ in \eqref{eqn:VVMFTransformationProperty} gives
 \begin{equation}
 \rho_\lambda \left(\modT\right) \vvptf_\lambda \left(u ,\tau \right)=
 \nu_{h_\lambda}(\modT)^{-1} j_{\operatorname{wt}[u]}(\modT ;\tau)^{-1} \vvptf_\lambda \left(u , \tau +1 \right) 
  = \epi\brac*{-\frac{h_\lambda}{12}-\frac{\mathbf{c}}{24}}
  \diag\set*{\epi\brac*{h_{\mu_1}},\dots,\epi\brac*{h_{\mu_{\delta(\lambda)}}}} \vvptf_\lambda \left(u ,\tau\right),
 \end{equation}
 and thus \eqref{eqn:RhoT} and \eqref{eq:Texponents}.
 Meanwhile, $\vvptf_\lambda (u,\tau)$ is
 holomorphic at $\ii\infty$ if and only if each component function \(\ptf^j\) is holomorphic at
 $\ii\infty$, and this is true if and only if $h_{\mu_j} - \mathbf{c}/24 \geq
 0$ for all $1\le j\le \delta(\lambda)$.
 \end{proof}

We stress that choices were made to construct the vector $\vvptf_\lambda (u,\tau)\in \mathcal{M}^!(\rho_\lambda ,\nu_{h_\lambda})$ from $\mathcal{C}_1(W^{\lambda})$.
 For example, for any $U\in \operatorname{GL}_d(\CC)$ the components of the vector $U\vvptf_\lambda (u,\tau)$ will also from a basis of $\mathcal{C}_1(W^{\lambda})$ and hence give rise to an equivalent representation $\rho_U$ related to the previously constructed representation via $\rho_U (\gamma)=U\rho_\lambda (\gamma)U^{-1}$ for all $\gamma \in \SLTZ$. Furthermore,
 $U\vvptf_\lambda (u,\tau)\in \mathcal{M}^!(\wt
 [u],\rho_U,\nu_{h_\lambda})$. The association of
 \(\mathcal{C}_1(W^{\lambda})\) to vector-valued modular forms is therefore
 only determined up to a choice of basis.

A natural question to ask, is if there is a discrepancy between the dimension
of the space of intertwining operators \(\mathcal{I}_\lambda\) and that of the
unevaluated torus 1-point functions \(\mathcal{C}_1(W^{\lambda})\), or equivalently,
if the trace map \(\operatorname{tr}^\lambda\) in \eqref{eq:intmaps} has a non-trivial
kernel. In \cite[\refThm
5.3.1]{Zhu} it was shown that the kernel is trivial
in the special case $W^\lambda
= V$. However, it is currently not known, if this is true or false for general \(W^\lambda\).
A sufficient condition for the kernel being trivial is the existence of a vector \(u\in W^\lambda\) such that
the image of the basis \(\intset_\lambda\) under the composition
\(\operatorname{ev}_u\circ \operatorname{tr}^\lambda\) is linearly independent, as in 
this case the inequalities \eqref{eq:DimOf1ptSpaceBound} are all equalities.
Such vectors will also be shown to exist below (not assuming $W^\lambda = V$)
in the example of the simple affine vertex operator algebra 
constructed from $\SLA{sl}{2}$ at any non-negative integer
level.

\section{General results on spaces of $1$-point functions}
\label{sec:1ptSpace}

In this section we consider the behaviour of evaluated torus 1-point functions, when
one lets the insertion vector \(u\in W^\lambda\) vary.
    Let
 \begin{equation}
 	\mathcal{V}\left(\rho_\lambda\right)_{n} = \operatorname{span}_{\CC} \left\{\vvptf_\lambda (u,\tau) \st u\in W^\lambda_{[h_\lambda +n]} \right\}
 \end{equation}
  for $n\in \mathbb{N}_0$ and
\begin{equation}
  \mathcal{V} \left(\rho_\lambda\right) = \bigoplus_{n=0}^\infty \mathcal{V}\left(\rho_\lambda\right)_{n}.
  \label{eq:spaceofvecmodforms}
\end{equation}
That is, $\mathcal{V} \left(\rho_\lambda\right)$ is the space of all
evaluations of the \(\vvptf(-,\tau)\) at any \(u\in W^\lambda\).

\correct{\begin{defn}
  Let \(u\in W^\lambda\) be an \(L_{[0]}\)-eigenvector.
  \begin{enumerate}
  \item We define the following
    subspaces of \(W^\lambda\).
    \begin{align}
      &U(\mathcal{L})u,\nonumber\\
      &\operatorname{Vir}(u)=U(\mathcal{L}_{[<0]})u=\spn{L_{[-n_1]}\cdots L_{[-n_\ell]}u \st n_1,\dots
        ,n_\ell \in \NN,  \ell\in \NN_0},\nonumber\\
      &\mathcal{N}(u)= U(\mathcal{L}_{[>0]})\mathcal{L}_{[>0]}u
        =\spn{L_{[n_1]}\cdots L_{[n_\ell]}u \st n_1,\dots
        ,n_\ell \in \NN,  \ell\in \NN},
    \end{align}
    these are, respectively, the Virasoro submodule generated by the
    \(L_{[n]}, n\in\ZZ\) from \(u\);
    \correct{spanned by words in the negative modes \(L_{[-n]}, n\in \NN\) applied to
      \(u\), including the empty word}; and the subspace 
    \correct{spanned by non-empty words in the positive modes \(L_{[n]}, n\in \NN\) applied to \(u\).}
    These subspaces give rise to the (non-direct) sum
    decomposition
    \begin{equation}
      U(\mathcal{L})u = \operatorname{Vir}(u) +
      U(\mathcal{L}_{[<0]})\mathcal{N}(u),
    \end{equation}
    where \(U(\mathcal{L}_{[<0]})\mathcal{N}(u)\) is the subspace generated by
    \(\mathcal{N}(u)\) under the action of negative Virsoro modes
    \(L_{[-n]},n\in \NN\).
  \item We say that \(u\) is \emph{torus primary} if
    \begin{equation}
      \vvptf_\lambda (u,\tau)\neq0,\qquad\text{and}\qquad \vvptf_\lambda (w,\tau)=0
    \end{equation}
    for all \(w\in \mathcal{N}(u)\).
  \item Set
 \begin{equation}
 \mathcal{V}^u \left(\rho_\lambda\right) = \spn{ \vvptf_\lambda (w,\tau) \mid
   w\in \operatorname{Vir}(u)}\subset\mathcal{V}
 \left(\rho_\lambda\right),\qquad
 \mathcal{V}^u \left(\rho_\lambda\right)_n =\mathcal{V}^u
 \left(\rho_\lambda\right) \cap \mathcal{V} \left(\rho_\lambda\right)_n.
 \end{equation}
That is, $\mathcal{V}^u (\rho_\lambda)$ is the subspace of
$\mathcal{V}(\rho_\lambda)$ consisting of all evaluations of $\vvptf_\lambda
(-,\tau)$ on \correct{vectors in  \(\operatorname{Vir}(u)\)}.
  \end{enumerate}
\end{defn}
Note that $\operatorname{Vir}(u)$ is a Virasoro module, that is, closed under the 
 action of the Virasoro algebra, if and only if $u$ is a singular vector (an
 $L_{[0]}$-eigenvector that is annihilated by all positive Virasoro
 modes). In this case \(\operatorname{Vir}(u)=U(\mathcal{L})u\) and \(\mathcal{N}(u)=0\).}

\correct{
\begin{lem}
  Let \(u\in W^\lambda\) be a torus primary vector.
  \begin{enumerate}
  \item The subspace \(U(\mathcal{L}_{[<0]})\mathcal{N}(u)\) lies in the kernel of
    \(\vvptf_\lambda(-,\tau)\), that is, for all \(m\in
    U(\mathcal{L}_{[<0]})\mathcal{N}(u)\), 
    \begin{equation}
      \vvptf_\lambda(m,\tau)=0.
    \end{equation}
    \label{itm:descsub1}
  \item For every \(w\in U(\mathcal{L})u\) there
    exists a vector \(\tilde{w}\in \operatorname{Vir}(u)\) such that
    \begin{equation}
      \vvptf_\lambda(w,\tau)=\vvptf_\lambda(\tilde{w},\tau).
    \end{equation}
    In particular, if \(w\) lies in an \(L_{[0]}\)-eigenspace, then
    \(\tilde{w}\) can be chosen to lie in the same eigenspace.
    \label{itm:descsub2}
  \end{enumerate}
  \label{thm:descsubspace}
\end{lem}}
\begin{proof}
  \correct{
  \cref{itm:descsub2} is clearly implied by \cref{itm:descsub1}, so we need only
  show \cref{itm:descsub1}.}

\correct{
  The subspace \(U(\mathcal{L}_{[<0]})\mathcal{N}(u)\) is spanned by elements
  of the form
  \begin{equation}
    L_{[-n_1]}\cdots L_{[-n_\ell]}y, \qquad y\in \mathcal{N}(u),\ n_i=1,2,
    \label{eq:genkerel}
  \end{equation}
  since \(L_{[-1]}\) and \(L_{[-2]}\) generate all other negative Virasoro
  modes. We prove the lemma by induction in the word length \(\ell\).}

\correct{
  For the base case with \(\ell=0\) we have that \(\wt[y]<\wt[u]\) for all
  homogeneous vectors \(y\in \mathcal{N}(u)\) and hence
  \begin{equation}
    \vvptf_\lambda(y,\tau)=0.
  \end{equation}
  For the induction step assume that \cref{itm:descsub1} holds on all elements
  of the form \eqref{eq:genkerel} for all \(\ell \le P\in \NN_0\) and consider
  an element of the form
  \begin{equation}
    w= L_{[-n_1]} x, \quad x=L_{[-n_2]}\cdots L_{[-n_{P+1}]}y,\quad y\in \mathcal{N}(u).
  \end{equation}
  If \(n_1=1\), then by \eqref{eq:1pt[0]result} (recall
  $\correct{\tilde{\omega}}_{[n+1]}=L_{[n]}$) it follows that $\vvptf_\lambda (w,\tau)=0$. If
  instead \(n_1=2\), \eqref{eq:1pt[-1]result} implies
  \begin{equation}
    \vvptf_\lambda (w,\tau) = \partial \vvptf_\lambda (x,\tau) +\sum_{j=2}^\infty \eis{2j}(\tau) \vvptf_\lambda (L_{[2j-2]}x,\tau).
  \end{equation}
  Since \(x\) satisfies the induction hypothesis \(\vvptf_\lambda
  (x,\tau)=0\), while the arguments of the \(\vvptf_\lambda(-,\tau)\) in the second summand can be rewritten as
  follows.
  \begin{equation}
    L_{[2j-2]}x=[L_{[2j-2]},L_{[-n_2]}\cdots L_{[-n_{P+1}]}]y+L_{[-n_2]}\cdots L_{[-n_{P+1}]}L_{[2j-2]}y.
  \end{equation}
  Thus the second summand above satisfies the induction hypothesis
  because \(L_{[2j-2]}y\in \mathcal{N}(u)\). Further, since \(j\ge2\)
  evaluating the commutator in the first summand will give a linear
  combination of vectors satisfying the induction hypothesis, hence \(\vvptf_\lambda (w,\tau)=0\).}
\end{proof}

Recall the ring of modular differential operators given in \eqref{eq:moddiffops}.

\begin{prop}
  Let \(u\in  W^\lambda_{[\wt[u]]}\) be a
  \correct{torus primary vector}.
    Then $\mathcal{V}^{u} \left(\rho_\lambda \right) \subseteq
   \mathcal{R}\vvptf_\lambda (u,\tau)$. Furthermore, if
   $-\operatorname{wt}[u] \notin \NN_0$
   then $\mathcal{V}^{u} \left(\rho_\lambda \right) 
   = \mathcal{R}\vvptf_\lambda (u,\tau)$, that is, \(\mathcal{V}^{u} \left(\rho_\lambda \right)\) is cyclic
    as an \(\mathcal{R}\)-module and is generated by \(\vvptf_\lambda (u,\tau)\).
    \label{prop:RmoduleResult}
  \end{prop}

  \begin{proof}
    This proof is a
    generalisation of \cite[\refProp 2(b)]{DM-MoonshineHigher}, \correct{where
      \(u\) being a Virsoro singular vector is replaced by \(u\) being
      torus primary}.
 
 We first prove $\mathcal{V}^{u} \left(\rho_\lambda \right) \subseteq \mathcal{R}\vvptf_\lambda (u,\tau)$.
  Consider $\vvptf_\lambda (w,\tau) \in \mathcal{V}^{u} \left(\rho_\lambda
  \right)$\correct{, with the weight of \(w\) being \(\wt[w]=h_{\lambda}+N\)} 
  for some $N \in \NN_0$. 
  We will prove by
  induction on $N$ that $\vvptf_\lambda (w,\tau) \in
  \mathcal{R}\vvptf_\lambda (u,\tau)$.
  \correct{The base case is when \(N=\wt[u]-h_\lambda\) (as \(\operatorname{Vir}(u)\) has
  no weight spaces of lesser conformal weight) and hence \(w\) is a scalar
  multiple of \(u\), then \(\vvptf_\lambda (w,\tau) \in
  \mathcal{R}\vvptf_\lambda (u,\tau)\) by assumption.}

\correct{Suppose the result holds up to some arbitrary $N \ge \wt[u]-h_\lambda\in \NN_0$ and consider the
  case $w\in \operatorname{Vir}(u)$ with \(\wt[w]= h_\lambda + N+1\).} Since 
   $w$ is not a scalar multiple of $u$, we may
  assume without loss of generality that $w=L_{[-n_1]}L_{[-n_2]}\cdots
  L_{[-n_t]}u$, where $n_j=1,2$ for all $1\leq j\leq t$ since
  $L_{[-1]}$ and $L_{[-2]}$ generate $L_{[-n]}$ for all $n>0$.
  In the case $n_1=1$, we have by \eqref{eq:1pt[0]result}
  that $\vvptf_\lambda (w,\tau)=0$. If $n_1=2$, then setting $x=L_{[-n_2]}\cdots L_{[-n_t]}u$ and using \eqref{eq:1pt[-1]result} we have
  \begin{equation}
  \vvptf_\lambda (w,\tau) = \partial \vvptf_\lambda (x,\tau) +\sum_{j=2}^\infty \eis{2j}(\tau) \vvptf_\lambda (L_{[2j-2]}x,\tau).
\end{equation}
\correct{
 Note $\wt [x]$ and $\wt [L_{[2j-2]}x]$ are both strictly less than $\wt
[w]=h_\lambda +N+1$. 
 Since \(x\in \operatorname{Vir}(u)\) this implies that the
induction hypothesis applies to \(x\) and hence \(\partial \vvptf_\lambda
(x,\tau)\in \mathcal{R}\vvptf_\lambda (u,\tau)\). The \(L_{[2j-2]}x\) need
not lie in \(\operatorname{Vir}(u)\), however, by \cref{thm:descsubspace}
there exist vectors \(x_j\in \operatorname{Vir}(u)\) of the same weight such
that \(\vvptf_\lambda (L_{[2j-2]}x,\tau)=\vvptf_\lambda (x_j,\tau)\) to which
the induction hypothesis therefore applies, hence \(\vvptf_\lambda
(L_{[2j-2]}x,\tau)=\vvptf_\lambda (x_j,\tau)\in \mathcal{R}\vvptf_\lambda (u,\tau)\).
It follows that $\vvptf_\lambda (w,\tau)\in \mathcal{R}\vvptf_\lambda (u,\tau)$.}

  We turn to showing that $\mathcal{R}\vvptf_\lambda (u,\tau)\subseteq \mathcal{V}^{u} \left(\rho_\lambda \right)$ 
  if $-\operatorname{wt}[u] \notin \NN_0$.
    Recall every element in $\mathcal{R}\vvptf_\lambda (u,\tau)$ is of the form
  \begin{equation}
  \left(\phi_0+\phi_1\partial +\cdots +\phi_t \partial^t\right)\vvptf_\lambda (u,\tau),
  \end{equation}
  for $\phi_{\ell} \in \mathcal{M}$ with $0\leq \ell \leq t$ and $t\in \mathbb{N}_0$. We also have
  that $\phi_{\ell}$ is a linear combination of terms $\eis{4}(\tau)^i
  \eis{6}(\tau)^j$ for some $i,j\in \NN_0$.
  Thus, by linearity, we need only show that \(\mathcal{V}^{u}
  \left(\rho_\lambda \right)\) is closed under taking modular derivatives
  \(\partial\), and multiplication by \(\eis{4}(\tau)\) and \(\eis{6}(\tau)\) 
  if $-\operatorname{wt}[u] \notin \NN_0$.
 For $z\in \operatorname{Vir}(u)$ and $r\in \mathbb{N}_0$, set $x_r(z) = L_{[-2]}L_{[-1]}^{2r}z$. 
 Then for $r\geq 1$ and using \eqref{eq:1pt[-1]result}
  \begin{align}
  \vvptf_\lambda \left(x_r(z),\tau\right) &=\partial\vvptf_\lambda\brac*{L_{[-1]}^{2r}z,\tau}+\sum_{j=2}^\infty \eis{2j}(\tau)\vvptf_\lambda\brac*{L_{[2j-2]}L_{[-1]}^{2r}z,\tau}\nonumber\\
  &= \alpha
  \eis{2r+2}(\tau)\vvptf_\lambda (z,\tau) + \sum_{\ell =r+2}^\infty
  \eis{2\ell}(\tau) \vvptf_\lambda \left(L_{[2\ell-2]}L_{[-1]}^{2r}z,\tau \right) 
  \end{align}
  \correct{where $\alpha=(2r+1)!\wt[z]$
  is \(0\) if and only if \(\operatorname{wt}[z] =  0\),
  and $\wt\left[L_{[2\ell-2]}L_{[-1]}^{2r}z\right] 
  < \wt [z]$ for $\ell \geq r+2$. In the second equality above the derivative
  term vanishes due to  \eqref{eq:1pt[0]result}. Similarly the terms in the
  sum also vanish by \eqref{eq:1pt[0]result} for \(2\le \ell \le r\) by using the
  Virasoro relations to show that \(L_{[2\ell-2]}L_{[-1]}^{2r}z\) lies in the
  image of \(L_{[-1]}\). The same calculation also shows
  \(L_{[2r]}L_{[-1]}^{2r}z=\alpha z+L_{[-1]}y=(2r+1)!\wt[z]z+L_{[-1]} y
  \), for some vector \(y\).
  By \cref{thm:descsubspace} 
  }
   there exist \(x_{\ell,r,z}\in
  \operatorname{Vir}(u)\) in the same weight space as
  \(L_{[2\ell-2]}L_{[-1]}^{2r}w\) such that
  \begin{equation}
    \vvptf_\lambda \left(x_r(z),\tau\right) =
    \alpha
  \eis{2r+2}(\tau)\vvptf_\lambda (z,\tau) + \sum_{\ell =r+2}^\infty
  \eis{2\ell}(\tau) \vvptf_\lambda \left(x_{\ell, r,z},\tau \right).
   \label{eq:DMproof01}
  \end{equation}
  
 We first claim that for any $r \geq 1$ and $\vvptf_\lambda (w,\tau)\in
 \mathcal{V}^{u}(\rho_\lambda )$ we have $\eis{2r+2}(\tau) \vvptf_\lambda
 (w,\tau)\in \mathcal{V}^{u}(\rho_\lambda )$. We prove this by induction on
 $\wt[w]$, and more specifically, by induction on the non-negative integer $N$
 such that $\wt [w]=\wt[u]+N$. In the case $N =0$, we have that $w$ is a
 scalar multiple of $u$ and hence we can choose \(w=u\) without loss of
 generality. Then \eqref{eq:DMproof01} gives
 \begin{equation}
 \vvptf_\lambda \left(x_r(u),\tau\right) = \alpha \eis{2r+2}(\tau)\vvptf_\lambda (u,\tau) + \sum_{\ell =r+2}^\infty \eis{2\ell}(\tau) \vvptf_\lambda \left(x_{\ell, r,u},\tau \right) = \alpha \eis{2r+2}(\tau)\vvptf_\lambda (u,\tau)
 \end{equation}
 since $\wt [x_{\ell ,r,w}]<\wt [u]$ so that $\vvptf_\lambda (x_{\ell,
   r,z},\tau ) = 0$ for $\ell \geq r+2$ by assumption. Here $\alpha$ is
   non-zero by our assumption $-\operatorname{wt}[u]\notin\NN_0$.
 Thus $\eis{2r+2}(\tau) \vvptf_\lambda (w,\tau)\in \mathcal{V}^{u} (\rho_\lambda
 )$ for $r\geq 1$ and the base case is established. Assume $\eis{2r+2}(\tau)
 \vvptf_\lambda (w,\tau)\in \mathcal{V}^{u} \left(\rho_\lambda \right)$  for
 any $r \geq 1$ when $\wt[w]= \wt[u]+m$ and
 $0\leq  m\leq N$ \correct{for some \(N\in\mathbb{N}_0\)}.
 Consider the $N+1$ case, i.e., $\wt [w]=\wt[u]+N+1$. Then \eqref{eq:DMproof01} becomes
 \begin{equation}
 \vvptf_\lambda \left(x_r(w),\tau\right) = \alpha
 \eis{2r+2}(\tau)\vvptf_\lambda (w,\tau) + \sum_{\ell =r+2}^\infty
 \eis{2\ell}(\tau) \vvptf_\lambda \left(x_{\ell, r,w},\tau \right)
 \end{equation}
 where $\wt [x_{\ell,r,w}] < \wt [w]=\wt[u]+N+1$ for $\ell \geq r+2$.
 Again, $\alpha$ is non-zero since $\operatorname{wt}[w]>\operatorname{wt}[u]$
 and \(-\operatorname{wt}[u]\notin \NN_0\), so \(\operatorname{wt}[w]\neq0\). As such, we have by our induction hypothesis that $\eis{2\ell}(\tau) \vvptf_\lambda (x_{\ell, r,w},\tau ) \in \mathcal{V}^{u} (\rho_\lambda )$ for $\ell \geq r+2$. It follows that $\eis{2r+2}(\tau)\vvptf_\lambda (w,\tau) \in \mathcal{V}^{u} (\rho_\lambda )$ for $r\geq 1$.
 
 Next, we establish that $\partial \vvptf_\lambda (w,\tau ) \in \mathcal{V}^{u} (\rho_\lambda )$ for any $w\in \operatorname{Vir}(u)$. Indeed, \eqref{eq:1pt[-1]result} gives
 \begin{equation}
  \vvptf_\lambda \left(x_0(w),\tau\right) = \partial \vvptf_\lambda (w,\tau) +
  \sum_{\ell =2}^\infty \eis{2\ell}(\tau) \vvptf_\lambda
  \left(\correct{L_{[2\ell-2]}w}, \tau \right)
 \end{equation}
 where $\wt [\correct{L_{[2\ell-2]}w}] < \wt [w]$. Since \correct{$L_{[2\ell-2]}w$ can be replaced by $x_{\ell,w} \in \operatorname{Vir}(u)$ of the same weight by Lemma \ref{thm:descsubspace}} and
 $\eis{2\ell}(\tau) \vvptf_\lambda (\correct{x_{\ell,w}},\tau ) \in \mathcal{V}^{u}
 (\rho_\lambda )$ for $\ell \geq 2$
 we have $\partial \vvptf_\lambda (w,\tau) \in \mathcal{V}^{u}
 (\rho_\lambda )$.
 Thus \(\mathcal{V}^{u} (\rho_\lambda )\) is closed under the action of
  \(\mathcal{R}\) and hence \(\mathcal{V}^{u} (\rho_\lambda )=\mathcal{R}\vvptf_\lambda(u,\tau)\).
\end{proof}

Since the above proposition gives the decomposition of
\(\mathcal{V}^u(\rho_\lambda)\) as an \(\mathcal{R}\)-module it is natural to
ask how it decomposes as an \(\mathcal{M}\)-module. To this end we require the
following \correct{proposition}.

\begin{prop}[{Marks-Mason\cite[\refThm 1]{MarksMason10}, Marks \cite[\refLem 2.7]{Marks2011}}]
  Let \(F\in\mathcal{H}\brac*{k,\rho,\nu}\) be a $d$-dimensional holomorphic \vvmf{} of weight
  \(k\) for a representation and multiplier system \((\rho,\nu)\), whose components form a fundamental set of solutions for a  monic modular differential equation. Then the set \(\set{F,\pd F,\dots \pd^{d-1}F}\) is
  an \(\mathcal{M}\)-basis for \(\mathcal{R}F\).
  Further, if
  \(c(d,n)\in \NN\) are the coefficients of the series expansion
  \begin{equation}
    \frac{1-t^{2d}}{\brac*{1-t^2}\brac*{1-t^4}\brac*{1-t^6}}=\sum_{n=0}^\infty
    c(d,n)t^n,
    \label{eq:HPseries}
  \end{equation}
  then the weight \(k+n\) homogeneous subspace \((\mathcal{R}F)_n = \{ (\phi_0
  +\phi_1\partial + \phi_2 \partial^2 + \cdots )F \mid \phi_j \in
  \mathcal{M}_{n-2j}, j\geq 0\}\) satisfies
  \begin{equation}
    \dim \left((\mathcal{R}F)_n\right) = c(d,n).
    \label{eq:gradeddim}
  \end{equation}
  \label{thm:linindep}
\end{prop}
\begin{proof}
 In \cite[\refThm 1]{MarksMason10} it was shown that
  \correct{\(\mathcal{H}\brac*{\rho,\nu}\)} is free and rank \(d\) over
  \(\mathcal{M}\). Further, in \cite[\refLem 2.7]{Marks2011} it is shown that \(\set{F,\pd F,\dots ,\pd^{d-1}F}\) is
   linearly independent over \(\mathcal{M}\). By assumption, the components of \(F\)
  form a fundamental set of solutions for a monic modular differential
  equation of degree \(d\), hence \(\pd^d F\) lies in the \(\mathcal{M}\)-span of \(\set{F,\pd F,\dots ,
    \pd^{d-1}F}\). Thus \(\set{F,\pd F,\dots ,
    \pd^{d-1}F}\) is an \(\mathcal{M}\)-basis for
  \(\mathcal{R}F\).
  The series \eqref{eq:HPseries} is the Hilbert-Poincar\'e series for
  \(\mathcal{R}F\), see  \cite[\refEq 2.9]{Marks2011}.
\end{proof}

 \begin{thm}
   \correct{Let \(u\in  W^\lambda_{[\wt[u]]}\) be a torus primary vector
     satisfying \(-\wt[u]\notin \NN_0\)}.
    Let \(\set{\mu_1,\dots,\mu_{\delta(\lambda)}}\) (recall $\delta (\lambda)=\dim (\mathcal{C}_1(W^\lambda))$) be the leading
    exponents of \(\vvptf_\lambda (u,\tau)\) and let
    \(\mu_{\min{}},\mu_{\max{}}\) be the least and greatest leading exponent,
    respectively. Suppose further that the following four conditions hold.
    \begin{enumerate}
    \item The subspace \(N=\set{w\in W^\lambda_{[\wt[u]]}\st \vvptf_\lambda
          (w,\tau)=0}\) has codimension \(1\) in \(W^\lambda_{[\wt[u]]}\).
      \label{itm:uniqvec}
    \item The exponents  \(\set{\mu_1,\dots,\mu_{\delta(\lambda)}}\) saturate the
      inequality \eqref{eq:bigweightineq}, that is, it is an equality.
      \label{itm:minimalwtass}
    \item The exponents \(\set{\mu_1,\dots,\mu_{\delta(\lambda)}}\) of
      \(\vvptf_\lambda (u,\tau)\) are minimal among all \vvmf{s} in
      \(\mathcal{V}(\rho_\lambda)\), that is, for any \(F\in
      \mathcal{V}(\rho_\lambda)\) the leading exponent of the \(j\)th
      component will be at least \(\mu_{j}\).
      \label{itm:minexp}
     \item The space \(\mathcal{H}(\rho_\lambda,\nu_{h_\lambda-12\mu_{min}})\) is
       cyclic over \(\mathcal{R}\).
      \label{itm:cyclicass}
    \end{enumerate}
    Then
    \begin{equation}
      \mathcal{V}^u(\rho_\lambda)\subset\mathcal{V}(\rho_\lambda)\subset\eta^{24\mu_{\min{}}}\mathcal{H}\brac*{\rho_\lambda,\nu_{h_\lambda-12\mu_{\min{}}}}
      \label{eq:incseq}
    \end{equation}
    and \(\set{\mu_1-\mu_{\min{}},\dots,\mu_{\delta(\lambda)}-\mu_{\min{}}}\) is an admissible set for $(\rho_\lambda , \nu_{h_\lambda-12\mu_{\min{}}})$.
    Moreover, let
     \(\set{\lambda_1,\dots,\lambda_{\delta(\lambda)}}\) be the minimal
    admissible set for $(\rho_\lambda , \nu_{h_\lambda-12\mu_{\min{}}})$ 
    and define
    \begin{equation}
      M =\sum_{i=1}^{\delta(\lambda)} \brac*{\mu_{i}-\mu_{\min{}}-\lambda_i}.
    \end{equation}
    Then
    \begin{equation}
      c(\delta(\lambda),n)=\dim\brac*{\mathcal{V}^u(\rho_\lambda)_{n\correct{+\wt[u]-h_\lambda}}}\le
      \dim\brac*{\mathcal{V}(\rho_\lambda)_{n\correct{+\wt[u]-h_\lambda}}} \le
      c\brac*{\delta(\lambda),n+12\frac{M}{\delta(\lambda)}},
      \label{eq:dimineq}
    \end{equation}
    where \(c(\delta(\lambda),n)\) are the graded dimensions of
    \eqref{eq:gradeddim}. 
    If
    additionally 
    \begin{enumerate}[resume]
        \item \(\mu_{\max{}}-\mu_{\min{}}<1\),
      \label{itm:expdiffass}
    \end{enumerate}
    then \correct{this is a necessary and sufficient condition for}
    \(\set{\mu_1-\mu_{\min{}},\dots,\mu_{\delta(\lambda)}-\mu_{\min{}}}\)
    \correct{to be} 
    the minimal admissible set for \((\rho_\lambda
      ,\nu_{h_\lambda-12\mu_{\min{}}})\).
      \correct{Further,}
    \begin{equation}
      \mathcal{V}^u(\rho_\lambda)=\mathcal{V}(\rho_\lambda)=\eta^{24\mu_{\min{}}}\mathcal{H}\brac*{\rho_\lambda,\nu_{h_\lambda-12\mu_{\min{}}}}
      \label{eq:holomchar}
    \end{equation}
    and
    \begin{equation}
      \dim\brac*{\mathcal{V}(\rho_\lambda)_{n\correct{+\wt[u]-h_\lambda}}} = c(\delta(\lambda),n).
      \label{eq:strongdimform}
    \end{equation}
    \label{thm:generaldimformula}
  \end{thm}

  \begin{proof}
    \correct{We begin by showing the sequence of inclusions \eqref{eq:incseq},
      where the first inclusion holds by definition, so only the second
      need be justified.}
    By \cref{prop:RmoduleResult}, \(\mathcal{V}^u(\rho_\lambda) \subset
    \mathcal{R}\vvptf_\lambda (u,\tau)\) and
    \(\mathcal{V}^u(\rho_\lambda)=\mathcal{R}\vvptf_\lambda (u,\tau)\) if
    $-\operatorname{wt}[u]\notin \NN_0$.
    Since the exponents of \(\vvptf_\lambda (u,\tau)\) saturate the
    inequality \eqref{eq:bigweightineq}, so do the exponents of
    \(\eta^{-24\mu_{\min{}}}\vvptf_\lambda (u,\tau)\) and the exponents of each component function are
    all non-negative, hence \(\eta^{-24\mu_{\min{}}}\vvptf_\lambda
    (u,\tau)\) is holomorphic.  
    \correct{A brief computation shows \(\eta\) is in the kernel of the modular
      derivative and hence, \(\eta\) commutes with
    modular derivatives. We therefore have}
    that \(\mathcal{R}\eta^{-24\mu_{\min{}}}\vvptf_\lambda
    (u,\tau)=\eta^{-24\mu_{\min{}}}\mathcal{R}\vvptf_\lambda
    (u,\tau)=\eta^{-24\mu_{\min{}}}\mathcal{V}^u(\rho_\lambda)\). 
    Additionally, the leading exponents of \(\vvptf_\lambda (u,\tau)\) are minimal among all
    \vvmf{s} in \(\mathcal{V}(\rho_\lambda)\), \correct{and thus}
    \(\eta^{-24\mu_{\min{}}}\mathcal{V}(\rho_\lambda)\subset\mathcal{H}(\rho_\lambda,\nu_{h_\lambda-12\mu_{\min{}}})\) \correct{(this inclusion uses the fact
    that \(\nu_{h_\lambda}\nu_{-12\mu_{\min{}}}=\nu_{h_\lambda-12\mu_{\min{}}}\))}
    and \eqref{eq:incseq} follows. 
    The bounds \correct{\eqref{eq:dimineq} on the graded dimensions} then
    follow from the fact that \(\mathcal{V}^u(\rho_\lambda)\) and
    \(\mathcal{H}(\rho_\lambda,\nu_{h_\lambda-12\mu_{\min{}}})\) are
    both cyclic \(\mathcal{R}\)-modules and the weight of the cyclic
    generator of
    \(\mathcal{H}(\rho_\lambda,\nu_{h_\lambda-12\mu_{\min{}}})\) \correct{(computed
    from \eqref{eq:weightineq})}
    differs from the weight of
    \(\correct{\eta^{-24\mu_{\min{}}}}\vvptf_\lambda  (u,\tau)\)
    by \(12\frac{M}{\delta(\lambda)}\).

    \correct{The set \(\set{\mu_1-\mu_{\min{}},\dots,\mu_{\delta(\lambda)}-\mu_{\min{}}}\) is
    minimal admissible for \((\rho_\lambda ,\nu_{h_\lambda-12\mu_{\min{}}})\)
    if and only if \(\mu_{\max{}}-\mu_{\min{}}<1\), because its members are
    non-negative by construction and bounded by \(1\) precisely when all
    \(\mu_i\) differ by less than \(1\).}
    Thus \(\eta^{-24\mu_{\min{}}}\vvptf_\lambda
    (u,\tau)\) has the same weight as the cyclic generator of
    \(\mathcal{H}(\rho_\lambda,\nu_{h_\lambda-12\mu_{\min{}}})\), which
    lies in a \(1\)-dimensional weight space, that is, \(\eta^{-24\mu_{\min{}}}\vvptf_\lambda
    (u,\tau)\) is a non-zero scalar multiple of the cyclic generator. Hence
    \(\eta^{-24\mu_{\min{}}}\mathcal{V}^u(\rho_\lambda)=\mathcal{H}(\rho_\lambda,\nu_{h_\lambda-12\mu_{\min{}}})\),
    which implies \eqref{eq:holomchar}\correct{. By \cite[Theorem
      1.3]{MarksMason10} it is known that the 
    components of the cyclic generator of \(\mathcal{H}(\rho_\lambda ,\nu_{h_\lambda -12\mu_{\operatorname{min}}}) \)
    form a fundamental system of solutions of a  monic modular differential
    equation of order $\delta (\lambda)$. 
    Proposition \ref{thm:linindep} now gives the dimension formula \eqref{eq:strongdimform}.
}
  \end{proof}

\cref{prop:RmoduleResult} and \cref{thm:generaldimformula} generalise Lemma
2.1, Theorem 3.5, and Corollary 3.3 of \cite{intvvmf18}. Additionally, we
note here that the necessary condition on the values of $\operatorname{wt}[u]$
is absent in the statement of Lemma 2.1 in \cite{intvvmf18}, and thus also
Theorem 3.5 and Corollary 3.3 in loc.\ cit. Indeed, in the notation of that
paper, these results and the relevant discussion should all include the
assumption that $-h_{m,n} \not \in \mathbb{N}_0$. Fortunately,
in the
application of these results to the Virasoro minimal models in
\cite[Section 3]{intvvmf18}, the analysis of small dimensions automatically excludes
$-h_{m,n}\in \mathbb{N}_0$ with one exception. This
  exception is in
the $1$-dimensional setting where the
trivial case of $(m,n)=(1,1)$ is included in the second statement of Theorem
3.7 when it should not be. However, this corresponds to the Virasoro
minimal model at central charge \(c=0\)
 (i.e., the trivial \voa{} isomorphic to \(\CC\)) acting on itself. In this case (up to rescaling) there is only one torus
\(1\)-point function and it is constant.

\section{Affine \texorpdfstring{\(\SLA{sl}{2}\)}{sl2}}
\label{sec:sl2}

In this section we apply the theory of the previous section to the simple affine vertex operator algebra at non-negative integral level \(k\) associated to the Lie algebra $\SLA{sl}{2}$,  which we denote $L(k,0)$. The central charge of this vertex operator algebra is
\begin{equation}
  \mathbf{c}=\frac{3k}{k+2}.
  \label{eq:ccharge}
\end{equation}
Denote the standard \(\SLA{sl}{2}\) Chevalley basis by $\{e,h,f\}$ subject
  to the standard commutation relations
\begin{equation}
  \comm{e}{f}=h,\quad \comm{h}{e}=2e,\quad \comm{h}{f}=-2f,
\end{equation}
with \(h\) spanning a choice of Cartan subalgebra $\mathfrak h \cong \mathbb Ch$. Further we identify the root lattice with even integers, \(Q\cong 2\ZZ\), and the weight lattice with the integers, \(P\cong \ZZ\). We normalise the Killing form such that the length squared of the highest root
is 2, that is, the non-vanishing pairings are \(\kappa (e,f)=\correct{\kappa(f,e)=}1\), \(\kappa
(h,h)=2\). As is well known, \(L(k,0)\) and its representation theory satisfy the assumptions
   given in the paragraph preceding \eqref{eq:iopspace}. 
In particular, the representation theory of \(L(k,0)\)
is semisimple and a complete set of representatives of simple modules is given
by the simple quotients of affine Verma modules \(L(k,\mu)\), where \(0\le
\mu\le k\) is the \(h_{0}\) eigenvalue of the
highest weight vector (also known as the finite part of the affine weight) \cite{Frenkel1992}.
The conformal weight of this highest weight vector is
\begin{equation}
  h_\mu=\frac{\mu(\mu+2)}{4(k+2)}.
  \label{eq:cwt}
\end{equation}

In order to obtain non-vanishing trace functions $\ptf^{\mathcal{Y}}(u,\tau)$, we 
need to find non-vanishing spaces of intertwining operators 
of type ${L(k,\mu) \choose L(k,\lambda) \, \, \, L(k,\mu)}$, 
and a suitable basis $\intset_\lambda$ as defined in \eqref{eq:intset}. 
Note that as intertwining operator spaces for triples of simple \(L(k,0)\)-modules are
always at most $1$-dimensional, we will identify the basis vectors in
\(\intset_\lambda\) with the index $\mu$ appearing in
${L(k,\mu) \choose L(k,\lambda) \, \, \, L(k,\mu)}$.
Further, let $\ptf^\mu$ denote the trace
of the basis intertwining operator
corresponding to the label \(\mu\) over the module $L(k,\mu)$, as in \eqref{eq:q-expansionOfTraceFunction}.
\begin{prop}
  Let \(0\le \lambda \le k\), then
\begin{equation}
  \intset_\lambda =
  \begin{cases}
    \left\{\mu \bigg\rvert \frac{\lambda}{2} \leq \mu \leq k - \frac{\lambda}{2}\right\}& \lambda \text{ even},\\
    \emptyset&\lambda \text{ odd}.
  \end{cases}
\end{equation}
In particular, if $\lambda$ is even then $\lvert \intset_\lambda \rvert = k -\lambda+1$.
\label{thm:formdim}
\end{prop}
\begin{proof}
  This is an immediate consequence of the \(L(k,0)\) fusion rules given by
  \(L(k,\lambda)\boxtimes L(k,\mu)=\bigoplus_\nu N_{\lambda\, \mu}^\nu L(k,\nu)\), where
  \begin{equation}
     N_{\lambda\, \mu}^\nu = 
         \begin{cases}
        1 & \text{if } |\lambda - \mu| \leq \nu \leq \mathrm{min}\{\lambda+\mu, 2k-\lambda-\mu\} \text{ and }  \lambda + \mu + \nu \equiv 0 \pmod{2},\\ 
        0 & \text{otherwise}.
         \end{cases}
         \label{eq:fusionrules}
    \end{equation}
  \end{proof}
    The fusion rules \eqref{eq:fusionrules} were originally presented in the
    physics literature 
    in \cite{ZamFat85,GepWit86} and predate \voa{s}. They were later proved in \cite[\refThm 1]{Tsuchiya1987} and \cite[\refCor 3.2.1]{Frenkel1992}.
    We next record many properties of the traces \(\ptf^{\mu}\), with the
remainder of the section dedicated to proving these properties.
Throughout the section for any $0\leq \lambda \leq k$ we denote the highest weight vector of \(L(k,\lambda)\) by \(\ket{\lambda}\). 

\begin{thm}
  Let \(0\le \lambda\le k\) with $\lambda$ even. 
  \begin{enumerate}
  \item Let \(0\le n \le \frac{\lambda}{2}-1\), then \(\ptf^{\mu}(v,\tau)=0\)
    for all \(v\in L(k,\lambda)_{[h_\lambda+n]}\) and $\mu\in
    \intset_\lambda$, where \(L(k,\lambda)_{[h_\lambda+n]}\) denotes a
      homogeneous space with respect to the square bracket grading \eqref{eq:ModulesqGrading}.
    \label{itm:zeros}
  \item For any $\mu \in \intset_\lambda$, the subspace \(N=\set{v \in
      L(k,\lambda)_{[h_\lambda+\frac{\lambda}{2}]}\st \ptf^{\mu}(v,\tau)=0}\)
    has codimension \(1\) in
    \(L(k,\lambda)_{[h_\lambda+\frac{\lambda}{2}]}\). Hence there is a unique
    \correct{torus primary} vector \(u\in L(k,\lambda)_{[h_\lambda+\frac{\lambda}{2}]}\), up to rescaling or addition by elements in \(N\), such that \(\ptf^{\mu}(u,\tau)\neq0\).
    The vector \(u\) can be chosen to be $f_{[-1]}^{\frac{\lambda}{2}}\ket{\lambda}$ or, more generally, as any element
    \begin{equation}
      u \in f_{-1}^{\frac{\lambda}{2}}\ket{\lambda}+\bigoplus_{n=0}^{\frac{\lambda}{2}-1}L(k,\lambda)_{h_\lambda+n},
      \label{eq:actingvec}
    \end{equation}
    where \(L(k,\lambda)_{h_\lambda+n}\) denotes a
      homogeneous space with respect to the standard conformal grading \eqref{eq:confgrading}.
    \label{itm:nonzeros}
  \item For any $\mu \in \intset_\lambda$ and \correct{torus primary}
    \(u\in L(k,\lambda)\) \correct{of the form} 
    \eqref{eq:actingvec}, the leading exponent of \(\ptf^{\mu}(u,\tau)\) is
    \begin{equation}
    \label{eqn:LeadingExp}
      h_\mu-\frac{\mathbf{c}}{24}=\frac{2\mu^2+4\mu -k}{8(k+2)}.
    \end{equation}
    These leading exponents saturate the inequality
      \eqref{eq:bigweightineq}, that is, it is an equality.
    \label{itm:leadingexp}
  \item For \correct{a torus primary} \(u\in L(k,\lambda)\) \correct{of
      the form} 
    \eqref{eq:actingvec}, the intertwining operator basis \(\intset_{\lambda}\) can be normalised such that all coefficients of the series expansion of \(\ptf^\mu(u,\tau)\) are rational for each \(\mu\in\intset_{\lambda}\).
    \label{itm:ratcoeff}
  \item Let  $u\in L(k,\lambda)$ be \correct{a torus primary}
    \correct{of the form} 
    \eqref{eq:actingvec}. The set
    \(\set{\ptf^\mu(u,\tau)\st \mu\in \intset_\lambda}\) is linearly
    independent in $\mathcal{C}_1^u(L(k,\lambda))$, the space of torus $1$-point
    functions evaluated at \(u\), and the set \(\set{\ptf^\mu(-,\tau)\st \mu\in
      \intset_\lambda}\) is linearly independent in
    $\mathcal{C}_1(L(k,\lambda))$, the space of (unevaluated) torus \(1\)-point functions. The dimension of both of these spaces is \(k-\lambda+1\), the cardinality of \(\intset_\lambda\).
    Thus the vector \(\vvptf_\lambda(u,\tau)=(\ptf^\mu(u,\tau)\st \mu\in
    \intset_\lambda)^t=(\ptf^{\frac{\lambda}{2}}(u,\tau),\dots,\ptf^{k-\frac{\lambda}{2}}(u,\tau))^t\)
    is a \(\lvert
    \intset_\lambda \rvert=k-\lambda+1\)-dimensional weakly holomorphic vector-valued modular form of  weight
    \(h_\lambda+\frac{\lambda}{2}\), representation \(\rho_\lambda\), and
      multiplier system \(\nu_{h_\lambda}\).
    Moreover,
  \begin{equation}
    \rho_\lambda(\modT)=\diag\left\{\mathbf{e}\left(r_{\frac{\lambda}{2}}\right), \dots ,\mathbf{e}\left(r_{k-\frac{\lambda}{2}}\right) \right\},
  \end{equation}
  where
  \begin{equation}
    r_\mu=h_\mu-\frac{\mathbf{c}}{24}-\frac{h_\lambda}{12},\qquad \mu\in \intset_\lambda.
  \end{equation}
    \label{itm:basis}
      \item For arbitrary $w\in
        L(k,\lambda)_{[\wt[w]]}$ we have that $\vvptf_\lambda (w,\tau)$ is a holomorphic
        vector-valued modular form if $\lambda \geq -2 +
        \sqrt{2k+4}$. Additionally, for \(u\in L(k,\lambda)\) \correct{a
          torus primary} \correct{of the form} 
        \eqref{eq:actingvec}, the following are equivalent.
       \begin{enumerate}
 \item $\vvptf_\lambda (u ,\tau) \in \mathcal{H}(\rho_\lambda
   ,\nu_{h_\lambda})$.
   \label{itm:hol}
 \item $\lambda \geq -2+\sqrt{2k+4}$.
   \label{itm:wtbound}
 \item $k\geq d-2+\sqrt{2d-1}$ or $k\leq d-2-\sqrt{2d-1}$, where $d=\lvert
   \intset_\lambda \rvert$ is the dimension of $\vvptf_\lambda (u ,\tau)$.
   \label{itm:lvlbound}
 \end{enumerate}
    \label{itm:holomcon}
\end{enumerate}
  \label{thm:sl2alldims}
\end{thm}

The proof of \cref{thm:sl2alldims}
requires some preparation. First note that the conformal weight $1$ space of \(L(k,0)\) is isomorphic to \(\SLA{sl}{2}\) and hence all conformal weight spaces of the simple modules \(L(k,\lambda)\) completely reduce into finite direct sums of finite-dimensional simple \(\SLA{sl}{2}\) modules. Thus, by
\cref{prop:trivialmodulerequiredfornonzerotrace}, for each $\mu \in
\intset_\lambda$ we have \(\ptf^{\mu}(-,\tau)\) vanishes when restricted to a
conformal weight space of \(L(k,\lambda)\) that does not contain the trivial
\(\SLA{sl}{2}\) module. Therefore, we need to find the conformal weight space in $L(k,\lambda)$ at which the trivial module first appears.

\begin{lem}
  \label{thm:trivialmodulelambdaovertwo}
  Let \(0\le \lambda \le k\), $\lambda$ even, and $n\in \mathbb{Z}$.
    \begin{enumerate}
    \item For \(\lambda\ge2\) and  \(n < \frac{\lambda}{2}\), the multiplicity of the trivial \(\SLA{sl}{2}\) module in the conformal weight spaces \(L(k,\lambda)_{h_\lambda+n}\) is \(0\).
      \label{itm:mult0}
    \item For \(\lambda\ge0\), the multiplicity of the trivial \(\SLA{sl}{2}\) module in the conformal weight spaces \(L(k,\lambda)_{h_\lambda+\frac{\lambda}{2}}\) is \(1\).
      \label{itm:mult1}
    \end{enumerate}
\end{lem}

\begin{proof}
If $n<0$, then $L(k,\lambda)_{h_\lambda +n} =\{0\}$.
 Meanwhile, for $n\geq 0$ and $\alpha \in \mathbb{Z}$ we define 
 \begin{equation}
 L(k,\lambda)_{h_\lambda +n\colon \alpha} = \left\{ w\in L(k,\lambda)_{h_\lambda +n} \mid h_0w=\alpha w  \right\}.
\end{equation}
We compute the multiplicity of the trivial \(\SLA{sl}{2}\) module by considering
character formulae obtained from the BGG resolution of
\(L(k,\lambda)\) in terms of Verma modules \(V(k,\sigma)\) (where \(\sigma\)
again denotes the \(\SLA{sl}{2}\) weight, that is, the \(h_0\) eigenvalue, of the
generating highest weight vector)
\begin{equation}
  \cdots \rightarrow V(k,\lambda_{i})\oplus V(k,\lambda_{-i})\rightarrow \cdots \rightarrow V(k,\lambda_{2})\oplus V(k,\lambda_{-2})\rightarrow V(k,\lambda_{1})\oplus V(k,\lambda_{-1})\rightarrow V(k,\lambda)\rightarrow L(k,\lambda)\rightarrow 0,
\end{equation}
where for \(j\in \ZZ\), 
\begin{equation}
  \lambda_{2j}=\lambda+2 \correct{j}(k+2),\quad \text{and}\quad \lambda_{2j-1}=-\lambda -2+2j(k+2).
  \label{eq:wtformulae}
\end{equation}
  The above resolution is given in \cite[\refSec 4]{BerFelsl290} using results
  from \cite{KacLieAlg85}.
  Recall that the character of a Verma module is given by
  \begin{equation}
    \mathrm{ch}[V(k,\mu)]=\mathrm{tr}_{V(k,\mu)}z^{h_0}q^{L_0-\frac{\mathbf{c}}{24}}=\frac{z^\mu q^{h_{\mu}-\frac{\mathbf{c}}{24}}}{\prod_{m=1}^\infty(1-z^2q^m)(1-q^m)(1-z^{-2}q^{m-1})}.
  \end{equation}
  The Verma character formulae, in turn, yield the character formula for simple modules via the BGG resolution above. We have
  \begin{align}
    q^{\frac{\mathbf{c}}{24}-h_\lambda}\mathrm{ch}[L(k,\lambda)]
    = q^{\frac{\mathbf{c}}{24}-h_\lambda} \mathrm{tr}_{L(k,\correct{\lambda})}z^{h_0}q^{L_0-\frac{\mathbf{c}}{24}}
    &=q^{\frac{\mathbf{c}}{24}-h_\lambda}\mathrm{ch}[V(k,\lambda)]
                                +q^{\frac{\mathbf{c}}{24}-h_\lambda}\sum_{i=1}^\infty (-1)^i\brac*{\mathrm{ch}[V(k,\lambda_i)]+\mathrm{ch}[V(k,\lambda_{-i})]}\nonumber\\
                              &= q^{\frac{\mathbf{c}}{24}-h_\lambda}\sum_{i=0}^\infty(-1)^i\brac*{\mathrm{ch}[V(k,\lambda_i)]-\mathrm{ch}[V(k,\lambda_{-i-1})]}\nonumber\\
    &=\sum_{i=0}^\infty(-1)^i
      q^{h_{\lambda_i}-h_\lambda}\frac{\sum_{n=0}^{\lambda_i}
      z^{\lambda_i-2n}}{\prod_{m=1}^\infty(1-z^2q^m)(1-q^m)(1-z^{-2}q^{m})},
      \label{eq:charform}
  \end{align}
  where multiplication by the factor \(q^{\frac{\mathbf{c}}{24}-h_\lambda}\)
  shifts the exponents of the above power series such that the coefficient of \(z^m q^n\) is the dimension of
  \(L(k,\lambda)_{h_\lambda +n\colon m}\). The last equality in the above
  character formula uses that \(\lambda_{-m}=-2-\lambda_{m-1}\), \(m\in \ZZ\), and hence
    the conformal weights corresponding to these \(\SLA{sl}{2}\) weights
    satisfy \(h_{\lambda_{-m}}=h_{\lambda_{m-1}}\). The expansion of
      \eqref{eq:charform} up to degree \(\frac{\lambda}{2}\) in \(q\) will
      allow us to conclude the lemma.

  Note that all even weight simple \(\SLA{sl}{2}\) modules have a
    \(1\)-dimensional weight \(0\) space. Further the weight \(2\) space
    vanishes for the trivial module, while it is \(1\)-dimensional for all
    other even weight simple \(\SLA{sl}{2}\) modules. Therefore
  the difference \(\dim (L(k,\lambda)_{h_\lambda+m:0}) - \dim
  (L(k,\lambda)_{h_\lambda+m:2})\) is the multiplicity of the trivial module
  in \(L(k,\lambda)_{h_\lambda+m}\). This difference is also equal to the difference of the coefficients of
  \(z^0 q^m\) and \(z^2q^m\) in the character formula \eqref{eq:charform}
  above. Further, note that $h_{\lambda_m }- h_{\lambda_0}$ increases monotonically in $m>0$, and in particular,
  \begin{equation}
    h_{\lambda_1}- h_{\lambda} = k-\lambda + 1, \quad h_{\lambda_2} - h_{\lambda} = k+\lambda + 3 > \frac{\lambda}{2}.
  \end{equation}
  Thus if we wish to
  expand $q^{\frac{\mathbf{c}}{24}-h_\lambda}\mathrm{ch}[L(k,\lambda)]$ up to degree $q^{\frac{\lambda}{2}}$ it
  is sufficient to only consider the summands coming from \(i=0,1\) in the
  character formula \eqref{eq:charform}. 
  To simplify formulae, we introduce the notation 
  \((q)_i=\prod_{m=1}^i(1-q^m), i\ge0\) and record the \(q\)-series
  identity \cite[\refEq 9.16]{Ridoutsl2} 
  \begin{equation}
    \frac{1}{\prod_{m\geq 1} (1-z^2 q^m)(1-z^{-2}q^m)} = \sum_{n \in \ZZ} z^{2n} \sum_{i \geq 0} \frac{q^{2i+|n|}}{(q)_i (q)_{i+|n|}}.
  \end{equation}
  This identity is a consequence of the identity
  \begin{equation}
    \frac{1}{\prod_{m \geq 1} (1-z^2q^m)}= \sum_{j=0}^\infty \frac{q^j}{(q)_j} z^{2j},
  \end{equation}
  in \cite[\refEq 2.2.5]{andrews84}.  Thus  
  \begin{align}
      q^{\frac{\mathbf{c}}{24}-h_\lambda}\mathrm{ch}[L(k,\lambda)] &= \frac{\sum_{n=0}^{\lambda} z^{\lambda-2n} - q^{k+1-\lambda} \sum_{n=0}^{2(k+1)-\lambda} z^{2(k+1)-\lambda-2n}}{\prod_{m\geq 1} (1-z^2q^m)(1-q^m)(1-z^{-2}q^m)} + \mathcal{O}(q^{\lambda + k +3})\nonumber\\
      &= \left[\sum_{n\in\ZZ} z^{2n} \sum_{i=0}^\infty \frac{q^{2i+|n|}}{(q)_i (q)_{i+|n|} (q)_ \infty} \right] \left[ \sum_{m=0}^\lambda z^{\lambda-2m} -q^{k-\lambda+1} \sum_{m=0}^{2(k+1)-\lambda} z^{2(k+1)-\lambda-2m}\right]+ \mathcal{O}(q^{\lambda + k +3}).
  \end{align}
  Collecting the $z^0$ terms gives
  \begin{equation}
    \sum_{m=0}^\lambda \sum_{i=0}^\infty \frac{q^{2i + |\frac{\lambda}{2}-m|}}{(q)_i (q)_{i + |\frac{\lambda}{2}-m|} (q)_\infty} - q^{k-\lambda +1} \sum_{m=0}^{2(k+1)-\lambda} \sum_{i=0}^\infty \frac{q^{2i + |k+1-\frac{\lambda}{2}-m|}}{(q)_i (q)_{i+|k+1-\frac{\lambda}{2} - m|} (q)_\infty}+\mathcal{O}(q^{\lambda + k +3}),
  \end{equation}
  while collecting the $z^2$ terms gives
  \begin{equation}
    \sum_{m=0}^\lambda \sum_{i=0}^\infty \frac{q^{2i + |\frac{\lambda}{2}-m-1|}}{(q)_i (q)_{i+|\frac{\lambda}{2}\correct{-m-1}|} (q)_\infty} - q^{k-\lambda+1} \sum_{m=0}^{2(k+1)-\lambda} \sum_{i=0}^\infty \frac{q^{2i + |k-\frac{\lambda}{2}-m|}}{(q)_i (q)_{i+|k-\frac{\lambda}{2}-m|} (q)_\infty}+\mathcal{O}(q^{\lambda + k +3}).
  \end{equation}
  The difference of the $z^0$ and $z^2$ terms is therefore
  \begin{align}
    &\sum_{m=0}^\infty \brac*{\dim \left(L(k,\lambda)_{h_\lambda+m:0}\right) - \dim \left(L(k,\lambda)_{h_\lambda+m:2}\right)}q^m\nonumber\\
     &\quad= \sum_{i=0}^\infty \left[\frac{q^{2i+\frac{\lambda}{2}}}{(q)_i (q)_{i+\frac{\lambda}{2}} (q)_\infty} - \frac{q^{2i+\frac{\lambda}{2}+1}}{(q)_i (q)_{i+\frac{\lambda}{2}+1} (q)_\infty} \right] - q^{k-\lambda+1} \sum_{i=0}^\infty \left[ \frac{q^{2i+k+1-\frac{\lambda}{2}}}{(q)_i (q)_{i+k+1-\frac{\lambda}{2}} (q)_\infty} - \frac{q^{i+k+2-\frac{\lambda}{2}}}{(q)_i (q)_{i + k +2 -\frac{\lambda}{2}} (q)_\infty}\right]+\mathcal{O}(q^{\lambda + k +3})\nonumber\\
      &\quad= q^{\frac{\lambda}{2}} + \mathcal{O}(q^{\frac{\lambda}{2}+1}),
  \end{align}
  where we have used that only the first term of the first summand at \(i=0\) contributes to \(q^{\frac{\lambda}{2}}\). Thus, both parts of the lemma follow.
\end{proof}

\begin{lem}
  \label{thm:leadingcoefficient}
  Let \(0\le\lambda\le k\), $\lambda$ even, and $u\in L(k,\lambda)$ be \correct{of the form} 
  \eqref{eq:actingvec}.
  \begin{enumerate}
  \item For all \(\mu\in \intset_\lambda\), the trace
    \(\ptf^\mu(u,\tau)\) is non-vanishing\correct{, hence \(u\) is torus primary,} and the leading
    exponent is
  \begin{equation}
    h_\mu-\frac{\mathbf{c}}{24}=\frac{2\mu^2+4\mu-k}{8(k+2)}.
    \label{eq:expformla}
  \end{equation}
  These leading exponents saturate the inequality
    \eqref{eq:bigweightineq}, that is, it is an equality.
  \label{itm:leadingterm}
\item The intertwining operator underlying the trace \(\ptf^\mu(u,\tau)\) can be normalised such that all coefficients of the series expansion are rational.
  \label{itm:scaletorat}
\end{enumerate}

\end{lem}

\begin{proof}
  Recall that we denote by \(o(u)\) the coefficient of \(z^{-\wt (u)}\) in the
  series expansion of \(u\) inserted into the intertwining operator underlying
  the trace \(\ptf^\mu\). Therefore,
  \begin{equation}
    \ptf^\mu(u,\tau) =\mathrm{tr}_{L(k,\mu)} o(u) q^{L_0-\frac{\mathbf{c}}{24}} = q^{h_\mu -\frac{\mathbf{c}}{24}} \sum_{n=0}^\infty q^n\mathrm{tr}_{L(k,\mu)_{h_\mu + n}} o(u).
  \end{equation}
  We compute the coefficient of the leading term corresponding to $n=0$ and show that it is non-zero, which in turn will imply that formula \eqref{eq:expformla} gives the leading exponent. Note that \(L(k,\mu)_{h_\mu}\) is a module over the finite-dimensional Lie algebra \(\SLA{sl}{2}\) by restriction, and it is isomorphic to the simple highest weight module \(L(\mu)\) of highest weight \(\mu\). We choose the basis \(\set{v_i = f_0^i\ket{\mu}}_{i=0}^\mu\) of
  \(L(k,\mu)_{h_\mu}\cong L(\mu)\) and the corresponding dual basis
  \(\set{\phi_i}_{i=0}^{\mu} \subset{L(\mu)^\ast\cong L(\mu)}\). Let \(\pair{\ }{\ }\) denote the standard
  pairing between \(L(\mu)\) and its dual space so that the standard left
  action of \(\SLA{sl}{2}\) on the dual space is characterised by \(\pair{x_0
    \psi}{w}=-\pair{\psi}{x_0w}\), \(x\in\SLA{sl}{2}\), \(\psi\in
  L(\mu)^\ast\), \(w\in L(\mu)\). With these conventions we have
  \begin{equation}
    \phi_i=(-1)^i\frac{(\mu-i)!}{i! \mu!}e_0^i\phi_0.
    \label{eq:dbasisnorm}
  \end{equation}
  Then 
  \begin{equation}
    \mathrm{tr}_{L(k,\mu)_{h_\mu}} o\left(f_{-1}^{\frac{\lambda}{2}}\ket{\lambda} \right)
    =\sum_{i=0}^{\mu} \pair{\phi_i}{o
      \left(f_{-1}^{\frac{\lambda}{2}}\ket{\lambda} \right)v_i}.
  \end{equation}
  To further evaluate this expression we recall the Jacobi identity 
  \eqref{eq:jacID}. In that identity we set \(v=f_{-1}\wun\),
  \(U_1=L(k,\lambda)\), \(U_2=U_3=L(k,\mu)\), multiply both sides by
  \(z_0^{-1}\), and take the residue in \(z_0\) and \(z_1\) to obtain the
  identity
  \begin{equation}
    \mathcal{Y}(f_{-1}u_1,z_2)u_2=\sum_{s=0}^\infty
    z_2^{s}f_{-s-1}\mathcal{Y}(u_1,z_2)u_2+z_2^{-s-1}\mathcal{Y}(u_1,z_2)f_{s}u_2.
    \label{eq:jacspec}
  \end{equation}
  Specialising further to \(u_2=v_i\) and noting that for \(s\ge1\) we have \(f_{s}v_i=0\) and \(\phi_i(f_{-s}w)=0\)
  for any \(w\in L(k,\lambda)\), we obtain
  \begin{equation}
     \mathrm{tr}_{L(k,\mu)_{h_\mu}} o\left(f_{-1}^{\frac{\lambda}{2}}\ket{\lambda} \right)
    =\sum_{i=0}^\mu \pair{\phi_i}{ o(\ket{\lambda}) f_0^{\frac{\lambda}{2}}
      v_i}
    = \sum_{i=0}^{\mu-\frac{\lambda}{2}} \pair{\phi_i}{ o(\ket{\lambda})
        f_0^{\frac{\lambda}{2}} v_i}
    =\sum_{i=0}^{\mu-\frac{\lambda}{2}}\correct{(-1)^i\frac{(\mu-i)!}{i!\mu!}}\pair{e_0^i\phi_0}{ o(\ket{\lambda})
        f_0^{\frac{\lambda}{2}} v_i}.
  \end{equation}
  Here, the second equality is due to \(f_0^{\frac{\lambda}{2}}v_i\) vanishing for
  \(i > \mu-\frac{\lambda}{2}\), while the third equality follows from the
  formula above for the basis and its dual.
  Evaluating the action on the dual space and using the identity
  \([e_0,o\brac*{\ket{\lambda}}]=o\brac*{e_0\ket{\lambda}}=0\) (which follows
  similarly to \eqref{eq:jacspec} by taking appropriate residues of the Jacobi
  identity \eqref{eq:jacID}) we obtain
  \begin{equation}
  \label{eqn:CoeffCalc}
    \mathrm{tr}_{L(k,\mu)_{h_\mu}}
    o\left(f_{-1}^{\frac{\lambda}{2}}\ket{\lambda} \right)
    = \sum_{i=0}^{\mu-\frac{\lambda}{2}} \frac{(\mu-i)!}{i!\mu!}\pair{\phi_0}{o(\ket{\lambda}) e_0^i f_0^{\frac{\lambda}{2}+i} v_0}.
  \end{equation}
  Observe that
  \begin{equation}
    e^i_0 f_0^{i+\frac{\lambda}{2}} |\mu\rangle = \frac{(\frac{\lambda}{2}+i)!}{(\frac{\lambda}{2})!} \frac{(\mu-\frac{\lambda}{2})!}{(\mu-\frac{\lambda}{2}-i)!} f_0^{\frac{\lambda}{2}}|\mu\rangle,
  \end{equation}
  which combined with the Jacobi identity
  \(\comm{f_0}{o\brac*{w}}=o\brac*{f_0w}\) to move $f_0^{\frac{\lambda}{2}}$
  back into the intertwining operator zero mode yields
  \begin{equation}
    \mathrm{tr}_{L(k,\mu)_{h_\mu}}
    o\left(f_{-1}^{\frac{\lambda}{2}}\ket{\lambda} \right)
    =\pair{\phi_0}{ o\left(f_{0}^{\frac{\lambda}{2}} \ket{\lambda} \right)v_0}(-1)^{\frac{\lambda}{2}}
    \sum_{i=0}^{\mu-\frac{\lambda}{2}} \frac{(\mu-i)!}{i!\mu!}
    \frac{\left(\frac{\lambda}{2}+i\right)!
      \left(\mu-\frac{\lambda}{2}\right)!}{\left(\frac{\lambda}{2}\right)!\left(\mu-\frac{\lambda}{2}-i\right)!}
  \end{equation}
  Since $\mu \geq \frac{\lambda}{2}$ and $i \leq \mu-\frac{\lambda}{2}$, the sum is strictly positive and rational. Thus, \(\mathrm{tr}_{L(k,\mu)_{h_\mu}} o(f_{-1}^{\frac{\lambda}{2}}\ket{\lambda})\) is non-vanishing if and only if
  \(\langle {\phi_0},{o(f_0^{\frac{\lambda}{2}} \ket{\lambda})v_0} \rangle\)
  is, which in turn must be non-zero because the intertwining operator is.
  
  We turn to showing that the leading exponents provide an equality in \eqref{eq:bigweightineq}.
  Set $\mu_n = \frac{\lambda}{2} + n$ for $n=0, 1, \dots, k-\lambda$. Then
   $\intset_\lambda = \{\mu_0,\mu_1,\dots ,\mu_{k-\lambda} \}$. 
   Recalling $\lvert \intset_\lambda \rvert = k-\lambda+1$, we have
 \begin{equation}
   \frac{\sum_{n=0}^{k-\lambda} \left(h_{\mu_n}
       -\frac{\mathbf{c}}{24}\right)}{k-\lambda+1} = \frac{4k(k+2) -
     2\lambda(k+1) + \lambda^2}{48(k+2)},
 \end{equation}
and thus
\begin{equation}
	\frac{12 \left( \sum_{n=0}^{k-\lambda} h_{\mu_n} -\frac{\mathbf{c}}{24}\right)}{k-\lambda+1} +\lambda - k = \frac{\lambda(\lambda+2k+6)}{4(k+2)}.
\end{equation}
Finally, we note that the above is equal to
\begin{equation}
  h_\lambda+\frac{\lambda}{2}= \frac{\lambda(\lambda+2k+6)}{4(k+2)}
\end{equation}
and hence the equality in \eqref{eq:bigweightineq} is obtained.

  Next we show that \correct{the} intertwining operator can be normalised such that the
  trace \(\ptf^\mu(u,\tau)\) has rational coefficients. 
 First note that since
  the level \(k\) is integral, the commutation relations of the affine
  generators \(e_n,h_n,f_n\) all have integral structure constants. Further,
  the two generating singular vectors of the maximal proper submodule of the
  Verma module \(V(k,\mu)\) can be normalised to have integral expansions in
  the standard \PBW{} (PBW) basis. Therefore, a basis of the simple quotient \(L(k,\mu)\) can be chosen such that its representatives in \(V(k,\mu)\) expand in the standard PBW basis with rational coefficients. See \cite{MerPrim96} for a description of such bases.
  
Finally, note that when expressing dual basis vectors in simple finite \(\SLA{sl}{2}\)-modules in terms of the dual of the highest vector (as in \eqref{eq:dbasisnorm}) all normalisation factors are again rational. Thus every computation of \(\mathrm{tr}_{L(k,\mu)_{h_\mu +m}} o(f_{-1}^{\frac{\lambda}{2}}\ket{\lambda})\) will reduce to 
\(\langle {\phi_0},{ o(f_0^{\frac{\lambda}{2}} \ket{\lambda})v_0} \rangle\) multiplied by a sum of products of rational numbers, hence to ensure that the coefficients of \(\ptf^\mu(u,\tau)\) are rational it is necessary and sufficient to normalise the intertwining operator such that 
\(\langle {\phi_0},{ o(f_0^{\frac{\lambda}{2}} \ket{\lambda})v_0} \rangle\)
is rational, which can always be done.
\end{proof}

We now have all results needed to prove \cref{thm:sl2alldims}.

\begin{proof}[Proof of \cref{thm:sl2alldims}.]
  Recall that for all non-negative integers $m$ we have \(L(k,\lambda)_{[{h_\lambda+m}]}\subset \bigoplus_{n=0}^m L(k,\lambda)_{h_\lambda+n}\), hence Part \ref{itm:zeros} follows from \cref{thm:trivialmodulelambdaovertwo}.\ref{itm:mult0}. 

  To conclude Part \ref{itm:nonzeros} note that
  \cref{thm:trivialmodulelambdaovertwo}.\ref{itm:mult1} bounds the codimension
  of the subspace \(N\) above by \(1\), while \cref{thm:leadingcoefficient}.\ref{itm:leadingexp} bounds it below by \(1\), hence the codimension is \(1\). 

  Part \ref{itm:leadingexp} is given in \cref{thm:leadingcoefficient}.\ref{itm:leadingterm}, while Part \ref{itm:ratcoeff} is \cref{thm:leadingcoefficient}.\ref{itm:scaletorat}.

  Finally, for Part \ref{itm:basis}, we show linear independence of the traces
  \(\ptf^\mu\) evaluated at \(u\). Recall that the leading exponents are
\begin{equation}
    h_\mu - \frac{\mathbf{c}}{24} = \frac{2\mu^2+4\mu-k}{8(k+2)},\qquad \frac{\lambda}{2}\le\mu\le k-\frac{\lambda}{2}.
  \end{equation}
  Observe that the numerator is quadratic in \(\mu\) with a minimum at \(\mu=-1\) which is below the range of \(\mu\) hence all exponents are distinct. Thus, the set \(\set{\ptf^\mu(u,\tau)\st \mu\in\intset_\lambda}\) is linearly independent and thus so is \(\set{\ptf^\mu (-,\tau)\st \mu\in\intset_\lambda}\). 
  
  Finally, we turn to Part \ref{itm:holomcon}. By Theorem \ref{thm:TraceFunctionIsVVMF} we know $\vvptf_\lambda (w ,\tau)$ is a weakly holomorphic vector-valued modular form. It remains to show that if $\lambda\geq -2+\sqrt{2k+4}$, then all exponents of each component of $\vvptf_\lambda (w ,\tau)$ are non-negative. By Proposition \ref{thm:formdim}, the smallest possible leading exponent in the $q$-expansions among all components is $h_\mu -\mathbf{c}/24$ for $\mu = \lambda /2$. All other exponents are larger since $h_{\mu_1}-\mathbf{c}/24\geq h_{\mu_2}-\mathbf{c}/24$ if $\mu_1\geq \mu_2$. Thus, we are assured all exponents will be non-negative if $h_{\lambda /2}-\mathbf{c}/24 \geq 0$. By \eqref{eqn:LeadingExp} this is equivalent to 
 \begin{equation}
 \frac{\lambda^2 +4\lambda -2k}{16(k+2)} \geq 0,
 \end{equation}
 which in turn amounts to $\lambda^2 +4\lambda -2k\geq 0$. This establishes the holomorphicity of $\vvptf_\lambda (w,\tau)$, as desired.
 
 As discussed above, the smallest exponent occurring in any $q$-expansion in the components of $\vvptf_\lambda (u,\tau)$ is 
$h_{\lambda/2} -\mathbf{c}/24$. 
The same argument as above now gives the equivalence between Parts \ref{itm:hol}
and \ref{itm:wtbound}. Meanwhile, we recall that
\ref{thm:sl2alldims}.\ref{itm:basis} gives $d= k-\lambda +1$, or $\lambda =
k-d+1$. Plugging this into the inequality $\lambda^2 +4\lambda -2k\geq 0$ we
obtained above and solving for $k$ gives the equivalence between Parts \ref{itm:hol} and \ref{itm:lvlbound}.
\end{proof}

\section{Analysing representations and categorising spaces of torus 1-point
  functions}
\label{sec:analyzing}
In this section we present a detailed analysis of vector-valued
modular forms arising from the $1$-point functions of the simple affine \voa{s} constructed from $\SLA{sl}{2}$ at non-negative integer levels. We determine congruence or non-congruence for large families of examples, and also whether the representations associated to these forms have finite or infinite image. Additionally, we give complete descriptions of the spaces of vector-valued modular forms of dimension at most three.
  
Recall the ring of integral weight modular
forms \(\mathcal{M}=\CC[\eis{4},\eis{6}]\) and the skew polynomial ring of modular
differential operators \(\mathcal{R}\). Throughout this section let $0\leq
\lambda \leq k$, with $\lambda$ even, and
\(\vvptf_\lambda (u ,\tau)\)
denote the vector-valued modular form defined in \cref{thm:sl2alldims}.\ref{itm:basis}, with
\(u\in L(k,\lambda)\) \correct{torus primary} \correct{of the form} 
\eqref{eq:actingvec}. Consider the cyclic
$\mathcal{R}$-submodule $\mathcal{R}\vvptf_\lambda (u ,\tau)$ of
$\mathcal{M}^!(\rho_\lambda ,\nu_{h_\lambda})$ generated by $\vvptf_\lambda
(u,\tau)$. In this section we will consider how the \(\mathcal{R}\)-modules
\(\mathcal{R}\vvptf_\lambda (u ,\tau),\ \mathcal{V}^{u} \left(\rho_\lambda
\right),\
\mathcal{V} \left(\rho_\lambda \right),\ \mathcal{H}(\rho_\lambda
,\nu_{h_\lambda})\) and \(\mathcal{M}^!(\rho_\lambda ,\nu_{h_\lambda})\) are
interrelated.

\begin{prop}
Let $2\leq \lambda \leq k$, $\lambda$ even, and $u\in L(k,\lambda)$
\correct{a torus primary of the form} 
\eqref{eq:actingvec}. Then $\mathcal{V}^{u} \left(\rho_\lambda \right) = \mathcal{R}\vvptf_\lambda (u ,\tau)$.
  \label{prop:RmoduleResultsl2}
\end{prop}

\begin{proof}
 Suppose $w\in L(k,\lambda)_{[h_\lambda +\ell]}$ for some $\ell \in
 \mathbb{Z}$. If $\ell < \lambda /2$ we have $\vvptf_\lambda (w,\tau)=0$ by
 \cref{thm:sl2alldims}.\ref{itm:zeros}. The result now follows from
 \cref{prop:RmoduleResult} after noting that \(\operatorname{wt}[u]>0\).
\end{proof}

Next, we prepare a sufficient
condition for concluding the irreducibility of a representation of \(\SLTZ\).

\begin{lem}
  Let \(d\in \NN\) and let \(\rho\colon \SLTZ \to \GL{d}\) be a \(d\)-dimensional
  representation of \(\SLTZ\) such that \(\rho(T)\) is diagonalisable with
  eigenvalues \(\{\lambda_1,\dots, \lambda_d\}\). If
  every non-empty proper subproduct of \(\det (\rho(T))=\prod_{i=1}^d \lambda_i\) is 
  not a \(12\)th root of unity, then \(\rho\) is irreducible.
  \label{thm:irredcond}
\end{lem}

\begin{proof}
  Recall that taking the determinant \(\det\brac*{\rho}\) of \(\rho\) yields
  an element in the group of characters \(\operatorname{Hom}(\SLTZ ,\CC^\times)\) and
  that this group is cyclic of order $12$
 (one choice of cyclic generator assigns \(\modT\mapsto
\epi\brac*{\frac{1}{12}}\) \(\modS\mapsto \epi\brac*{\frac{3}{4}}\)).
In particular, every element in \(\operatorname{Hom}(\SLTZ ,\CC^\times)\) maps \(\modT\)
to some \(12\)th root of unity. Note
further that any invariant subspace of the representation \(\rho\) admits a
basis of \(\rho(\modT)\) 
eigenvectors so taking the determinant of the representation restricted to
this subspace will map \(T\) to a product of \(\rho(\modT)\)-eigenvalues with
as many factors as the dimension of the subspace and this product would need
to be a \(12\)th root of unity. So if no non-empty product
of \(\rho(\modT)\)-eigenvalues is a \(12\)th root of unity, then \(\rho\)
admits no non-trivial invariant subspace.
\end{proof}

\subsection{Dimension One}

We begin by considering \(1\)-dimensional \vvmf{s}.

 \begin{thm}
   Let $0\leq \lambda \leq k$, \(\lambda\) even, $u\in L(k,\lambda)$
   \correct{a torus primary of the form} 
   \eqref{eq:actingvec}, and $\rho_\lambda$ be the representation associated
   to $\vvptf_\lambda(u,\tau)$. The dimension of the vector-valued modular
   form \(\vvptf_\lambda (u,\tau)\) 
   is \(1\) if and only
   if $\lambda =k$ and hence the level \(k\) is even. Moreover, in this case the following hold.
 \begin{enumerate}
 \item The representation $\rho_{k}$ is irreducible and congruence. 
 \label{itm:1dimirreducibility}
 \item The representation $\rho_k$ satisfies
   \(\rho_k(\modS)=\epi\brac*{\frac{\correct{-}k}{8}}\) and \(\rho_{k}(\modT) =
   \epi\brac*{\frac{k}{24}}\). In particular, \(\rho_k\) is trivial if and
   only if $k$ is a multiple of $24$.
   \label{itm:1dimcongruence}
 \item We have the inclusion \(\mathcal{V}\left(\rho_{k}\right) \subset
   \mathcal{H}(\rho_{k},\nu_{h_k})\) for all even \(k\ge0\). The
   inclusion is an equality if \(2\le k\le 14\) and it is proper if \(k=0\)
   or \(k\ge 16\).
   \label{itm:1dholomcrit}
\item There exists a normalisation of the intertwining operator underlying
  \(\vvptf_k(u,\tau)\) such that
  \(\vvptf_k(u,\tau)=\eta^{\frac{3k}{2}}\) for all even
    \(k\ge0\). Further, for \(k\ge2\) we have the identity of
  \(\mathcal{R}\)-modules
  \begin{equation}
    \mathcal{V}^u\left(\rho_{k}\right) =\mathcal{V}\left(\rho_{k}\right) =
    \mathcal{R}\eta^{\frac{3k}{2}}.
    \label{eq:1dimV}
  \end{equation}
  As \(\mathcal{M}\)-modules each of the above is free of rank \(1\) with
  basis \(\set{\eta^{\frac{3k}{2}}}\).
  For $n\in \NN_0$, 
   \begin{equation}
     \dim \left(\mathcal{V}\left(\rho_{k}\right)_{\correct{\frac{k}{2}+}n} \right)
     = 
     \begin{cases}
       0&\text{ if } n\equiv 1\pmod2\\
       \floor*{\frac{n}{12}} & \text{ if } n \equiv 2\pmod{12} \\
       \floor*{\frac{n}{12}} +1 & \text{ otherwise,}
     \end{cases}
     \label{eq:onedimdimformula}
   \end{equation}
    where $\lfloor x
\rfloor$ denotes the floor
function of a real number $x$. 
\label{itm:onedimdimformula}
 \end{enumerate}
  \label{thm:1dimcongruence}
\end{thm}

\begin{proof}
 Theorem \ref{thm:sl2alldims}.\ref{itm:basis} gives that $\rho_{k}$ is a
 $1$-dimensional representation for even \(k\).
  Further, $\rho_k$ is irreducible as it is $1$-dimensional and we note any
  irreducible $1$-dimensional representation
  must also be congruence (see, for example,
  \cite[\refSec 3.1]{intvvmf18} for more
  details).
  This gives Part \ref{itm:2dimirreducibility}.

  For Part \ref{itm:2dimcongruence}, note that the representation \(\rho_k\) is trivial if and only \(\rho_k(\modT)=1\).
  By \cref{thm:sl2alldims}.\ref{itm:basis}, \(\rho_{k}(\modT) = \epi\brac*{\frac{k}{24}}\), thus
  the representation is trivial if and only if $k$ is a multiple of $24$.
  
  We consider Parts \ref{itm:1dholomcrit} and \ref{itm:onedimdimformula} together. 
  In the case $k=0$, it is well known that $L(0,0)\cong \mathbb{C}$, and it
  follows that \(\mathcal{V}\left(\rho_{k}\right)=\CC\), which is strictly
  contained in $\mathcal{H}(\rho_0,\nu_{h_0})$. Therefore, we consider the
  case $k>0$, and thus $\operatorname{wt}[u]>0$ by Theorem \ref{thm:sl2alldims}.
  Note for \correct{irreducible representations of} dimensions three or less
  spaces of holomorphic \vvmf{s} are always cyclic \(\mathcal{R}\)-modules
  \correct{\cite{Marks2011}}, and 
  hence all \correct{five} 
  conditions in \cref{thm:generaldimformula} are
  satisfied, which implies
  \begin{equation}
    \mathcal{V}^u\left(\rho_k\right)=\mathcal{V}\left(\rho_k\right)=\eta^{\frac{3k}{2}}\mathcal{H}\left(\rho_k,\nu_{-\frac{k}{2}}\right)
  \end{equation}
  and the dimension
  formula \eqref{eq:onedimdimformula} follows. Further, the cyclic generator of
  \(\mathcal{H}(\rho_k,\nu_{-\frac{k}{2}})\) has weight \(0\) and thus can be
  chosen to be \(1\), which gives the formula
  \(\vvptf_k(u,\tau)=\eta^{\frac{3k}{2}}\).
  \correct{Finally, the leading exponent of  \(\vvptf_k(u,\tau)\) is
  \begin{equation}
    h_{\frac{k}{2}}-\frac{\mathbf{c}}{24}=\frac{k}{16}.
  \end{equation}
  For \(k\ge2\) the inclusion \(\mathcal{V}\left(\rho_{k}\right) \subset
   \mathcal{H}(\rho_{k},\nu_{h_k})\) is therefore proper if and only if the
   leading exponent is at least \(1\) (and hence no longer minimal admissible), that is, if and only if \(k\ge16\).}
 \end{proof}

\subsection{Dimension Two}

We turn to describing the $2$-dimensional setting and prepare some notation.
Let \(j(\tau)\) be Klein's \(j\)-invariant (normalised so that the leading
term is \(q^{-1}\)) and $J(\tau) = j(\tau)/1728$.
Additionally, for $a,b,c\in \mathbb{C}$, $c$ not a negative integer, and a variable $z$, let ${}_2F_1 (a,b;c;z)$ denote the Gaussian hypergeometric function, which is given by
\begin{equation}
{}_2F_1 (a,b;c;z) =  1+ \sum_{n=1}^\infty \frac{(a)^n(b)^n}{(c)^n}\frac{z^n}{n!}
\end{equation}
with $(x)^n$ being the (rising) factorial for $x\in \mathbb{C}$ given by $(x)^n = x(x+1)\cdots (x+n-1)$ for $n\in \mathbb{N}$.
Finally, set
 \begin{equation}
   \Phi  =
   \begin{pmatrix}
     J^{\frac{1}{24}}{}_2F_1\left(\frac{-1}{24}, \frac{7}{24}; \frac34; J^{-1} \right) \\
     J^{-\frac{5}{24}} {}_2F_1\left(\frac{5}{24}, \frac{13}{24}; \frac{5}{4}; J^{-1} \right)
   \end{pmatrix},
   \label{eq:2dhgeom1}                                  
 \end{equation}
  where \(J^{-1}=1/J\).

\begin{thm}
 \label{thm:2dimresults}
 Let $0\leq \lambda \leq k$, \(\lambda\) even, $u\in L(k,\lambda)$
 \correct{torus primary of the form} 
 \eqref{eq:actingvec}, and let $\rho_\lambda$ be the representation associated
 to $\vvptf_\lambda (u,\tau)$. 
 The dimension of the vector-valued modular form \(\vvptf_\lambda
 (u,\tau)\) is
 \(2\) if and only if $\lambda =k-1$ and hence the level \(k\) is odd.
 \begin{enumerate}
 \item The representation $\rho_{k-1}$ is irreducible. Moreover, among
   all indecomposable representations $\rho'$ of 
    \(\SLTZ\) satisfying $\rho' (T)=\operatorname{diag}(\epi \left(
     \frac{k-2}{24}\right), \epi \left( \frac{k+4}{24}\right))$,
   \(\rho_{k-1}\) is the unique (up to isomorphism) one that is irreducible.
 \label{itm:2dimirreducibility}
 \item The representation $\rho_{k-1}$ is congruence with congruence level
   \(N=8\) for \(k\equiv 2\pmod{3}\) and \(N=24\) otherwise.
 \label{itm:2dimcongruence}
\item We have an inclusion \(\mathcal{V}\left(\rho_{k-1}\right)\subset
  \mathcal{H}(\rho_{k-1},\nu_{h_{k-1}})\) for all odd \(k\ge3\). This
  inclusion is an equality if \(3\le k\le 13\) and proper if \(k\ge 15\).
  \label{itm:2dimholomcrit}
\item There exists a normalisation of the intertwining operator underlying
  \( \vvptf_{k-1}(u,\tau)\) such that
  \begin{equation}
    \vvptf_{k-1}(u,\tau)=\eta^{\frac{3k^2+2k-5}{2(k+2)}}\Phi
    \label{eq:2dimcycgenformula}
  \end{equation}
  for all odd \(k\ge1\).
  Further, for \(k\ge3\) we have the identity of \(\mathcal{R}\) modules
  \begin{equation}
    \mathcal{V}^u\left(\rho_{k-1}\right) = \mathcal{V}\left(\rho_{k-1}\right) = 
    \mathcal{R}\eta^{\frac{3k^2+2k-5}{2(k+2)}}\Phi.
  \end{equation}
  As \(\mathcal{M}\)-modules each of the above is free of rank \(2\) with
  basis
  \(\set{ \vvptf_{k-1}(u,\tau),\pd \vvptf_{k-1}(u,\tau)}=\set{\vvptf_{\correct{k-1}}(u,\tau), \vvptf_{k-1}(L_{[-2]}u,\tau)}\). For
  \(n\in\NN_0\),
  \begin{equation}
    \dim\brac*{\mathcal{V}\left(\rho_{k-1}\right)_{\correct{\frac{k-1}{2}+}n}}=
    \begin{cases}
      0&\text{ if } n\equiv 1\pmod{2}\\
      \floor*{\frac{n}{6}}+1&\text{ otherwise}.
    \end{cases}
    \label{eq:2dimdimformula}
  \end{equation}
 \label{itm:2dimdimformula}
\end{enumerate}
\label{itm:2dimHolomResults}
\end{thm}

\begin{proof}
 That $\rho_{k-1}$ is a $2$-dimensional representation if and only if $\lambda =k-1$ is even follows directly from Theorem \ref{thm:sl2alldims}.\ref{itm:basis}.
 To show Part \ref{itm:2dimirreducibility} we use the criterion in
 \cref{thm:irredcond} and note by \cref{thm:sl2alldims}.\ref{itm:basis} that
  \begin{equation}
    \rho_{k-1}(\modT) = \diag \left(
      \epi \left(  \frac{k-2}{24}\right),
      \epi \left( \frac{k+4}{24}\right)
    \right) .
    \label{eq:2dTformula}
  \end{equation}
  Thus, by \cref{thm:irredcond}, \(\rho_{k-1}\) is irreducible if neither
  of the 
  \(\rho_{k-1}(\modT)\)-eigenvalues are a \(12\)th root of unity. This is
  clearly the case, since \(k\) is odd.
  The fact that \(\rho_{k-1}\) is the unique irreducible representation among
  indecomposable representations with \(\modT\) given by the formula
  \eqref{eq:2dTformula} is due to \cite[\refThm 3.1]{Mason2007}. 
  
  Part \ref{itm:2dimcongruence} follows from \cite[\refThm 3.7]{Mason2007},
    \correct{where all} $2$-dimensional irreducible finite image representations are
    classified. They turn out to all be congruence representations. The
    congruence levels are recorded in the tables following that 
    theorem, where each representation is characterised by the fractions
    (or rather the smallest non-negative representative of their integer
    coset) that appear in the formula \eqref{eq:2dTformula} for
    \(\rho_{k-1}(\modT)\). The congruence level is then always the order of
    \(\rho_{k-1}(\modT)\), that is, \(8\) if \(k\equiv
    2\pmod{3}\) and \(24\) otherwise. A simple
    calculation reveals that for all odd \(k\)
    each \(\rho_{k-1}(\modT)\) corresponds to a case in \cite[\refTab 3]{Mason2007}.
    
 We consider Parts \ref{itm:2dimholomcrit} and \ref{itm:2dimdimformula} together. By Theorem
 \ref{thm:sl2alldims}.\ref{itm:basis}
 the leading exponents of \(\vvptf_{k-1} (u,\tau)\) are
 \begin{equation}
   \mu_{\min{}}=h_{\frac{k-1}{2}}-\frac{\mathbf{c}}{24}=\frac{k^2-3}{16(k+2)},\qquad
   \mu_{\max{}} = h_{\frac{k+1}{2}}-\frac{\mathbf{c}}{24} = \frac{k^2+4k+5}{16(k+2)}=\frac{k^2-3}{16(k+2)}+\frac{1}{4}.
 \end{equation}
 Clearly these exponents are non-negative if and only if \(k\ge3\), which
 proves the inclusion in Part \ref{itm:2dimholomcrit}. Since spaces of
 holomorphic \vvmf{s} \correct{for irreducible representations of} of dimension three or less are always cyclic
 \(\mathcal{R}\) modules \correct{\cite{Marks2011}}, and since the two exponents above differ by
 \(\frac{1}{4}\), all \correct{five} assumptions of \cref{thm:generaldimformula}
 apply if \(k\ge3\).  Hence
 \begin{equation}
   \mathcal{V}^u\left(\rho_{k-1}\right) = \mathcal{V}\left(\rho_{k-1}\right) =
   \eta^{3\frac{k^2-3}{2(k+2)}}\mathcal{H}\brac*{\rho_{k-1},\nu_{\frac{2-k}{2}}}
 \end{equation}
  for \(k\ge 3\). For \(k=1\)\correct{, we have $\operatorname{wt}[u]=0$ and cannot invoke Theorem \ref{thm:generaldimformula}. However}, by \cref{prop:RmoduleResult},
   we can still assert
 \(\mathcal{V}^u\left(\rho_{k-1}\right)\subset
 \eta^{3\frac{k^2-3}{2(k+2)}}\mathcal{H}(\rho_{k-1},\nu_{\frac{2-k}{2}})\).
 Further, for all odd \(k\ge 1\) we have that
 \(\eta^{-3\frac{k^2-3}{2(k+2)}}\vvptf(u,\tau)\) is a cyclic generator for
 \(\mathcal{H}(\rho_{k-1},\nu_{\frac{2-k}{2}})\), since it has the right
 weight and spans a \(1\)-dimensional weight space.
 By construction, a cyclic generator of
 \(\mathcal{H}(\rho_{k-1},\nu_{\frac{2-k}{2}})\) has leading exponents
 \(\set{0,\frac{1}{4}}=\set{\lambda_1,\lambda_2}\) (these exponents are a minimal admissible set) and the
   weight of this cyclic generator is \(\frac{1}{2}\), in particular it
   satisfies the equality in \eqref{eq:weightineq} and so its components form
   a fundamental system of solutions to a monic modular differential equation.
   This allows us to use \cite[\refThm 3.1]{Marks2011}
   and its proof (which additionally requires $\rho (\modS^2)$ be a scalar
   matrix, but this is automatic due to \(\rho_{k-1}\) being irreducible; see
   the discussion below \eqref{eq:rhoT}).
   This gives the existence of a non-zero vector-valued modular form $F$
   of weight \(\mdwt=6(\lambda_1+\lambda_2)-1=\frac{1}{2}\) such that
     $\mathcal{H}(\rho_{k-1} ,\nu_{h_{k-1}})=\mathcal{R}F$ (as
     $\mathcal{R}$-modules) and $\mathcal{H}(\rho_{k-1}
     ,\nu_{h_{k-1}})=\mathcal{M}F \oplus \mathcal{M}\partial F$ (as
     $\mathcal{M}$-modules). 
     Moreover, the component functions of $F$ form a fundamental system of solutions of a second order monic modular differential equation of the form
 \begin{equation}
  \left(\partial_{\mdwt}^2  +\phi \right) f =0,
\label{eqn:2ndOrderMMDE01}
\end{equation}
 where $\phi \in \mathcal{M}_4$. Note that since $\mathcal{M}_4 =
 \operatorname{span}_{\mathbb{C}}\{\eis{4}\}$, up to a scalar $\mdwtt$ we have that \eqref{eqn:2ndOrderMMDE01} can be rewritten as
  \begin{equation}
 \label{eqn:2ndOrderMMDE02}
 \left(\partial_{\mdwt}^2  - 
   \mdwtt \eis{4}\right) f =0.
 \end{equation}
This equation is characterised by its indicial roots being \(\lambda_1\) and
\(\lambda_2\) (these are related to \(\mdwtt\) via
  \(\mdwtt=180\brac*{\lambda_1-\lambda_2}^2-5\))
which are the exponents of the first and second component functions of $F$, respectively, as detailed in the proof of \cite[Theorem 3.1]{Marks2011} (cf.\ \eqref{eq:VVMFq-expansion} and \eqref{eq:MinimalAdmissibleRelation}).
Meanwhile, \cite[\refProp
 2.2]{FrancMason2014} (see also \cite[Section 4.1]{FrancMason-Hypergeometric})
 gives that the functions
  \begin{align}
 f_1 &= \eta^{2\mdwt} J^{-\frac{6\left(\lambda_1 -\lambda_2\right)+1}{12}} {}_2F_1 \left(\frac{6\left(\lambda_1 -\lambda_2\right)+1}{12}, \frac{6\left(\lambda_1 -\lambda_2\right)+5}{12}; \lambda_1-\lambda_2+1 ; J^{-1} \right)
       \nonumber\\
    f_2 &= \eta^{2\mdwt} J^{-\frac{6\left(\lambda_2
          -\lambda_1\right)+1}{12}} {}_2F_1
          \left(\frac{6\left(\lambda_2 -\lambda_1\right)+1}{12},
          \frac{6\left(\lambda_2 -\lambda_1\right)+5}{12};
          \lambda_2-\lambda_1+1 ; J^{-1} \right)
   \label{eq:2dhypergeom}
 \end{align}
 form a fundamental set of solutions for \eqref{eqn:2ndOrderMMDE02}. We make some notes pertaining to our use of \cite{FrancMason2014, FrancMason-Hypergeometric}. First, loc.\ cit.\ assumes integral weights, however, a careful examination
 of the proof shows that it also holds for real weights.
 Second, there is a difference of normalisations of Eisenstein series between the
 $\eis{2k}$ in this paper and the $E_{2k}$ of \cite{FrancMason2014, FrancMason-Hypergeometric} given by
 $E_{2\ell} = -\frac{(2\ell)!}{B_{2\ell}}\eis{2\ell}$ for $\ell \in
 \mathbb{N}$.
 Since $f_1$ and $f_2$ form a fundamental set of solutions for
 \eqref{eqn:2ndOrderMMDE02}, up to a matrix $A\in
 \operatorname{GL}(2,\mathbb{C})$, we have $A(f_1,f_2)^t =F$. That is,
 $(f_1,f_2)^t$ is a vector-valued modular form of weight $\mdwt$, but with
 representation $A\rho_{k-1}A^{-1}$. However, the leading exponents of $f_1$
 and $f_2$ are $\lambda_1$ and $\lambda_2$, respectively (we note that this
 disagrees with \cite[\refRmk 2.3]{FrancMason2014}, where there is a minor
 typographical error listing the exponents in reverse order).
 Thus, it must be that $A = \operatorname{diag}(\alpha ,\beta)$ for some $\alpha ,\beta \in \mathbb{C}^\times$ and we have $A\rho_{k-1}A^{-1}= \rho_{k-1}$. In particular, $F_{\alpha,\beta} = (\alpha f_1,\beta f_2)^t$.
 It remains to show that $(\alpha f_1,\beta f_2)^t$ is equal to \eqref{eq:2dimcycgenformula}. This
 follows immediately from \eqref{eq:2dhypergeom} by specialising
 \(\mdwt=\frac{1}{2}\), \(\lambda_1=0\), \(\lambda_2=\frac{1}{4}\). The
 dimension formula \eqref{eq:2dimdimformula} then follows from the evaluation of
 \eqref{eq:strongdimform} \cite[\refCor 3.2]{Marks2011}. Finally the inclusion
 of Part \ref{itm:2dimholomcrit} is an equality if and only the leading
 exponents of \(\vvptf\brac*{u,\tau}\) lie in the interval \([0,1)\) which
 happens if and only of \(3\le k\le 13\).
\end{proof}

See \cref{tab:2dtable} for explicit expansions of \(\vvptf_{k-1}(u,\tau)\) for
the first few values of the level \(k\).

\begin{table}[ht]
\renewcommand{\arraystretch}{1.5}
\begin{tabular}{|c||l|l|l|}
\hline
  Level \(k\)& Cyclic generator \(\vvptf_{k-1}(u,\tau)\)\\ \hline\hline
	\(3\)&  \(\begin{matrix*}[l]
		q^{3/40}\left( 1 + \frac15 q - \frac{117}{25} q^2 - \frac{84}{125} q^3 +\frac{3659}{625} q^4 +\cdots \right)\\
		q^{13/40}\left(1 - \frac95 q - \frac{2}{25} q^2 -\frac{39}{125} q^3 -\frac{126}{625} q^4 +\cdots \right)
	\end{matrix*}\)\\ \hline
	\(5\) &  \(\begin{matrix*}[l]
		q^{11/56} \left( 1 - \frac{19}{7}q - \frac{264}{49}q^2 + \frac{6061}{343}q^3 + \frac{22963}{2401} q^4 + \cdots\right)\\
		q^{25/56} \left(1 -\frac{33}{7}q + \frac{247}{49}q^2 + \frac{1672}{343}q^3 -\frac{18183}{2401}q^4 + \cdots \right)
	\end{matrix*}\)\\ \hline
	\(7\) & \(\begin{matrix*}[l]
		q^{23/72}\left( 1 -\frac{17}{3}q + \frac{23}{9}q^2 + \frac{3128}{81}q^3 - \frac{13429}{243}q^4+ \cdots\right)\\
		q^{41/72}\left(1 - \frac{23}{3}q + \frac{170}{9}q^2 - \frac{391}{81}q^3 - \frac{10948}{243}q^4 + \cdots \right)
	\end{matrix*}\)\\ \hline
	\(9\) & \(\begin{matrix*}[l]
		q^{39/88} \left( 1 - \frac{95}{11}q + \frac{2340}{121}q^2 + \frac{48165}{1331}q^3 - \frac{2895523}{14641}q^4 + \cdots\right)\\
		q^{61/88} \left( 1 - \frac{117}{11}q + \frac{5035}{121}q^2 - \frac{74100}{1331}q^3 - \frac{1011465}{14641}q^4+\cdots\right)
	\end{matrix*}\)\\ \hline
	\(11\) & \(\begin{matrix*}[l]
		q^{59/104} \left(1 - \frac{151}{13}q + \frac{7611}{169}q^2- \frac{35636}{2197}q^3 - \frac{9959957}{28561}q^4 + \cdots \right)\\
		q^{85/104} \left(1 - \frac{177}{13}q +\frac{12382}{169}q^2 -\frac{383087}{2197}q^3 + \frac{1229442}{28561}q^4 + \cdots \right)
	\end{matrix*}\)\\ \hline
	\(13\) & \(\begin{matrix*}[l]
		q^{83/120} \left( 1 - \frac{73}{5}q + \frac{1992}{25}q^2 - \frac{18177}{125}q^3 - \frac{224261}{625}q^4 + \cdots\right)\\
		q^{113/120} \left( 1 - \frac{83}{5}q + \frac{2847}{25}q^2 - \frac{48472}{125}q^3 + \frac{309009}{625}q^4+ \cdots\right)
	\end{matrix*}\)\\ \hline
\end{tabular}
\caption{
  The first five terms of the \(q\)-series expansions for
 \(\vvptf_{k-1}(u,\tau)\) for all levels $k$ at which \(\vvptf_{k-1}(u,\tau)\) generates \(\mathcal{H}(\rho_{k-1},\nu_{h_{k-1}})\).
 In each case the series have been normalised so that
 the leading coefficient is 1. This can always be achieved by an appropriate
 choice of normalisation of the intertwining operators in
 \(\vvptf_{k-1}(u,\tau)\).}
\label{tab:2dtable}
\end{table}

\subsection{Dimension Three}

Here we consider the $3$-dimensional case. 
Recall, ${}_3F_2$, the generalised hypergeometric function, which for
 $a,b,c,d,e\in \mathbb{C}$, $d,e$ not negative integers, and a variable $z$,
 is given by
\begin{equation}
{}_3F_2 (a,b,c;d,e;z) =  1+ \sum_{n=1}^\infty \frac{(a)^n(b)^n(c)^n}{(d)^n(e)^n}\frac{z^n}{n!}.
\end{equation}
We set 
\begin{equation}
 \Phi_{k}
 = \begin{pmatrix}
     J^{\frac{(k+1)}{12(k+2)}} {}_3F_2\left(-\frac{(k+1)}{12(k+2)}, \frac{11k+14}{24(k+2)}, \frac{19k+30}{24(k+2)}; \frac{3k+7}{4(k+2)}, \frac{1}{2}; J^{-1} \right)\\
     J^{-\frac{k+1}{6(k+2)}} {}_3F_2\left(\frac{k+1}{6(k+2)}, \frac{3k+5}{6(k+2)}, \frac{5k+9}{6(k+2)}; \frac{5k+9}{4(k+2)}, \frac{5}{8}; J^{-1} \right)\\
     J^{-\frac{5k+11}{12(k+2)}} {}_3F_2\left(\frac{5k+11}{12(k+2)}, \frac{9k+19}{12(k+2)}, \frac{13k+27}{12(k+2)}; \frac{3}{2}, \frac{5k+11}{4(k+2)}; J^{-1} \right)
   \end{pmatrix},
   \label{eq:threedimvoavvmfformula}
\end{equation}
where \(J=j/1728\) is the same renormalisation of Klein's
\(j\)-invariant as in the previous section.

\begin{thm}
 \label{thm:3dimresults}
 Let $0\leq \lambda \leq k$, \(\lambda\) even, $u\in L(k,\lambda)$
 \correct{torus primary of the form} 
 \eqref{eq:actingvec}, $\rho_\lambda$ the representation associated to
 $\vvptf_\lambda (u,\tau)$. The dimension of the vector-valued modular form \(\vvptf_\lambda
 (u,\tau)\) is \(3\) if and only if $\lambda =k-2$ and hence the level \(k\) is even.
 \begin{enumerate}
 \item The representation $\rho_{k-2}$ is irreducible. 
 \label{itm:3dimirreducibility}
\item The representation $\rho_{k-2}$ has finite image. 
   Additionally, the order of \(\rho_{k-2}(\modT)\) is $12(k+2)$ if $k\equiv 4\pmod{6}$, and is $4(k+2)$ otherwise.
 \label{itm:3dimorder}
 \item If the order of \(\rho_{k-2}(\modT)\) does not divide $25,401,600=2^8
   \cdot 3^4 \cdot 5^2 \cdot 7^2$, then the representation is non-congruence,
  in particular, this gives an infinite family of non-congruence
  representations and a finite bound on the number of congruence representations.
 \label{itm:3dimcongruence}
\item We have an inclusion \(\mathcal{V}(\rho_{k-2}) \subset
  \mathcal{H}(\rho_{k-2},\nu_{h_{k-2}})\) for all even
  \(k\ge4\). This inclusion is an equality if \(4\le k \le 10\) and it is
  proper if \(k\ge 12\).
  \label{itm:3dimholomcrit}
\item There exists a normalisation of the intertwining operators underlying
  \(\vvptf_{k-2}\brac*{u,\tau}\) such that
  \begin{equation}
    \vvptf_{k-2}\brac*{u,\tau}=\eta^{\frac{3k^2-2k-8}{2(k+2)}}\Phi_k
    \label{eq:3dimptfformula}
  \end{equation}
  for all even \(k\ge2\). Further, for \(k\ge4\) we have the identity of
  \(\mathcal{R}\)-modules
  \begin{equation}
    \mathcal{V}^u\left(\rho_{k-2}\right) = \mathcal{V}\left(\rho_{k-2}\right)
    = \mathcal{R}\eta^{\frac{3k^2-2k-8}{2(k+2)}}\Phi_k.
    \label{eq:3dimRmod}
  \end{equation}
  As \(\mathcal{M}\)-modules each of the above is free of rank $3$
  with basis
  \begin{equation}
    \set*{\vvptf_{k-2}\brac*{u,\tau},\pd\vvptf_{k-2}\brac*{u,\tau},\pd^2\vvptf_{k-2}\brac*{u,\tau}}
    =\set*{\vvptf_{k-2}\brac*{u,\tau},\vvptf_{k-2}\brac*{L_{[-2]}u,\tau},\vvptf_{k-2}\brac*{L_{[-2]}^2u+\delta
        L_{[-4]}u,\tau}},
    \label{eq:3dimbasis}
  \end{equation}
  where
  \(\delta=\frac{2\left(16+k-6k^2\right)}{3(k-2)(4+3k)}\). For \(n\in \NN_0\),
  \begin{equation}
    \dim\brac*{\mathcal{V}\left(\rho_{k-2}\right)_{\correct{\frac{k-2}{2}+}n}}=
    \begin{cases}
      0,&\text{ if } n\equiv 1\pmod2,\\
      \floor*{\frac{n}{4}}+1,&\text{ otherwise.}
    \end{cases}
  \end{equation}
  \label{itm:3dimdimformula}
 \end{enumerate}
\end{thm}

\begin{proof}
 Note that Theorem \ref{thm:sl2alldims}.\ref{itm:basis} gives that $\rho_{\lambda}$ is a $3$-dimensional representation if and only if $\lambda =k-2$ is even.
 
 To establish Part \ref{itm:3dimirreducibility} we use the irreducibility criterion
 in \cref{thm:irredcond}. Note that \cref{thm:sl2alldims}.\ref{itm:basis}
 yields the formula 
  \begin{equation}
 	\rho_{k-2}(\modT) = \diag\left(
 	\epi \left(\frac{k(k-2)-6}{24(k+2)} \right) ,
 	\epi \left(\frac{k(k+4)}{24(k+2)} \right),
 	\epi \left(\frac{k(k+10)+18}{24(k+2)} \right)
 	\right).
 	\label{eq:threedim}
 \end{equation}
The fractions
in the formula for $\rho_{k-2}(\modT)$ above multiplied by \(12\) are,
respectively,
\begin{equation}
  \frac{k}{2}-2 + \frac{1}{k+2} \equiv \frac{1}{k+2} \pmod{1},\quad
  \frac{k}{2}+\frac{k}{k+2}\equiv \frac{k}{k+2} \pmod{1},\quad
  \frac{k}{2}+4+\frac{1}{k+2}\equiv \frac{1}{k+2} \pmod{1}.
\end{equation}
These are never integral since the numerators on the \rhs{s} of the above
identities are always strictly less than the denominators. \correct{Proper sub}products of two
eigenvalues correspond to sums of two of the above fractions. The numerator of
such sums cannot exceed $k+1$, so they are never integral for any \(k\in 2\NN\),
concluding no twelfth root of unity arises as a product of one or two
\(\rho_{k-2}(\modT)\) eigenvalues. Thus \(\rho_{k-2}\) is irreducible.

Moving to Part \ref{itm:3dimorder}, we show that $\rho_{k-2}$ has finite image
by using the criterion \cite[\refProp 5.1]{Marks2015}, which states that a
\(3\)-dimensional irreducible representation with diagonalisable
\(\rho(\modT)\) has finite image if there exist two eigenvalues of
\(\rho(\modT)\) whose ratio is \(-1\). To this end, observe that
  \begin{equation}
    \frac{k(k+10)+18}{24(k+2)} - \frac{k(k-2)-6}{24(k+2)} = \frac12.
  \end{equation}
  Hence \(\rho_{k-2}\) has finite image.
  To determine the order of \(\rho_{k-2}(T)\), note that it is the least common
  multiple of the 
  denominators of the reductions of the fractions in the exponents in \eqref{eq:threedim}. Since $k$ is even, we take $k = 2\ell$ for $\ell \geq 0$ and
  consider all three exponents
  \begin{equation}
    \frac{k(k-2)-6}{24(k+2)} = \frac{2\ell(\ell-1)-3}{24(\ell+1)},\qquad
    \frac{k(k+4)}{24(k+2)}=\frac{\ell(\ell+2)}{12(\ell+1)},\qquad
    \frac{k(k+10)+18}{24(k+2)}=\frac{2\ell(\ell+5)+9}{24(\ell+1)}.
    \label{eq:partredfrac}
  \end{equation}
  Next we compute the greatest common divisor of the
  numerator and denominator of the first and third fractions above to reduce them.
  Denote $a = 2\ell(\ell-1)-3 = 2(\ell-2)(\ell+1) +1$,
  \(b=2\ell(\ell+5)+9=2(\ell+1)(\ell+4)+1\), and $c = 24(\ell+1)$ so that the
  first and third fractions are equal to \(\frac{a}{c}\) and
\(\frac{b}{c}\), respectively. Note that \(a,b\) are odd while \(c\) is even, so
\(\gcd\brac*{a,c}\) and \(\gcd\brac*{b,c}\) will both be odd. Further, if a prime \(p\ge5\)
divides \(c\), then it must divide \((\ell+1)\), but then \(a \equiv 1
\pmod{p} \equiv b\). So \(p\) does not divide \(\gcd\brac*{a,c}\) or
\(\gcd\brac*{b,c}\). Both of these greatest common divisors are therefore a power of
\(3\). We have that \(3\) divides \(a\) if and only if \(\ell \equiv 0,1
\pmod{3}\) and the same is also true for \(b\). Further, if \(\ell \equiv 0,1
\pmod{3}\), then \(\ell +1 \equiv 1,2\pmod{3}\), and so \(9\) does not divide
\(c\). Thus,
\begin{equation}
  \gcd\brac*{a,c}=\gcd\brac*{b,c}=
  \begin{dcases}
    1& \ell\equiv 2\pmod{3}\\
    3& \ell\equiv 0,1\pmod{3}.
  \end{dcases}
\end{equation}
Therefore, after reduction, the denominators of the first and third fractions in
\eqref{eq:partredfrac} are \(24(\ell+1)=12(k+2)\) if \(k\equiv 4\pmod{6}\), and
\(8(\ell+1)=4(k+2)\) if \(k\equiv 0,2\pmod{6}\). Next we see that the reduced
denominator of the middle fraction divides \(12(\ell+1)\) if \(\ell
\equiv2\pmod{3}\), and that it divides \(4(\ell+1)\) if \(\ell\equiv
0,1\pmod{3}\). It follows that the least common multiple of the reduced divisors in
\eqref{eq:partredfrac}, the order of \(\rho_{k-2}(\modT)\), is as claimed.

  Part \ref{itm:3dimcongruence} follows from \cite[\refCor 3.5]{Marks2015},
  which states that $\rho_{k-2}$ is non-congruence if there exists a prime
  dividing
  $\frac{\correct{o(\rho_{k-2}(T))}}{\left(\correct{o(\rho_{k-2}(T))},
      2^8 \cdot 3^4 \cdot 5^2 \cdot 7^2\right)}$\correct{, where
    \(\correct{o(\rho_{k-2}(T))}\) is the order of \(\rho_{k-2}(T)\)}. 
  \correct{Equivalently, $\rho_{k-2}$ is non-congruence if 
    \(o(\rho_{k-2}(T))\) does not divide \(2^8 \cdot 3^4 \cdot 5^2 \cdot 7^2\).}
  
  We consider Parts \ref{itm:3dimholomcrit} and \ref{itm:3dimdimformula} together. By
  \cref{thm:sl2alldims}.\ref{itm:basis} the leading exponents of
  \(\vvptf_{k-2}\brac*{u,\tau}\) are
  \begin{equation}
    \mu_{\min{}}=h_{\frac{k-2}{2}}-\frac{\mathbf{c}}{24}=\frac{k^2-2k-4}{16(k+2)},\quad
    h_{\frac{k}{2}}-\frac{\mathbf{c}}{24}=\frac{k}{16}=\mu_{\min{}}+\frac{k+1}{4(k+2)},\quad
    \mu_{\max{}}=h_{\frac{k+2}{2}}-\frac{\mathbf{c}}{24}=\frac{k^2+6k+12}{16(k+2)}=\mu_{\min{}}+\frac12.
  \end{equation}
  Clearly these exponents are non-negative if and only if \(k\ge4\), which
  proves the inclusion in Part  \ref{itm:3dimholomcrit}. Spaces of holomorphic
  \vvmf{s} \correct{for irreducible representations} of dimension three or less are always cyclic over \(\mathcal{R}\) \correct{\cite{Marks2011}},
  and the maximal and minimal exponents above differ by \(\frac12\). Therefore, all
  \correct{five} 
  assumptions of \cref{thm:generaldimformula} apply if \(k\ge4\). Hence
  \begin{equation}
    \mathcal{V}^u\left(\rho_{k-2}\right) = \mathcal{V}\left(\rho_{k-2}\right)
    =\eta^{3\frac{k^2-2k-4}{2(k+2)}}\mathcal{H}\brac*{\rho_{k-2},\nu_{\frac{3+k-k^2}{2(k+2)}}}
  \end{equation}
  for \(k\ge 4\). For \(k=2\)\correct{, we have $\operatorname{wt}[u]=0$ and cannot invoke Theorem \ref{thm:generaldimformula}. However}, by \cref{prop:RmoduleResult},
   we can  assert
 \(\mathcal{V}^u\left(\rho_{k-\correct{2}}\right)\subset
 \eta^{3\frac{k^2-2k-4}{2(k+2)}}\mathcal{H}(\rho_{k-\correct{2}},\nu_{\correct{\frac{3+k-k^2}{2(k+2)}}})\),
 however, for all even
 \(k\ge 2\) we still have that
 \(\eta^{-3\frac{k^2-2k-4}{2(k+2)}}\vvptf(u,\tau)\) is a cyclic generator for
 \(\mathcal{H}(\rho_{k-2},\nu_{\frac{3+k-k^2}{2(k+2)}})\), since it has the right
 weight and spans a \(1\)-dimensional weight space.
  By construction, the leading exponents of the cyclic generator of
  \(\mathcal{H}(\rho_{k-2},\nu_{\frac{3+k-k^2}{2(k+2)}})\) are
  \(\set{0,\frac{k+1}{4(k+2)},\frac{1}{2}}\) (they form a minimal admissible
  set) and the weight of this cyclic generator is \(\frac{k+1}{k+2}\). This is
  implies the components of this \vvmf{} form a set of fundamental
  solutions for a third order monic modular linear differential equation of
  the form
  \begin{equation}
   \left(\partial_{\mdwt}^3 + \mdwtt \eis{4} \partial_{\mdwt} + \mdwttt \eis{6} \right)f =0,
   \label{eqn:3rdOrderMMDE02}
 \end{equation}
 see \cite[\refEq 15]{FrancMason-Hypergeometric} and the surrounding text.  While \cite{FrancMason-Hypergeometric}
   works in the context of integral weight vector-valued modular forms, a
   careful analysis of their construction of solutions \cite[\refEq 16]{FrancMason-Hypergeometric} to
   \eqref{eqn:3rdOrderMMDE02} shows that it is valid for real weight as
   well. In terms of the leading exponents
   \(\set{\lambda_1,\lambda_2,\lambda_3}\) a set of fundamental solutions is
   given by
   \begin{align}
     f_1 &=\eta^{2\mdwt} J^{-\frac{4\lambda_1 -2\lambda_2 -2\lambda_3 +1}{6}}
     {}_3F_2\left(\frac{4\lambda_1 -2\lambda_2 -2\lambda_3 +1}{6},
     \frac{4\lambda_1 -2\lambda_2 -2\lambda_3 +3}{6}, \frac{4\lambda_1
     -2\lambda_2 -2\lambda_3 +5}{6}; \lambda_1 -\lambda_2 +1, \lambda_1
     -\lambda_3+1; J^{-1} \right), \nonumber\\
     f_2 &=\eta^{2\mdwt} J^{-\frac{4\lambda_2 -2\lambda_1 -2\lambda_3 +1}{6}}
           {}_3F_2\left(\frac{4\lambda_2 -2\lambda_1 -2\lambda_3 +1}{6},
           \frac{4\lambda_2 -2\lambda_1 -2\lambda_3 +3}{6}, \frac{4\lambda_2
           -2\lambda_1 -2\lambda_3 +5}{6}; \lambda_2 -\lambda_1 +1, \lambda_2
           -\lambda_3+1; J^{-1} \right), \nonumber\\
     f_3 &=\eta^{2\mdwt} J^{-\frac{4\lambda_3 -2\lambda_1 -2\lambda_2 +1}{6}} {}_3F_2\left(\frac{4\lambda_3 -2\lambda_1 -2\lambda_2 +1}{6}, \frac{4\lambda_3 -2\lambda_1 -2\lambda_2 +3}{6}, \frac{4\lambda_3 -2\lambda_1 -2\lambda_2 +5}{6}; \lambda_3 -\lambda_1 +1, \lambda_3 -\lambda_2+1; J^{-1} \right),
   \end{align}
   with \(\mdwt=4(\lambda_1+\lambda_2+\lambda_3)-2\).
   Specialising to
   \(\set{\lambda_1,\lambda_2,\lambda_3}=\{0,\frac{k+1}{4(k+2)},\frac{1}{2}\}\)
   gives
   \begin{align}
     f_1&=\eta^{\frac{2k}{k+2}}J^{\frac{(k+1)}{12(k+2)}} {}_3F_2\left(-\frac{(k+1)}{12(k+2)}, \frac{11k+14}{24(k+2)}, \frac{19k+30}{24(k+2)}; \frac{3k+7}{4(k+2)}, \frac{1}{2}; J^{-1} \right),\nonumber\\
         f_2&=\eta^{\frac{2k}{k+2}}J^{-\frac{k+1}{6(k+2)}} {}_3F_2\left(\frac{k+1}{6(k+2)}, \frac{3k+5}{6(k+2)}, \frac{5k+9}{6(k+2)}; \frac{5k+9}{4(k+2)}, \frac{5}{8}; J^{-1} \right),\nonumber\\
         f_3&=\eta^{\frac{2k}{k+2}}J^{-\frac{5k+11}{12(k+2)}} {}_3F_2\left(\frac{5k+11}{12(k+2)}, \frac{9k+19}{12(k+2)}, \frac{13k+27}{12(k+2)}; \frac{3}{2}, \frac{5k+11}{4(k+2)}; J^{-1} \right).
   \end{align}
   The components of the cyclic generator are therefore linear combinations of
   the above fundamental solutions. More specifically, since the leading
   exponents of the respective components are
   \(\{0,\frac{k+1}{4(k+2)},\frac{1}{2}\}\), the cyclic generator must be
   of the form \((\alpha f_1,\beta f_2,\gamma f_3)^t\),
   \(\alpha,\beta,\gamma\in \CC^\times\). Hence
   \begin{align}
     \vvptf\brac*{u,\tau}=\eta^{3\frac{k^2-2k-4}{2(k+2)}}
     \begin{pmatrix}
       \alpha f_1\\
       \beta f_2\\
       \gamma f_3
     \end{pmatrix}
   \end{align}
   and the underlying intertwining operators can be normalised such that
   \(\alpha=\beta=\gamma=1\), so \eqref{eq:3dimptfformula} and
   \eqref{eq:3dimRmod} follow. That the \lhs{} of \eqref{eq:3dimbasis} is an \(\mathcal{M}\)-basis follows
   from \cref{thm:generaldimformula}. So all that remains is to relate powers
   of the modular derivative \correct{to} the action of Virasoro generators.   
   By the $a=\tilde{\omega}$ case of \eqref{eq:1pt[-1]result} we
 have $\vvptf_{k-2} (L_{[-2]}u,\tau)=\partial \vvptf_{k-2}
 (u,\tau)$. \correct{Furthermore, using the $a=\tilde{\omega}$ case of both \eqref{eq:1pt[0]result} and \eqref{eq:1pt[-1]result},
 the commutator relations among Virasoro modes, and Part \ref{itm:2dimirreducibility} of \cref{thm:sl2alldims}, we find
 that the components \(\ptf^{\mu}\), \(\frac{k-2}{2}\leq \mu \leq \frac{k+2}{2}\) of \(\vvptf(w,\tau)\),
 \(w\in L(k,k-2)\) satisfy 
 \begin{align}
 \ptf^{\mu} (\tilde{\omega}_{[-3]}u,\tau) &= \frac{1}{2} \ptf^{\mu} (L[-2](L[-1]^2 u),\tau) 
 = \frac{1}{2} \sum_{m=2}^\infty \eis{2m} \ptf^{\mu} (L[2m-2](L[-1]^2 u),\tau )
 \nonumber
 \\ &=\frac{1}{2} \eis{4} \ptf^{\mu} (L[2](L[-1]^2u),\tau) =3\operatorname{wt}[u]\eis{4}(\tau)\ptf^{\mu} (u,\tau).
 \label{eqn:Partial2Part1}
 \end{align}
 }
 We also have
 \correct{
 \begin{align}
  \ptf^{\mu} \left(\tilde{\omega}^2_{[-1]}u,\tau \right) &= 			
                                                           \mathrm{tr}\big\rvert_{L(k,\mu)} \left( L_0 - \frac{\mathbf{c}}{24}\right) o(\tilde{\omega}_{[-1]}u) q^{L_0-\frac{\mathbf{c}}{24}} + \sum_{m=1}^\infty  \eis{2m}(\tau) \ptf^{\mu} (L_{[2m-2]} L_{[-2]}u,\tau)\nonumber\\
		&= \frac{1}{2\pi i } \frac{\mathrm d}{\mathrm d\tau} 
			\mathrm{tr}\big\rvert_{L(k,\mu)} o(\tilde{\omega}_{[-1]}u) q^{L_0-\frac{\mathbf{c}}{24}} + \left(\wt [u]+2\right)\eis{2}(\tau)\ptf^{\mu} (\tilde{\omega}_{[-1]}u,\tau)+  \sum_{m=2}^\infty \eis{2m}(\tau) \ptf^{\mu} (L_{[2m-2]} L_{[-2]}u,\tau)\nonumber\\
		&= \partial^2 \ptf^{\mu} (u,\tau) + \sum_{m=2}^\infty \eis{2m}(\tau) \ptf^{\mu} (L_{[2m-2]} L_{[-2]}u,\tau) = \partial^2 \ptf^{\mu} (u,\tau) + \eis{4}(\tau) \ptf^{\mu} (L_{[2]} L_{[-2]}u,\tau)  \nonumber\\
		&=\partial^2 \ptf^{\mu} (u,\tau)+ \left( 4\mathrm{wt}[u] + \frac{\mathbf{c}}{2}\right) \eis{4}(\tau)\ptf^{\mu}  (u,\tau),
 \label{eqn:Partial2Part2}
 \end{align}
 where the first and third equalities follow from \eqref{eq:1pt[-1]result}.}
 Combining \eqref{eqn:Partial2Part1} and \eqref{eqn:Partial2Part2} we
 find $\partial^2 \vvptf_{k-2} (u,\tau) = \vvptf_{k-2} ((L^2_{[-2]} +
 \delta L_{[-4]})u,\tau)$ with $\delta = -(4\wt
 [u]+\mathbf{c}/24)/(3\wt [u])$. Plugging in the formula above for $\wt
 [u]$ and $\mathbf{c}$ gives the stated formula for $\delta$.
\end{proof}

See \cref{tab:3dtable} for explicit expansions of \(\vvptf_{k-2}(u,\tau)\) for
the first few values of the level \(k\).

\begin{table}[ht]
  \renewcommand{\arraystretch}{1.5}
  \begin{tabular}{|c||l|l|l|}
    \hline
    Level \(k\)
    & \(\vvptf_{k-2}(u,\tau)\)\\ \hline\hline
    \(4\)
    & \(\begin{matrix*}[l]
            q^{1/24} \left( 1 - \frac{6991}{171}q - \frac{1462930981}{198531}q^2 -\frac{11520966474250}{5360337}q^3 - \frac{467661528323716250}{627159429}q^4 +\cdots\right)\\
            q^{1/4}\left(1 + \frac{134}{9}q + \frac{167509}{81}q^2 + \frac{24672291010}{45927}q^3 +\frac{2054193740460070}{11986947}q^4 \cdots\right)\\
            q^{13/24} \left(1 -\frac{31}{27}q + \frac{473}{1215}q^2 -\frac{27056}{32805}q^3-\frac{1533931}{2657205}q^4 + \cdots \right)
          \end{matrix*}\)\\ \hline
      \(6\)
      &\(\begin{matrix*}[l]
           q^{5/32} \left(1-\frac{1041}{20}q-\frac{28822341}{3040}q^2-\frac{34699584029}{12160}q^3-\frac{2170275413391777}{2140160}q^4+ \cdots\right)\\
           q^{3/8} \left( 1 +\frac{74}{5}q +\frac{317943}{130}q^2+\frac{8423595}{13} q^3 + \frac{ 21692516271}{104}q^4 + \cdots\right)\\
           q^{21/32} \left(1 -\frac{31}{8}q + \frac{423}{128}q^2 +\frac{14247}{7168}q^3-\frac{485683}{229376}q^4 + \cdots \right)
         \end{matrix*}\)\\ \hline
      \(8\)
      & \(\begin{matrix*}[l]
            q^{11/40} \left(1-\frac{46803}{775}q-\frac{14944931541}{1375625}q^2-\frac{574656427747084}{171953125}q^3 -\frac{782261040133149248781}{649123046875}q^4 + \cdots\right)\\
            q^{1/2} \left(1 +\frac{1704}{125}q +\frac{108483138}{40625}q^2 + \frac{5094872662288}{7109375}q^3 +\frac{5893213005533601}{25390625}q^4+ \cdots\right)\\
            q^{31/40} \left( 1 - \frac{503}{75}q +\frac{ 44149}{3125}q^2-\frac{206842}{78125}q^3-\frac{420276376}{17578125}q^4 + \cdots \right)
          \end{matrix*}\)\\ \hline
      \(10\)
      & \(\begin{matrix*}[l]
            q^{19/48} \left(1- \frac{44717}{666}q -\frac{14421863479}{1222776}q^2-\frac{243672512766437}{66029904}q^3-\frac{68038738170466662661}{50617747584} q^4 + \cdots\right)\\
            q^{5/8} \left( 1  + \frac{107}{9}q + \frac{2963152}{1053}q^2 +\frac{64959522367}{85293}q^3 + \frac{5516615806491181}{22261473}q^4 + \cdots\right)\\
            q^{43/48} \left(1 -\frac{259}{27}q  +\frac{8110}{243}q^2-\frac{251140}{6561}q^3-\frac{25036652}{531441}q^4 + \cdots \right)
          \end{matrix*}\)\\ \hline
    \end{tabular}
    \caption{
      The first five terms of the \(q\)-series expansions for
\(\vvptf_{k-2}(u,\tau)\) for all levels $k$ at which \(\vvptf_{k-2}(u,\tau)\) generates \(\mathcal{H}(\rho_{k-2},\nu_{h_{k-2}})\). 
  In each case
  the series have been normalised so that
the leading coefficient is 1. This can always be achieved by an appropriate
choice of normalisation of the intertwining operators underlying
\(\vvptf_{k-2}(u,\tau)\).}
    \label{tab:3dtable}
\end{table}

\subsection{Select higher dimensions}

 We conclude this section by providing some results concerning non-congruence
 representations in some higher dimensions.

\begin{thm}
  Let the level be $k=p^t-2$, where $p>3$ is prime and $t$ is a positive integer.
  For $2\le \lambda \le k$ with $\lambda$ even,  if the representation
  \(\rho_\lambda\) is irreducible, then it is non-congruence if $t=1$ or if
  $t>1$ and $\lambda+1 > p^{t-2}$.\\
  \label{thm:gendimcongruence}
\end{thm}
\begin{proof}
  By \cref{thm:sl2alldims}.\ref{itm:basis},
  \begin{equation}
    \rho_\lambda(\modT) = \diag\set*{\epi(r_0), \dots, \epi(r_{k-\lambda})}
  \end{equation}
  where 
  \begin{equation}
    r_j = \frac{6j^2+6j(\lambda+2) + \lambda(\lambda+5)-3k}{24(k+2)}.
  \end{equation}
  Evaluating at $k = p^t-2$, the \(r_j\) become
  \begin{equation}
    r_j = \frac{6j^2 + 6j(\lambda+2) + \lambda(\lambda+5) -3p^t+6}{2^3 \cdot 3 p^t}.
    \label{eq:specialexp}
  \end{equation}
  Since $p>3$, the numerator in
  \eqref{eq:specialexp} is odd and hence indivisible by 2. Furthermore, the
  numerator is divisible by $3$ if and only if 3 divides $\lambda(\lambda+5)$,
  which is the case if and only if $\lambda \equiv 0,1 \pmod{3}$. The
  order of the $\rho_{\lambda}(\modT)$ is given by the least common multiple
  of the reduced denominators. To ensure this includes the factor of $p^t$, it
  suffices that $p$ does not divide the numerator for all $j = 0, \dots,
  k-\lambda = p^t-2-\lambda.$ As $j$ increases, the $j$th numerator is
  incremented by $18+12j+6\lambda$ which must be divisible by $p$ for all
  numerators to be divisible by $p$. This is only the case if $p$ divides
  $12$, i.e., $p=2,3$ which we have excluded. Thus the order of
  $\rho_{\lambda}(\modT)$ is $N = 2^3 \cdot 3 \cdot p^t$ if $\lambda
  \equiv 2 \pmod{3}$ and $ N = 2^3 \cdot p^t$ otherwise.
  
  To ascertain non-congruence based on the level and dimension of a
  representation, we follow the argument in \cite[\refSec 3.4]{intvvmf18}. Namely,
  for a $d$-dimensional congruence representation $\rho \colon \SLTZ \to
  \GL{d}$ of level $N$, its image is isomorphic to a quotient of 
  $\grp{SL}(2,\mathbb{Z}_{N})$. If $N = \prod_i p_i^{s_i}$ is the
  factorisation of \(N\) into distinct primes \(p_i\), then
  $\SLTZp{N} \cong \prod_i \grp{SL}(2,\mathbb{Z}_{p_i^{s_i}})$  
  and hence any irreducible
  representation of \(\grp{SL}(2,\mathbb{Z}_{N})\) can be constructed by tensoring irreducible
  representations of the \(\grp{SL}(2,\mathbb{Z}_{p_i^{s_i}})\) factors.
  In \cite{Nobs19761,Nobs19762}, all irreducible representations of
  $\grp{SL}(2,\mathbb{Z}_{p^t})$, for $p$ prime and \(t\in \NN\), were
  classified and their 
  dimensions were determined. In particular, tables summarising the
  classification are given in \cite[\refSec 9]{Nobs19762} (a summary of the
  minimal dimensions of non-trivial representations in English is given in
  \cite[\refThm 3.14]{intvvmf18}). Specifically the minimal dimensions of a
  \correct{representation} of level \(2^3\) or \( 3\) are \(2\) and \(1\),
  respectively. While for representations of level \(p^t\) the minimal
  dimension is \(\frac{1}{2}\brac*{p-1}\), if \(t=1\) and
  \(\frac{1}{2}(p^t-t^{t-2})\), if \(t\ge2\). Thus, the minimal dimension among
    representations of level \(N\) may 
    be found as a product of the minimal dimensions of the representations of
    level \(p_i^{s_i}\). Note that requiring each tensor factor to
    have respective level \(p_i^{s_i}\) precludes any of the tensor factors
    from being trivial. Thus for \(t=1\), we get that the minimal dimension is \(p-1\) for both
    \(N=2^3 \cdot 3\cdot p^t\) and \(N=2^3\cdot p^t\). Comparing this to
    \cref{thm:sl2alldims}.\ref{itm:basis} at level \(k=p^t-2\) we see that the
    dimension formula becomes \(p-1-\lambda\) which is less than \(p-1\) if
    \(\lambda\ge2\). Similarly, if \(t\ge2\) the minimal congruence dimension
    is \(p^t-p^{t-2}\). Hence we have non-congruence, if
    \(p^t-p^{t-2}>p^t-1-\lambda\), or equivalently, \(\lambda+1> p^{t-2}\).
\end{proof}

\begin{prop}
  Let $0\leq \lambda \leq k$, \(\lambda\) even, $u\in L(k,\lambda)$
  \correct{torus primary of the form} 
 \eqref{eq:actingvec}, and let $\rho_\lambda$ be the representation associated
 to $\vvptf_\lambda (u,\tau)$. 
 The dimension of the vector-valued modular form \(\vvptf_\lambda
 (u,\tau)\) is
 \(4\) if and only if $\lambda =k-3$ and hence the level \(k\) is
 odd. If
 \(\mathcal{H}(\rho_{k-3},\nu_{\frac{-k^2+4k+6}{2\brac*{k+2}}})\) is
 cyclic, then
 \begin{equation}
   \mathcal{V}^u(\rho_\lambda)=\mathcal{V}(\rho_\lambda)=\eta^{\frac{3}{2}\frac{k^2-4k-3}{k+2}}\mathcal{H}\brac*{\rho_{k-3},\nu_{\frac{-k^2+4k+6}{2\brac*{k+2}}}}.
 \end{equation}
\end{prop}
\begin{proof}
  Theorem \ref{thm:sl2alldims}.\ref{itm:basis} gives that $\rho_{\lambda}$ is
  a $4$-dimensional representation if and only if $\lambda =k-3$ is even. Next
  note that as \(h_{\frac{k+3}{2}}-h_{\frac{k-3}{2}}=\frac{3}{4}\) all conditions
  in \cref{thm:generaldimformula} other than condition \ref{itm:cyclicass}
  (cyclicity over \(\mathcal{R}\)) obviously hold. Thus the proposition
  follows for those levels where
  \(\mathcal{H}(\rho_{k-3},\nu_{\frac{-k^2+4k+6}{2\brac*{k+2}}})\) is cyclic.
\end{proof}

Note that in the case of general non-negative integral levels \(k\) and even
weight \(0\le \lambda\le k\), where one obtains \vvmf{s} of dimension
\(d=k-\lambda +1\), we have $\mu_{\max{}}-\mu_{\min{}}=(d-1)/4$. Thus, the
\correct{fifth} condition of \cref{thm:generaldimformula} holds only for those levels
and weights chosen so that $d\leq 4$ or equivalently $k-\lambda \leq 3$. For
\(d\ge5\) we therefore have that the inclusion
\begin{equation}
  \eta^{-\frac{3}{2}\frac{(k-d+1)(k-d+5)-2k}{k+2}}\mathcal{V}^u(\rho_\lambda)\subset
  \mathcal{H}\brac*{\rho_{k-d+1},\nu_{\frac{(k-d+1)(k-d+3)}{16\brac*{k+2}}}}
\end{equation}
is proper.

\section{Modular actions from categorical data}
\label{sec:mtcs}

So far we have studied the properties of traces of intertwining operators
directly, that is, using results from analytic number theory on modular
forms. However, since categories of modules over rational \voa{s} are modular
tensor categories, and additionally the categorical and number theoretic
notions of modularity coincide \correct{\cite{Huang-Rigidity,HuaVer08}}, we can repeat the above
analysis using categorical data. Let \(\catC\) be a \mtc{}, that is, a
monoidal category with many additional structures and properties (linear, abelian, semi-simple, finite, rigid, and braided with a non-degenerate braiding, etc.; see \cite{BKtensor2001,EGNO2015} for details).
To compute the action of the modular group, we will need the graphical calculus (also known as string diagram calculus, see \cite[\refSec 2.3]{BKtensor2001} for an introduction). To convert this abstract action of the modular group into actual matrices we will need to make explicit choices of bases (see \cite[\refSec 2]{FRSRCFT2002} for an introduction to working in such bases and some helpful identities), just as we needed to choose bases of intertwining operators in \cref{sec:intmodprops} above to obtain vector-valued modular forms.
Let \(\mathcal{I}\) be a complete set of representatives of simple isomorphism classes of objects in \(\catC\), with \(0\in \mathcal{I}\) denoting the tensor unit (that is, the \voa{} itself, if \(\catC\) is a category of \voa{} modules). The rigid dual of a simple object \(i\in\mathcal{I}\) will be denoted \(i^\ast\). For every triple \(i,j,k\in \mathcal{I}\), consider the vector space \(\Homgrp{\catC}{i\otimes j}{k}\) (called a 3-point coupling space) and pick a basis \(\{\coup{\alpha}{i}{j}{k}\}_{\alpha=0}^{\dim (\Homgrp{\catC}{i\otimes j}{k})-1}\) and denote its dual basis by \(\{\dcoup{\alpha}{i}{j}{k}\}_{\alpha}\subset \Homgrp{\catC}{k}{i\otimes j}\), where the evaluation of dual vector on vectors is given by
\begin{equation}
  \coup{\alpha}{i}{j}{k}\circ \dcoup{\beta}{i}{j}{k} =\delta_{\alpha,\beta}
  \id_k\in \Homgrp{\catC}{k}{k}=\CC \id_k,
  \label{eq:dualpairing}
\end{equation}
where $\delta_{\alpha ,\beta}$ is the Kronecker \(\delta\).
The 3-point coupling space \(\Homgrp{\catC}{i\otimes j}{k}\) is the
categorical counterpart to the space of intertwining operators of type
\(\binom{k}{i\ j}\) and picking a basis of 3-point couplings is equivalent to
picking a basis of intertwining operators. Therefore, a natural categorical
question is to ask how to characterise the subcategory of objects which
correspond to intertwining operators that map an object to itself (and hence
admit a trace as in \eqref{eq:q-expansionOfTraceFunction}). 
Recall the adjoint category \(\addcat{C}\) is defined to be the smallest full
subcategory of \(\catC\) containing all objects \(i\otimes i^\ast,\ i\in
\catC\) and all of their subquotients. Note that this category is closed under
taking duals. Another characterisation of \(\addcat{C}\) is as the centraliser
of the subcategory of invertible objects \cite[\refSec 4.14]{EGNO2015}.
\begin{lem}
  Let \(\catC\) be a modular tensor category. A simple object \(p\in\catC\) admits a non-vanishing \(\Homgrp{\catC}{p\otimes i}{i}\) for some \(i\in \catC\) if and only if \(p\) is in the adjoint subcategory \(\addcat{C}\).
\end{lem}
\begin{proof}
  Recall that Hom spaces in a modular tensor category satisfy natural isomorphisms
  \begin{equation}
    \Homgrp{\catC}{p\otimes i}{i}\cong \Homgrp{\catC}{i^\ast \otimes i}{p^\ast}.
  \end{equation}
  Thus the Hom spaces on the \rhs{} of the above identification are non-vanishing if and only \(p^\ast\) lies in the adjoint subcategory, which is the case if and only if \(p\) does.
\end{proof}

We continue fixing conventions. Note that if \(i=0\) or \(j=0\) then \(\dim
(\Homgrp{\catC}{0\otimes j}{k})=\delta_{j,k}\) and \(\dim (
\Homgrp{\catC}{i\otimes 0}{k})=\delta_{i,k}\). The non-vanishing 3-point
coupling spaces \(\Homgrp{\catC}{0\otimes j}{j}\) and
\(\Homgrp{\catC}{j\otimes 0}{\correct{j}}\) are spanned by the left and right
unitors, respectively, and so we choose these as our basis elements, that is,
\correct{\(\coup{0}{0}{j}{j}=\ell_j\)} and \correct{\(\coup{0}{i}{0}{i}=r_i\)}.
In our conventions for the graphical calculus we will always read diagrams from bottom to top (also called the optimistic direction).
The 3-point couplings and their duals are thus displayed as
\begin{equation}
  \coup{\alpha}{i}{j}{k}=
  \begin{tikzpicture}[baseline=(centre)]
    \coordinate (centre) at (0,0);
    \coordinate (top) at (0,1.5);
    \coordinate (bleft) at (-1.5,-1.5);
    \coordinate (bright) at (1.5,-1.5);
    \draw[black] (centre) -- (bleft);
    \draw[black] (centre) -- (bright);
    \draw[black] (centre) -- (top);
    \filldraw[black] (centre) circle (2pt) node [anchor=west] {$\alpha$};
    \node[anchor=south] at (top) {$k$};
    \node[anchor=east] at (bleft) {$i$};
    \node[anchor=west] at (bright) {$j$};
\end{tikzpicture},
\qquad  \dcoup{\alpha}{i}{j}{k}=
\begin{tikzpicture}[baseline=(centre)]
  \coordinate (centre) at (0,0);
  \coordinate (bot) at (0,-1.5);
  \coordinate (tleft) at (-1.5,1.5);
  \coordinate (tright) at (1.5,1.5);
  \draw[black] (centre) -- (tleft);
  \draw[black] (centre) -- (tright);
  \draw[black] (centre) -- (bot);
  \filldraw[black] (centre) circle (2pt) node[anchor=west] {$\alpha$};
  \node[anchor=north] at (bot) {$k$};
  \node[anchor=east] at (tleft) {$i$};
  \node[anchor=west] at (tright) {$j$};
\end{tikzpicture}
\end{equation}
and duality property \eqref{eq:dualpairing} is then expressed as
\begin{equation}
  \begin{tikzpicture}[baseline=(centre)]
	\node (centre) at (0,0) {$=$};
	\node at (-3,0) {$i$};
	\node at (-1,0) {$j$};
	\node at (1, 0) {$\delta_{\alpha, \beta}$};
	\filldraw[black] (-2,2) circle (0pt) node[anchor=south]{$k$};
	\filldraw[black] (-2,-2) circle (0pt) node[anchor=north]{$k$};
	\filldraw[black] (-2, 1) circle (2pt) node[anchor=north]{$\beta$};
	\filldraw[black] (-2, -1) circle (2pt) node[anchor=south]{$\alpha$};
	\filldraw[black] (2,-2) circle (0pt) node[anchor=north]{$k$};
	\draw[black] (-2,2) -- (-2, 1);
	\draw[black] (-2, -1) -- (-2, -2);
	\draw[black] (-2, 1) .. controls (-3,0.5) and (-3,-0.5) .. (-2,-1);
	\draw[black] (-2, 1) .. controls (-1,0.5) and (-1,-0.5) .. (-2,-1);
	\draw[black] (2, 2) -- (2,-2);
\end{tikzpicture}.
\end{equation}
With these choices of bases of 3-point couplings the associator structure morphisms (and their inverses) of \(\catC\)  can be expressed as the matrices
\begin{equation}
  \begin{tikzpicture}[baseline=(p)]
    \coordinate (b) at (3/2,3/4);
    \coordinate (i) at (0,0);
    \coordinate (j) at (1,0);
    \coordinate (k) at (2,0);
    \coordinate (a) at (1,3/2);
    \coordinate (l) at (1,2);
    \draw[-] (i) node[below=1mm] {$i$} -- (a) -- (l) node[above=1mm] {$l$};
    \draw[-] (j) node[below=1mm] {$j$} -- (b)  -- node[midway,below left=-1mm] (p) {$p$} (a);
    \draw[-] (k) node[below=1mm] {$k$}-- (b);
    \draw[fill] (a) circle (2pt) node[right=1mm] {$\alpha$};
    \draw[fill] (b) circle (2pt) node[right=1mm] {$\beta$};
  \end{tikzpicture}=
  \begin{tikzpicture}[baseline=(q)]
    \coordinate (g) at (-3/2,3/4);
    \coordinate (i) at (0,0);
    \coordinate (j) at (-1,0);
    \coordinate (k) at (-2,0);
    \coordinate (d) at (-1,3/2);
    \coordinate (l) at (-1,2);
    \draw[-] (i) node[below=1mm] {$k$} -- (d) -- (l) node[above=1mm] {$l$};
    \draw[-] (j) node[below=1mm] {$j$} -- (g)  -- node[midway,below right=-1mm] (q) {$q$} (d);
    \draw[-] (k) node[below=1mm] {$i$}-- (g);
    \draw[fill] (d) circle (2pt) node[left=1mm] {$\delta$};
    \draw[fill] (g) circle (2pt) node[left=1mm] {$\gamma$};
    \node[left =8mm of q] {$\displaystyle \sum_q\sum_{\gamma,\delta} F^{(i\,j\,k)l}_{\alpha p\beta,\, \gamma q \delta}$};
  \end{tikzpicture},\qquad
  \begin{tikzpicture}[baseline=(q)]
    \coordinate (g) at (-3/2,3/4);
    \coordinate (i) at (0,0);
    \coordinate (j) at (-1,0);
    \coordinate (k) at (-2,0);
    \coordinate (d) at (-1,3/2);
    \coordinate (l) at (-1,2);
    \draw[-] (i) node[below=1mm] {$k$} -- (d) -- (l) node[above=1mm] {$l$};
    \draw[-] (j) node[below=1mm] {$j$} -- (g)  -- node[midway,below right=-1mm] (q) {$p$} (d);
    \draw[-] (k) node[below=1mm] {$i$}-- (g);
    \draw[fill] (d) circle (2pt) node[left=1mm] {$\alpha$};
    \draw[fill] (g) circle (2pt) node[left=1mm] {$\beta$};    
  \end{tikzpicture}
  =
  \begin{tikzpicture}[baseline=(p)]
    \coordinate (b) at (3/2,3/4);
    \coordinate (i) at (0,0);
    \coordinate (j) at (1,0);
    \coordinate (k) at (2,0);
    \coordinate (a) at (1,3/2);
    \coordinate (l) at (1,2);
    \draw[-] (i) node[below=1mm] {$i$} -- (a) -- (l) node[above=1mm] {$l$};
    \draw[-] (j) node[below=1mm] {$j$} -- (b)  -- node[midway,below left=-1mm] (p) {$q$} (a);
    \draw[-] (k) node[below=1mm] {$k$}-- (b);
    \draw[fill] (a) circle (2pt) node[right=1mm] {$\gamma$};
    \draw[fill] (b) circle (2pt) node[right=1mm] {$\delta$};
    \node[left =8mm of p] {$\displaystyle \sum_p\sum_{\alpha,\beta} G^{(i\,j\,k)l}_{\alpha p\beta,\, \gamma q \delta}$};
  \end{tikzpicture},
  \label{eq:FGsymbols}
\end{equation}
while the braiding isomorphisms are expressed as the matrix
\begin{equation}
  \begin{tikzpicture}[baseline=(eq)]
    \coordinate (a) at (0,1); 
    \coordinate (b) at (1.8,2);
    \coordinate (c) at (1,3);
    \coordinate (d) at (0.2,2);
    \coordinate (e) at (2,1); 
    \coordinate (f) at (1,3.8);
    \coordinate (g) at (4,1); 
    \coordinate (h) at (5,3);
    \coordinate (i) at (6,1); 
    \coordinate (j) at (5,3.8);
    \coordinate (eq) at (3,2);
    \begin{knot}[consider self intersections]
      \strand (a) to [out=90, in=240] (b) to [out=60,in=0] (c) to [out=180,in=120] (d) to [out=300, in=90] (e);
    \end{knot}
    \node[below=1mm] at (a) {$i$};
    \node[below=1mm] at (e) {$j$};
    \draw[fill] (c) circle (2pt) node[above left=1mm] {$\alpha$};
    \draw[-] (c) -- (f) node[above=1mm] {$k$};
    \node at (eq) {$\displaystyle=\sum_{\beta} \Rmat^{(i\, j)k}_{\alpha\, \beta}\cdot$};
    \begin{knot}
      \strand (g) to [out=90, in=180] (h) to [out=0, in=90] (i);
    \end{knot}
    \node[below=1mm] at (g) {$i$};
    \node[below=1mm] at (i) {$j$};
    \draw[fill] (h) circle (2pt) node[above left=1mm] {$\beta$};
    \draw[-] (h) -- (j) node[above=1mm] {$k$};
  \end{tikzpicture}.
\end{equation}

\begin{thm}[Bakalov-Kirillov {\cite[\refThm 3.1.17]{BKtensor2001}}]
  Let \(p\in \mathcal{I}\) and consider the vector spaces \(W_{p,i}=\Homgrp{\catC}{p}{i\otimes i^\ast},\ i\in\mathcal{I}\), and their direct sum \(W_p=\bigoplus_{i\in\mathcal{I}} W_{p,i}\).
  Define linear maps \(\modS^{p},\modT^{p}\colon W_{p}\to W_{p}\) via the diagrams
  \begin{align}
      \begin{tikzpicture}[baseline=(eq)]
        \coordinate (eq) at (-2,0.5);
	\draw (0,0) -- (0,-1);
	\filldraw[black] (0,0) circle (2pt) node[anchor=south]{$\alpha$};
	\node at (0, -1.2) {$p$};
	\node at (-1.2,1) {$i$};
	\node at (-1.2, 2.2) {$j$};
	\begin{knot}[flip crossing = 2]
		\strand[
		only when rendering/.style={
			postaction=decorate,
		},
		decoration={
			markings,
			mark=at position 0.25 with {\arrowreversed{To}}
		}] (0,1) circle (1);
		\strand[only when rendering/.style={
			postaction=decorate,
		},
		decoration={
			markings,
			mark=at position 0.5 with {\arrow{To}}
		}] (1,2.2) arc[start angle = 0, end angle = -180, radius = 1];
	\end{knot}
	\node at (eq) {$\displaystyle\mapsto\quad
          \sum_{j\in\mathcal{I}}\frac{d_j}{D} \, \cdot\ $};
	\node at (-6.2, 0.5) {${\modS}^{(p)}\colon$};
	\draw[black] (-4.5,0.5) -- (-5.8,2) node[above=1mm]{$i$};
	\draw[black] (-4.5,0.5) -- (-3.3,2) node[above=1mm]{$i^*$};
	\draw[black] (-4.5,0.5) -- (-4.5,-1) node[below=1mm]{$p$};
	\filldraw[black] (-4.5,0.5) circle (2pt) node[anchor=west]{$\alpha$};
      \end{tikzpicture},\nonumber\\
    \begin{tikzpicture}
      [baseline=(eq)]
      \coordinate (eq) at (-2,.5);
	\draw[-] (0,.5) -- (-1.3,2) node[above=1mm]{$i$};
	\draw[-] (0,.5) -- (1.2,2) node[above=1mm]{$i^*$};
	\draw[-] (0,.5) -- (0,-1)  node[below=1mm]{$p$};
	\filldraw[black] (0,.5) circle (2pt) node[anchor=west]{$\alpha$}; 
	\node at (eq) {$\displaystyle\mapsto\qquad \frac{\theta_i}{\zeta}\cdot$}; 
	\node at (-6.2, 0.5) {${\modT}^{(p)}\colon$};
	\draw[black] (-4.5,0.5) -- (-5.8,2) node[above=1mm]{$i$};
	\draw[black] (-4.5,0.5) -- (-3.3,2) node[above=1mm]{$i^*$};
	\draw[black] (-4.5,0.5) -- (-4.5,-1) node[below=1mm]{$p$};
	\filldraw[black] (-4.5,0.5) circle (2pt) node[anchor=west]{$\alpha$};
    \end{tikzpicture},
  \end{align}
  where \(d_i\) is the quantum dimension of \(i\in\mathcal{I}\), \(D=\sum_i
  d_i^2\), and \(\zeta=\brac*{\frac{\sum_i\theta_i d_i}{\sum_i \theta_{i}^{-1}d_i}}^{\frac16}\). Then \({\modS}^{(p)}, {\modT}^{(p)}\) satisfy
  the relations \(\brac*{{\modS}^{(p)} {\modT}^{(p)}}^3=
  \brac*{{\modS}^{(p)}}^2\) and
  \(\brac*{{\modS}^{(p)}}^4=\theta_p^{-1}\). That is,
  \({\modS}^{(p)}, {\modT}^{(p)}\) satisfy the 
  defining relations 
  of the braid group \(\mathsf{B}_3\) on three strands (the modular group of
  the torus with one marked point) with the additional relation \(\brac*{{\modS}^{(p)}}^4=\theta_p^{-1}\) being the Dehn twist about the marked point.
  \label{thm:BKmodaction}
\end{thm}
The above theorem is a specialisation of \cite[\refThm 3.1.17 and
5.5.1]{BKtensor2001}, where Theorem 3.1.17 gives the action of the modular
group on duals of 3-point coupling spaces and Theorem 5.5.1 gives the action on
marked tori (which are the geometric interpretation of traces of intertwining
operators). Note that in order to be closer to the conventions of \voa{}
literature, we have rescaled the definition of \(\modT^{(p)}\) by a factor of
\(\zeta\) relative to the conventions of \cite{BKtensor2001}.
Note further that the action given in \cref{thm:BKmodaction} and above is an
action of \(\grp{B}_3\). To deprojectify and obtain an action of \(\SLTZ\) one
needs to include a multiplier system, which we shall do
a posteriori in the \(\alg{sl}(2)\) example below.
If \(\catC\) is a category of modules over a rational \voa{} (one satisfying
all of the assumptions in the paragraph preceding \eqref{eq:iopspace}), then the numbers
appearing in the theorem above can be expressed in terms of \voa{} data as
\(D=\modS_{0,0}^{-1}\), \(d_i=\frac{\modS_{i,0}}{\modS_{0,0}}\), where
\(\modS_{i,j}\) is the modular \(\modS\)-matrix of characters, and
\(\theta_p=\epi(h_p)\), where \(h_p\) is the conformal weight of the
simple module \(p\) and \(\zeta=\epi(\frac{\mathbf{c}}{24})\), where
\(\mathbf{c}\) is the central charge of the \voa{.}
\begin{thm}
  Let \(\catC\) be a \mtc{} with a
  set of representatives of simple isomorphism classes \(\mathcal{I}\), twist, braiding and fusing matrices given in a choice of basis of 3-point couplings, as described above.
  Let \(p\in\mathcal{I}\). 
  Then the pull back of \(\modS^{(p)},\modT^{(p)}\) to spaces of 3-point couplings via evaluation and co-evaluation, expanded in the basis \(\set{\coup{i}{p}{i}{\alpha}}\) is given by
  \begin{align}
    \Smat{i}{\alpha}{j}{\beta}{p}&=D^{-1} \cdot
    \begin{tikzpicture}[baseline=0pt]
        \begin{knot}[flip crossing = 2]
          \strand[only when rendering/.style={postaction=decorate,},decoration={markings,mark= between positions 0.25 and 0.75 step 0.5 with {\arrowreversed{To}}}]
          (0,0)  circle (1); 
          \strand[only when rendering/.style={postaction=decorate,},decoration={markings,mark=between positions 0.25 and 0.75 step 0.5 with {\arrowreversed{To}}}]
           (1.2,-.5) circle (1); 
         \end{knot}
         \draw[-] ($(1.2,-.5)+ (240:1)$) to [out=160, in=260] node[midway, below left] {$p$} (200:1);
         \filldraw[black] (200:1) circle (2pt) node[anchor=east]{$\alpha$};
         \filldraw[black] ($(1.2,-.5)+ (240:1)$) circle (2pt) node[anchor=north]{$\overline{\beta}$};
         \node at (135:1.2) {$i$};
         \node at ($ (1.2,-.5)+(135:1.2)$) {$j$};
      \end{tikzpicture} \nonumber\\
     &=
       \frac{d_i d_j}{D} \sum_{r\in i^\ast \otimes j}\sum_{\gamma\correct{,}\delta,\epsilon} 
       \frac{\twist{r}}{\twist{i}\twist{j}} \Gmat^{(ii^\ast j)j}_{0\,;\,\delta r \gamma}\Fmat^{(ii^\ast j)j}_{\epsilon r \gamma\,;\,0}\Gmat^{(pir)j}_{\delta i\alpha\,;\,\beta j \epsilon},
       \label{eq:smatformula}
  \end{align}
  where \(\gamma\) runs over the basis of
  \(\Homgrp{\catC}{i^\ast\otimes j}{r}\), and \(\epsilon,\delta\) both run over the
  basis of \(\Homgrp{\catC}{i \otimes r}{j}\). The modular \(\modT\)-matrix is given by
  \begin{equation}
    \Tmat{i}{\alpha}{j}{\beta}{p}=\delta_{i,j}\delta_{\alpha,\beta}\frac{\theta_i}{\zeta}.
  \end{equation}
  \label{thm:categoricalformulae}
\end{thm}
The above \(\modS\) and \(\modT\) matrices are the categorical counterpart
  to the analytic number theoretic ones discussed previously. Note however, that
  the multiplier system \(\nu_{h_p}\) has not yet been included in these formulae.
A similar diagrammatic formula for \(\Smat{i}{\alpha}{j}{\beta}{p}\) already
appeared \correct{in \cite[Ex. 7.1.f]{MSRCFT}}, but with some additional assumptions on
dimensions of 3-point coupling spaces, which we do away with here.
\begin{proof}
  The formula for \(\modT^{(p)}\) follows immediately from
  \cref{thm:BKmodaction}, so we focus on the formula for
  \(\modS^{(p)}\). The transferal of \(\modS^{(p)}\) as it is given in \cref{thm:BKmodaction} to 3-point couplings via evaluation and co-evaluation is
\begin{equation}
\begin{tikzpicture}[baseline=(basepoint)]
  \coordinate (basepoint) at (-4,0);
  \node (s) at (-5.5, 0) {$\modS^{(p)}\colon$};
	\node at (-2.6,0) {$\displaystyle\mapsto \sum_{j\in\mathcal{I}}\frac{d_j}{D}\, \cdot $};
	\draw[black] (-4,0) -- (-4,1.5);
	\draw[black] (-4,0) -- (-5.5,-1.5);
	\draw[black] (-4,0) -- (-2.5,-1.5);
	\filldraw[black] (-4,0) circle (2pt) node[anchor=west] (centre) {$\alpha$};
	\filldraw[black] (-4,1.5) circle (0pt) node[anchor=south]{$i$};
	\filldraw[black] (-5.5,-1.5) circle (0pt) node[anchor=west]{$p$};
	\filldraw[black] (-2.5,-1.5) circle (0pt) node[anchor=west]{$i$};
	\filldraw[black] (-1.7,0) circle (2pt) node[anchor=east]{};
	\node at (-1.3, 1.7) {$j$}; 
	\node at (-1.7,0.8) {$i$}; 
	\draw (-1.7,0) -- (-1.8,-2);
	\node at (-2.1, -2) {$p$};
	\begin{knot}[end tolerance=.01pt, flip crossing = 2]
		\strand[only when rendering/.style={
			postaction=decorate,
		},
		decoration={
                  markings,
                  mark=between positions 0.25 and 0.75 step 0.5 with {\arrowreversed{To}}
		}] (-.7,0) circle (1);
		\strand[only when rendering/.style={
			postaction=decorate,
		},
		decoration={
			markings,
			mark=at position 0.5 with {\arrow{To}}
		}] (.3,1) arc[start angle = 0, end angle = -180, radius = 0.8];
		\strand (-1.3, 1.5) -- (-1.3, 1);
		\strand[only when rendering/.style={postaction=decorate,},decoration={markings,	mark=at position 0.5 with {\arrowreversed{To}}}]
                (.3,1) arc[start angle = 0, end angle = -180, radius = -0.4];
		\strand (1.1, 1) -- (1.1, -2);
	\end{knot}
	\node at (2,0) {$\displaystyle=\sum_{j\in \mathcal{I}}\frac{d_j}{D} \, \cdot $};
	\draw (3.2,0) -- (2.9, -1.5);
	\filldraw[black] (2.9, -1.5) circle (0pt) node[anchor=north]{$p$};
	\filldraw[black] (3.2,0) circle (2pt) node[anchor=east]{$\alpha$};
	\node at (5.2, -1) {$j$}; 
	\node at (3.4, 1) {$i$}; 
	\begin{knot}[flip crossing = 1]
		\strand[only when rendering/.style={
			postaction=decorate,
		},
		decoration={
			markings,
			mark=at position 0.5 with {\arrowreversed{To}}
		}] (5.2,0) arc[start angle = 0, end angle = 180, radius = 1];
		\strand[only when rendering/.style={
			postaction=decorate,
		},
		decoration={
			markings,
			mark=at position 0.5 with {\arrow{To}}
		}] (5.2,0) arc[start angle = 0, end angle = -180, radius = 1];
		\strand (5,1.5) -- (5, -1.5);
	\end{knot}
      \end{tikzpicture}
      \begin{tikzpicture}[baseline=(basepoint)]
        \coordinate (basepoint) at (-4,0);
        \node (s) at (-5.5, 0) {$\displaystyle=\sum_{j, \beta}\modS^{(p)}_{i\alpha\, j\beta}\cdot$};
	\draw[black] (-4,0) -- (-4,1.5);
	\draw[black] (-4,0) -- (-5.5,-1.5);
	\draw[black] (-4,0) -- (-2.5,-1.5);
	\filldraw[black] (-4,0) circle (2pt) node[anchor=west] (centre) {$\beta$};
	\filldraw[black] (-4,1.5) circle (0pt) node[anchor=south]{$j$};
	\filldraw[black] (-5.5,-1.5) circle (0pt) node[anchor=west]{$p$};
	\filldraw[black] (-2.5,-1.5) circle (0pt) node[anchor=west]{$j$};
      \end{tikzpicture},
    \end{equation}
    where the first identity uses the straightening axiom of evaluation and co-evaluation to yield a morphism in \(\Homgrp{\catC}{p\otimes j}{j}\) and the second identity is the expansion of this morphism in our chosen basis, which defines the coefficients \(\modS^{(p)}_{i\alpha\, j\beta}\).
To extract the coefficient in front of each basis vector we pair with the dual
basis by attaching \(\dcoup{\beta}{p}{j}{j}\) to the diagram from below, which
will yield a morphism in \(\Homgrp{\catC}{j}{j}=\CC \id_j\), proportional to
the identity, that is,
\begin{equation}
  \modS^{(p)}_{i\alpha\, j\beta}\cdot
  \begin{tikzpicture}[baseline = {(0,0)}]
    \draw (0,1) -- (0,-1) node[right] {\(j\)};
  \end{tikzpicture}
  =\frac{d_j}{D}\cdot
  \begin{tikzpicture}[baseline = {(0,0)}]
    \draw (0,0) to [out=250, in=135] node[midway,below left] {\(p\)} (1, -2)
     -- (1,-2.5) node[right]{\(j\)};
    \filldraw[black] (0,0) circle (2pt) node[left]{$\alpha$};
    \node at (0.2, 1) {$i$}; 
    \begin{knot}[flip crossing = 1]
      \strand[only when rendering/.style={
        postaction=decorate,
      },
      decoration={
        markings,
        mark=at position 0.5 with {\arrowreversed{To}}
      }] (2,0) arc[start angle = 0, end angle = 180, radius = 1];
      \strand[only when rendering/.style={
        postaction=decorate,
      },
      decoration={
        markings,
        mark=at position 0.5 with {\arrow{To}}
      }] (2,0) arc[start angle = 0, end angle = -180, radius = 1];
      \strand (1.8,1.5) node[below right] {\(j\)} -- (1.8, -0.7) to [out=270, in=45] (1,-2) ;
    \end{knot}
    \filldraw[black] (1,-2) circle (2pt) node[above]{$\beta$};   
  \end{tikzpicture}.
\end{equation}
We can then take the trace over \(j\) to get a morphism in
\(\Homgrp{\catC}{0}{0}=\CC\id_0\), that is, we connect the \(j\)-strands at the
top and bottom of the diagram using evaluation and co-evaluation. The \lhs{} of
the identity is then just a circle labelled by \(j\) and this evaluates to
\(d_j\). The \rhs{} then becomes string diagram in
\eqref{eq:smatformula}.
To evaluate the string diagram 
we need the well known identities
\begin{align}
  \begin{tikzpicture}[decoration={markings, mark=at position 0.5 with {\arrow{>}}}, baseline=(eq)]
	\draw[postaction={decorate}] (-1,.5) to [out=30, in=150] (1,.5);
	\draw[postaction={decorate}] (1,-.5) to [in=-30, out=-150] (-1,-.5);
        \fill[gray] (-1.2,.5) rectangle (1.2,-.5);
	\filldraw[black] (-1,.5) circle (0pt) node[anchor=south east]{$i$};
	\filldraw[black] (1,.5) circle (0pt) node[anchor=south west]{$i^*$};
	\filldraw[black] (-1,-.5) circle (0pt) node[anchor=north east]{$i$};
	\filldraw[black] (1,-.5) circle (0pt) node[anchor=north west]{$i^*$};
	\node (eq) at (1.9,0) {$=\; d_i \, \cdot$};
        \fill[gray] (2.6,.5) rectangle (5,-.5);
	\draw (2.8,.5) -- (3.8,1) -- (4.8,.5);
	\draw (2.8,-.5) -- (3.8,-1) -- (4.8,-.5);
	\draw[dashed] (3.8,1) -- (3.8,1.5);
	\draw[dashed] (3.8,-1) -- (3.8,-1.5);
	\filldraw[black] (2.8,.5) circle (0pt) node[anchor=south east]{$i$};
	\filldraw[black] (4.8,.5) circle (0pt) node[anchor=south west]{$i^*$};
	\filldraw[black] (2.8,-.5) circle (0pt) node[anchor=north east]{$i$};
	\filldraw[black] (4.8,-.5) circle (0pt) node[anchor=north west]{$i^*$};
	\filldraw[black] (3.8,1.5) circle (0pt) node[anchor=south]{$0$};
	\filldraw[black] (3.8,-1.5) circle (0pt) node[anchor=north]{$0$};
      \end{tikzpicture}
  \quad\text{and}\quad
      \begin{tikzpicture}[baseline=(eq)]
        \begin{knot}[flip crossing = 2]
          \strand (0,0) ..controls +(0,.5) and +(0,-.5) .. (.5,1) .. controls +(0,.5) and +(0,-.5) .. (0,2);
          \strand (.5,0) .. controls +(0,.5) and +(0,-.5) .. (0,1) .. controls +(0,.5) and +(0,-.5) .. (.5,2);
        \end{knot}
        \node at (0,-.3) {$i^\ast$};
        \node at (.5,-.3) {$j$};
        \node (eq) at (2,.8) {$\displaystyle=\ \sum_{r\in i^\ast\otimes j}\sum_{\gamma}\dfrac{\theta_r}{\theta_i\theta_j}\cdot$};
        \draw (3,0) -- (3.5,.5) -- (4,0);
        \node at (3,-.3) {$i^\ast$};
        \node at (4,-.3) {$j$};
        \draw (3.5,.5) -- (3.5,1.5);
        \draw (3,2) -- (3.5,1.5) -- (4,2);
        \node at (3,2.3) {$i^\ast$};
        \node at (4,2.3) {$j$};
        \filldraw[black] (3.5,.5) circle (2pt) node[anchor=north]{$\gamma$};
         \filldraw[black] (3.5,1.5) circle (2pt) node[anchor=south]{$\overline{\gamma}$};
      \end{tikzpicture},
\end{align}
where the grey boxes in the left identity can contain any diagram and where we
have also used that \(\theta_{i^\ast}=\theta_i\). Then
\begin{align}
  \begin{tikzpicture}[baseline={(0,0)}]
        \begin{knot}[flip crossing = 2]
          \strand[only when rendering/.style={postaction=decorate,},decoration={markings,mark= between positions 0.25 and 0.75 step 0.5 with {\arrowreversed{To}}}]
          (0,0)  circle (1); 
          \strand[only when rendering/.style={postaction=decorate,},decoration={markings,mark=between positions 0.25 and 0.75 step 0.5 with {\arrowreversed{To}}}]
           (1.2,-.5) circle (1); 
         \end{knot}
         \draw[-] ($(1.2,-.5)+ (240:1)$) to [out=160, in=260] node[midway, below left] {$p$} (200:1);
         \filldraw[black] (200:1) circle (2pt) node[anchor=east]{$\alpha$};
         \filldraw[black] ($(1.2,-.5)+ (240:1)$) circle (2pt) node[anchor=north]{$\overline{\beta}$};
         \node at (135:1.2) {$i$};
         \node at ($ (1.2,-.5)+(135:1.2)$) {$j$};
       \end{tikzpicture}&=\sum_{r\in i^\ast \otimes j}\sum_\gamma\frac{\theta_r}{\theta_i\theta_j}\cdot
       \begin{tikzpicture}[baseline={(0,0)}]
         \coordinate (a) at (0,0);
         \coordinate (b) at (1.5,-1.7);
         \coordinate (gd) at (2,0);
         \coordinate (gu) at (2,1);
         \coordinate (ti) at (1,2);
         \coordinate (tj) at (3,2);
         \coordinate (bi) at (1,-1);
         \coordinate (bj) at (3,-2.5);
         \filldraw[black]  (a) circle (2pt);
         \filldraw[black]  (b) circle (2pt);
         \filldraw[black]  (gd) circle (2pt);
         \filldraw[black]  (gu) circle (2pt);
         \draw[->] (bi) to [out=180,in=315] node[midway,right] {\(i\)} (a) to [out=120, in=180]
         node[midway,above] {\(i\)} (ti);
         \draw[->] (ti) to [out=0,in=135] (gu) to [out=45,in=180]
         node[midway,right] {\(j\)} (tj);
         \draw[->] (tj) to [out=0, in=0] (bj);
         \draw[->] (bj) to [out=180,in=300] (b) to [out=60, in=315]
         node[midway,right] {\(j\)} (gd) to
         [out=225, in=0] (bi);
         \draw[-] (b) to [out=160,in=260] node[midway,below left] {\(p\)} (a);
         \draw[-] (gd) -- node[midway,right] {\(r\)} (gu);
         \node[left=3pt of a]  {\(\alpha\)};
         \node[below left=2pt of b] {\(\overline{\beta}\)};
         \node[above=3pt of gu] {\(\gamma\)};
         \node[below=3pt of gd] {\(\gamma\)};
       \end{tikzpicture}\nonumber\\
       &=d_i\sum_{r\in i^\ast \otimes
         j}\sum_\gamma\frac{\theta_r}{\theta_i\theta_j}\cdot
       \begin{tikzpicture}[baseline={(0,0)}]
         \coordinate (a) at (0,0);
         \coordinate (b) at (1.3,-1.9);
         \coordinate (gd) at (2,0);
         \coordinate (gu) at (2,1);
         \coordinate (ti) at (1,2);
         \coordinate (tj) at (3,2);
         \coordinate (bi) at (1,-1);
         \coordinate (bj) at (3,-2.5);
         \filldraw[black]  (a) circle (2pt);
         \filldraw[black]  (b) circle (2pt);
         \filldraw[black]  (gd) circle (2pt);
         \filldraw[black]  (gu) circle (2pt);
         \draw[-] (bi) to [out=135,in=315] node[midway,right] {\(i\)} (a) to [out=120, in=225]
         node[midway,above] {\(i\)} (ti);
         \draw[->] (ti) to [out=315,in=135] node[midway,below] {\(i^\ast\)} (gu) to [out=45,in=180]
         node[midway,right] {\(j\)} (tj);
         \draw[-,dashed] (ti) -- (1,2.5);
         \draw[->] (tj) to [out=0, in=0] (bj);
         \draw[-] (bj) to [out=180,in=300] (b) to [out=60, in=315]
         node[midway,right] {\(j\)} (gd) to
         [out=225, in=45] node[midway,above] {\(i^\ast\)} (bi);
         \draw[-,dashed] (bi) -- (1,-1.5);
         \draw[-] (b) to [out=160,in=260] node[midway,below left] {\(p\)} (a);
         \draw[-] (gd) -- node[midway,right] {\(r\)} (gu);
         \node[left=3pt of a]  {\(\alpha\)};
         \node[below left=2pt of b] {\(\overline{\beta}\)};
         \node[above=3pt of gu] {\(\gamma\)};
         \node[below=3pt of gd] {\(\gamma\)};
       \end{tikzpicture}
       .
  \label{eq:diagrameval}
\end{align}

Applying the well known identities
\begin{equation}
  \begin{tikzpicture}[baseline={(2,-1)}]
    \draw[-] (0,0) node[above] {\(i\)}  --  (1,-2) -- node[midway,above]
    {\(i^\ast\)} (2,-1)  node[circle,fill=black, inner sep=0pt,minimum
    size=4pt] {}
    node [right] {\(\gamma\)}--   (2,0) node[above] {\(r\)}; 
    \draw[-,dashed] (1,-2) -- (1,-3) node[below] {\(0\)};
    \draw[-] (2,-1) -- (3,-3) node[below] {\(j\)};
  \end{tikzpicture}=\sum_\epsilon\Fmat^{(i\,i^\ast\,j)j}_{\epsilon r \gamma,0}\cdot
  \begin{tikzpicture}[baseline={(1,-1)}]
    \draw[-] (0,0) node[above] {\(i\)} -- (1,-1) node[circle,fill=black, inner sep=0pt,minimum
    size=4pt] {} node [right] {\(\epsilon\)} -- (2,0) node[above] {\(r\)};
    \draw[-] (1,-1) -- (1,-2) node[below] {\(j\)};
  \end{tikzpicture},\quad
  \begin{tikzpicture}[baseline={(1,2)}]
    \draw[-] (0,0) node[below] {\(i\)} -- (1,2) --  node[midway,above]
    {\(i^\ast\)} (2,1)  node[circle,fill=black, inner sep=0pt,minimum
    size=4pt] {} node [right] {\(\gamma\)} -- (3,3) node[above] {\(j\)};
    \draw[-,dashed] (1,2) -- (1,3) node[above] {\(0\)};
    \draw[-] (2,1) -- (2,0) node[below] {\(r\)};
  \end{tikzpicture}
  =\sum_\delta \Gmat^{(i\,i^\ast\,j)j}_{0,\delta r \gamma}\cdot
  \begin{tikzpicture}[baseline = {(1,1)}]
    \draw[-] (0,0) node[below] {\(i\)} -- (1,1) node[circle,fill=black, inner sep=0pt,minimum
    size=4pt] {} node [right] {\(\delta\)} -- (2,0) node[below] {\(r\)};
    \draw[-] (1,1) -- (1,2) node[above] {\(j\)};
  \end{tikzpicture}
\end{equation}
and
\begin{equation}
  \begin{tikzpicture}[baseline={(1,1.5)}]
    \draw[-] (0,0) node[below] {\(\ell\)} --
    (0,1) node[circle,fill=black, inner sep=0pt,minimum size=4pt] {}
    node[left] {\(\gamma\)} --
    node[midway,left] {\(i\)}
    (-1,2.5) node[circle,fill=black, inner sep=0pt,minimum size=4pt] {}
    node[left] {\(\beta\)} --
    node[midway,above left] {\(p\)}
    (0,3)  node[circle,fill=black, inner sep=0pt,minimum size=4pt] {} node[right] {\(\alpha\)} --
    (0,4) node[above] {\(\ell\)};
    \draw[-] (0,1) -- node[midway,below right] {\(q\)}
    (1,1.5) node[circle,fill=black, inner sep=0pt,minimum size=4pt] {}
    node[right] {\(\delta\)} --
    node[midway,below left] {\(j\)}
    (-1,2.5);
    \draw[-] (1,1.5) -- node[midway,above right] {\(k\)} (0,3);
  \end{tikzpicture}\quad =\quad
  \Gmat^{(i\,j\,k)\ell}_{\alpha p\beta,\gamma q \delta}\cdot
  \begin{tikzpicture}[baseline={(0,0)}]
    \draw[-] (0,-1) node[right] {\(\ell\)} -- (0,1);
  \end{tikzpicture}
\end{equation}
to the last diagram in \eqref{eq:diagrameval} and again using the fact that the
circle labelled by \(j\) evaluates to
\(d_j\) yields the formula in \eqref{eq:smatformula}.
\end{proof}

\subsection{The \texorpdfstring{\(\SLA{sl}{2}\)}{sl2} example}
Explicit formulae for the Moore-Seiberg data of the \mtc{} for affine
\(\SLA{sl}{2}\) at level \(k\in \mathbb{N}_0\) are known and we reproduce them
here to compute some examples. We take the label set of simple modules to be
their highest weights \(\mathcal{I}=\{0,\dots,k\}\) with respective conformal
weights \(h_n=\frac{n(n+2)}{4(k+2)}\), \(n\in\mathcal{I}\) and central charge
\(\mathbf{c}=\frac{3k}{k+2}\). Since for \(\SLA{sl}{2}\) the 3-point coupling
spaces are always at most \(1\)-dimensional, no labels are needed for basis
vectors (the Greek indices
above).
The character \(\modS\)-matrix entries are given by
\begin{equation}
  \modS_{i,j}=\sqrt{\frac{2}{k+2}}\sin\brac*{\frac{\pi(i+1)(j+1)}{k+2}},\quad
  D=\modS_{0,0}^{-1},\quad d_i=\frac{\modS_{i,0}}{\modS_{0,0}}.
\end{equation}
The twist, \(\zeta\), braiding and fusing matrices are, respectively, given by
\begin{equation}
  \theta_r = \epi \left(h_r\right),
  \quad \zeta= \epi \left(\frac{\mathbf{c}}{24}\right),\quad
  \Rmat^{(r\, s)t}= (-1)^{r+s-t} \epi\left(\frac{h_r+h_s-h_t}{2} \right),\quad
  \Fmat^{(r\, s\, t\,)u}_{p\, q}=
    \begin{Bmatrix}
      t/2&s/2&p/2\\
      r/2& u/2&q/2
    \end{Bmatrix},
\end{equation}
where
\begin{align}
  \begin{Bmatrix}
      a&b&e\\
      d&c&f
    \end{Bmatrix}
    &=
      (-1)^{a+b-c-d-2e}\sqrt{\sqbrac*{2e+1}\sqbrac{2f+1}}\Delta(a,b,e)\Delta(a,c,f)\Delta(c,e,d)\Delta(d,b,f)\nonumber\\
       &\quad
        \times\hspace{-15mm} \sum_{z=\min\{a+b+c+d,a+d+e+f,b+c+e+f\}}^{\max\{a+b+e,a+c+f,b+d+f,c+d+e\}}\hspace{-15mm}
         (-1)^z\sqbrac*{z+1}!\left(\sqbrac*{z-a-b-e}!\sqbrac*{z-a-c-f}!\sqbrac*{z-b-d-f}!\sqbrac*{z-d-c-e}!\right. \nonumber\\
       &\hspace{20mm}
         \left.\sqbrac*{a+b+c+d-z}!\sqbrac*{a+d+e+f-z}!\sqbrac*{b+c+e+f-z}!\right)^{-1}
\end{align}
are quantum group \(6j\)-symbols and
\begin{align}
  \Delta(a,b,c)&=\sqrt{\frac{\sqbrac*{-a+b+c}!\sqbrac*{a-b+c}!\sqbrac*{a+b-c}!}{\sqbrac*{a+b+c+1}!}},\nonumber\\
  \sqbrac*{n}&=\frac{\sin\brac*{\frac{\pi n}{k+2}}}{\sin\brac*{\frac{\pi}{k+2}}},\qquad
  \sqbrac*{n}!=\prod_{m=1}^n\sqbrac*{m},\qquad \sqbrac*{0}!=1.
\end{align}
The above formulae can be found in \cite{KirResh89,Hou6j90}. Note that the
\(\Gmat\) can be computed from the above data
using the identity
\begin{equation}
  \Gmat^{(i\, j\, k)\ell}_{p\, q}=\frac{\Rmat^{(j\, k)q}\Rmat^{(i\, q)\ell}}{\Rmat^{(i\, j)p}\Rmat^{(p\,
      k)\ell}}\Fmat^{(k\, j\, i)\ell}_{p\, q}.
\end{equation}

Consider now the example when the label \(p\) of the acting object is equal to \(k\). Then we are in the $1$-dimensional case with \(\modS\) and \(\modT\) given by
\begin{equation}
  \modS^{(k)}=\ee^{-\frac{\pi}{2} \ii\brac*{h_k+\frac{k}{2}}}=\ee^{-\frac{\pi}{2} \ii\frac{3k}{8}},\qquad \modT^{(k)}=\ee^{\frac{\pi}{4}\ii \frac{k}{2}}.
\end{equation}
In particular, \(\modS^{(k)}\) and \(\modT^{(k)}\) are equal to the evaluation
of the multiplier system \(\nu_{h_k+\frac{k}{2}}\), so the categorical data appear to detect that the natural vector \(u\) to use for torus 1-point functions has conformal weight \(h_k+\frac{k}{2}\), as we have seen in the sections above. If we divide the above formulae by the \(\nu_{h_k}\) multiplier system then we recover the formulae in \cref{thm:1dimcongruence}.

The above formulae for \(\modS\) and \(\modT\) can also be used to show that
\cref{thm:gendimcongruence} admits examples of non-congruent representations
of dimension greater than three.
\begin{prop}
  Choose \(p=7\), \(t=1\) in \cref{thm:gendimcongruence} and
  hence \(k=5\). Then for \(\lambda=2\), the representation \(\rho_2\) is
  \(4\)-dimensional, irreducible and non-congruence. 
\end{prop}
\begin{proof}
  All of the conditions of \cref{thm:gendimcongruence} except for
  irreducibility hold by construction. So we only need to show
  irreducibility. Note that the existence or absence of a non-trivial
  invariant subspace does not depend on whether a multiplier system is
  included in the formulae for \(\modS\) or \(\modT\). Therefore, we
  work
  directly with the formulae in \cref{thm:categoricalformulae}.
  Note in
  particular that all of the eigenvalues of
  \begin{equation}
    \modT^{(2)}=
    \diag\brac*{\epi\brac*{\frac{1}{56}},\epi\brac*{\frac{11}{56}},\epi\brac*{\frac{25}{56}},\epi\brac*{\frac{43}{56}}}
    \label{eq:4dTmat}
  \end{equation}
  are distinct. Therefore a non-trivial invariant subspace would need to admit
  a basis \(B\) that is a proper non-empty subset of
  \(\modT^{(2)}\)-eigenvectors. The columns in \(\modS^{(2)}\) corresponding
  to these basis vectors would hence need to contain entries that are \(0\) in
  those rows which correspond to \(\modT^{(2)}\)-eigenvectors not in
  \(B\). However, 
  \begin{equation}
    \modS^{(2)}\approx
    \begin{pmatrix}
      -0.16 - 0.33\ii& -0.26 - 0.55\ii& -0.26 - 0.55\ii& -0.16 - 0.33\ii\\
      -0.26 - 0.55\ii& -0.16 - 0.33\ii&  0.16 + 0.33\ii&  0.26 + 0.55\ii\\
      -0.26 - 0.55\ii&  0.16 + 0.33\ii&  0.16 + 0.33\ii& -0.26 - 0.55\ii\\
      -0.16 - 0.33\ii&  0.26 + 0.55\ii& -0.26 - 0.55\ii&  0.16 + 0.33\ii
    \end{pmatrix}
  \end{equation}
  has no entries that are \(0\) and hence the representation must be
  irreducible. Here we have chosen to give a numerical approximation of
  \(\modS^{(2)}\) to two significant digits for simplicity, as the exact
  expression in terms of radicals is impractically large to present.
\end{proof}
Note that \(\nu_{h_2}\brac*{\modT}=\epi (\frac{1}{42})\) and so if we
divide the diagonal entries in \eqref{eq:4dTmat} by \(\nu_{h_2}\) and take the
product of the first two diagonal entries, we obtain
\(\epi (\frac{1}{6})\), which is a \(12\)th root of unity. Hence
  \cref{thm:irredcond} does not apply and we were only able to conclude
  irreducibility because of the categorical formulae.


\providecommand{\bysame}{\leavevmode\hbox to3em{\hrulefill}\thinspace}
\providecommand{\MR}{\relax\ifhmode\unskip\space\fi MR }
\providecommand{\MRhref}[2]{%
  \href{http://www.ams.org/mathscinet-getitem?mr=#1}{#2}
}
\providecommand{\href}[2]{#2}

\end{document}